\documentclass[12pt]{amsart}
\usepackage{macros}

\usepackage{todonotes}

\newcommand{\sara}[1]{}

\title[The metric space of limit laws for $q$-hook formulas]{The metric space of limit laws for $q$-hook formulas}
\author{Sara C. Billey and Joshua P. Swanson}
\date{\today}

\begin{document}

\begin{abstract}
In \cite{BKS-asymptotics}, Billey--Konvalinka--Swanson studied the
asymptotic distribution of the coefficients of Stanley's $q$-hook
length formula, or equivalently the major index on standard tableaux
of straight shape and certain skew shapes.  We extend those
investigations to Stanley's $q$-hook-content formula related to
semistandard tableaux and $q$-hook length formulas of Bj\"orner--Wachs
related to linear extensions of labeled forests.  We show that, while
their coefficients are ``generically'' asymptotically normal, there
are uncountably many non-normal limit laws. More precisely, we
introduce and completely describe the compact closure of the metric
space of distributions of these statistics in several regimes. The
additional limit distributions involve generalized uniform sum
distributions which are topologically parameterized by certain
decreasing sequence spaces with bounded $2$-norm. The closure
of these distributions in the L\'evy metric gives
rise to the space of DUSTPAN distributions.  As an application,
we completely classify the limiting distributions of the size
statistic on plane partitions fitting in a box.
\end{abstract}

\maketitle

\setcounter{tocdepth}{2}
\tableofcontents

\section{Introduction}
\label{sec:intro}

The famed Frame--Robinson--Thrall \textit{hook length formula} is a
rational product formula for counting the number of \textit{standard
Young tableaux} of a given partition shape $\lambda$
\cite{Frame-Robinson-Thrall.1954}, denoted $\SYT(\lambda)$.  Stanley's
$q$-analogue of the hook length formula \cite[Cor.~7.21.5]{ec2} is a remarkably simple
generalization for the polynomial generating function of the
\textit{major index} statistic on $\SYT(\lambda)$.  His $q$-\textit{hook length formula}
replaces each integer $n$ with the corresponding $q$-integer $[n]_q
\coloneqq 1+q+\cdots+q^{n-1}$, times an overall shift of
$q^{r(\lambda)}$ where $r(\lambda) \coloneqq \sum_{i \geq 1} (i-1) \lambda_i$:
\begin{equation}\label{eq:q_SYT}
\sum_{T \in \SYT (\lambda )} q^{\maj (T)} =  q^{r(\lambda)} \frac{[n]_q!}{\prod_{u \in \lambda} [h_u]_q}.
\end{equation}
Consequently, \eqref{eq:q_SYT} encodes probabilistic information
concerning the distribution of the major index statistic when sampling
from $\SYT(\lambda)$ uniformly at random.

\cite{BKS-asymptotics} considered the
distribution of $\maj$ on $\SYT(\lambda)$.
Given a sequence of partitions, \cite{BKS-asymptotics} completely determined when
the corresponding sequence of standardized random variables converges in distribution.
Equivalently, this determines the asymptotic distribution of the coefficients of
Stanley's $q$-hook length formula. For these random variables, countably many
continuous limit laws are possible:
one gets the normal distribution ``generically'' and, in certain degenerate regimes,
the \textit{Irwin--Hall distributions}. A key technical tool in \cite{BKS-asymptotics}
is an exact formula for the \textit{cumulants} of the underlying random variables,
which follows easily from work of Chen--Wang--Wang \cite{Chen-Wang-Wang.2008} and
Hwang--Zacharovas \cite{MR3346464} together with Stanley's $q$-hook length formula
\eqref{eq:q_SYT}.

The present work generalizes the explorations of
\cite{BKS-asymptotics} to the next most famous $q$-analogues of the
hook length formula: Stanley's $q$-hook-content formula for
semistandard tableaux, and formulae of Bj\"orner--Wachs for linear
extensions of labeled forests. See \Cref{tab:formulas} for a summary
of the $q$-hook-type formulas we use.  The limit laws in these cases
turn out to be much more intricate than in \cite{BKS-asymptotics},
with uncountably many rather than countably many possible limits.

Typical central limit theorems are based on an integer sequence so
they  ``let $n \to \infty$,'' even when the limit laws are
complicated such as in the work of Chatterjee--Diaconis
\cite{Chatterjee.Diaconis.2014}. By contrast, the combinatorial
statistics considered here and in \cite{BKS-asymptotics} have much
more complex indexing sets involving objects like integer partitions
and forests. We address this complication by considering sets of
standardized distributions as metric spaces under the L\'evy metric on
all distributions, together with a corresponding space of
parameters. Our overarching goal is to describe the closure of these
metric spaces and to completely classify which sequences tend to which
limit points in terms of the relevant parameter spaces.

\begin{table}
  \begin{center}
  \begin{tabu}{p{2cm}|p{2cm}|p{5cm}|p{5cm}}
  Statistic(s) & set & $q$-hook formula(s) & cumulant expression(s) \\
  \toprule\toprule
      $\maj$
      & $\SYT(\lambda)$
      & $q^{r(\lambda)} \frac{[n]_q!}{\prod_{u \in \lambda} [h_u]_q}$
      & $\sum_{i=1}^n j^d - \sum_{u \in \lambda} h_u^d$ \\
      \midrule
      $\rank$
      & $\SSYT_{\leq m}(\lambda)$
      & $q^{r(\lambda)} \prod_{u \in \lambda} \frac{[m+c_u]_q}{[h_u]_q}$
      & $\sum_{u \in \lambda} (m+c_u)^d - h_u^d$ \\
      \ 
      & \ 
      & $q^{r(\lambda)} \prod_{1 \leq i < j \leq m}
          \frac{[\lambda_i - \lambda_j + j - i]_q}{[j-i]_q}$
      & ${\scriptstyle \sum_{1 \leq i < j \leq m} (\lambda_i - \lambda_j + j - i)^d - (j-i)^d}$ \\
      \midrule
      $\size$
      & ${\scriptstyle \PP(a \times b \times c)}$
      & ${\scriptstyle \prod_{i=1}^a \prod_{j=1}^b
            \prod_{k=1}^c \frac{[i+j+k-1]_q}{[i+j+k-2]_q}}$
      & ${\scriptstyle \sum_{i,j,k} (i+j+k-1)^d - (i+j+k-2)^d}$ \\
      \midrule
      $\maj$
      & $\cL(P, w)$
      & $q^{\maj(P, w)} \frac{[n]_q!}{\prod_{u \in P} [h_u]_q}$
      & $\sum_{i=1}^n j^d - \sum_{u \in \lambda} h_u^d$ \\
      $\inv$
      & \ 
      & $q^{\inv(P, w)} \frac{[n]_q!}{\prod_{u \in P} [h_u]_q}$
      & \  \\
          \bottomrule
  \end{tabu}
  \end{center}
  \caption{Summary of combinatorial objects, statistics, $q$-hook formulas, and cumulant expressions used in this paper. Cumulants are obtained from cumulant expressions by multiplying by $\frac{B_d}{d}$ for $d > 1$. See \Cref{sec:background} for details.}
  \label{tab:formulas}
\end{table}


A key step in our approach is the introduction of a new family of continuous univariate distributions which we call \textit{DUSTPAN distributions}\footnote{A ``\textit{\underline{d}istribution associated to a \underline{u}niform \underline{s}um for $\mathbf{\underline{t}}$ \underline{p}lus \underline{a}n independent \underline{n}ormal distribution}.''}. These distributions involve convolutions of the normal law with a countable family of uniform measures supported on some intervals. More precisely, we have the following abstract description. See \Cref{def:uniform+nornmal} for the concrete version.

\begin{Theorem}
    The family of DUSTPAN distributions with variance $1$ is uniquely characterized as the smallest family $\cF$ of standardized real-valued distributions such that:
\begin{enumerate}[(i)]
  \item $\cU[0, 1]^* \in \cF$
  \item If $X, Y \in \cF$, then the standardized independent sum random variable $\frac{\alpha X + \beta Y}{\sqrt{\alpha^2 + \beta^2}}$ belongs to $\cF$ for any $\alpha, \beta \in \bR$ not both $0$.
  \item $\cF$ is closed under convergence in distribution.
\end{enumerate}
\end{Theorem}

The general strategy of our arguments is as follows. First, we convert formulas involving ratios of $q$-integers into explicit expressions for the cumulants.   In most cases these expressions involve significant cancellation. Next comes the difficult step where we find an asymptotically cancellation-free approximation to the cumulants in a suitable regime; see for instance \Cref{lem:forest_dc}. Finally, in all cases considered in this paper, we use the approximate cumulants to identify the limiting standardized distributions pertaining to SSYT's and linear extensions of trees as some particular DUSTPAN distribution. While the first step is quite generic, the combinatorial arguments and inequalities underlying the second step are highly domain-specific and expand on the corresponding approach from \cite{BKS-asymptotics}.



In \Cref{ssec:intro:SYT}, we summarize the results of
\cite{BKS-asymptotics} and reframe them in terms of metric spaces as a
prelude to our new, more technical results on semi-standard tableaux
and forests.  To keep this introduction to a manageable length and
avoid frequent digressions, we assume familiarity with tableaux
combinatorics and cumulants.  Detailed background on these topics is
provided in \cite[\S2]{BKS-asymptotics} or \cite[Ch.7]{ec2}.  The main
new results in this paper are outlined in \Cref{ssec:intro:SSYT} and
\Cref{ssec:intro:forests}.  See \Cref{sec:background} for background
necessary for the new material.

\subsection{Standard tableaux}\label{ssec:intro:SYT}

Let $\cX_\lambda[\maj]$ denote the random variable associated with
$\maj$ on $\SYT(\lambda)$, sampled uniformly at random.  Then the
probability $\mathbb{P}(\cX_\lambda[\maj]=k) =
a_{k}^{\lambda}/f^{\lambda }$ where $ \SYT(\lambda)^{\maj}(q) = \sum
a_{k}^{\lambda }q^{k}$ and $f^{\lambda}= \SYT(\lambda)^{\maj}(1)$ is
the number of standard Young tableaux of shape $\lambda$.  Hence,
studying the distribution of the random variable $\cX_\lambda[\maj]$
and the sequence of coefficients $\{a_{k}^{\lambda}: k\geq 0 \}$ for $
\SYT(\lambda)^{\maj}(q)$ are essentially equivalent.  Furthermore, any
polynomial in $q$ with nonnegative integer coefficients can be
associated to a random variable in a similar way.

For the sake of understanding limiting distributions, we typically
standardize the random variables involved so they have mean 0 and
variance 1.  In general, given any random variable $\cX$ with mean
$\mu$ and standard deviation $\sigma>0$, let $\cX^* \coloneqq (\cX -
\mu)/\sigma$ denote the corresponding \textit{standardized random
variable} with mean $0$ and variance $1$.  To avoid overemphasizing
trivialities, we implicitly ignore degenerate distributions with
$\sigma=0$ throughout the paper without further comment, so every
distribution we consider does have a standardization.  Write $\cX_n
\Rightarrow \cX$ to mean that the sequence $\cX_n$ converges in
distribution to $\cX$. Let $\cN(\mu, \sigma^{2})$ denote a normal
distribution, and let $\cIH_M$ denote the $M$th \textit{Irwin--Hall
distribution}, obtained by summing $M$ independent continuous uniform
$[0, 1]$ random variables.  These distributions are also referred to
as \textit{uniform sum distributions} in the literature. Note that the
normal and Irwin--Hall distributions are continuous, while each of the
random variables coming from $q$-hook formulas below determine
discrete distributions.

We may completely describe the possible limit distributions of $\cX_\lambda[\maj]^*$
using a simple auxiliary statistic on partitions, $\aft$. In particular,
let $\aft(\lambda) \coloneqq |\lambda| - \max\{\lambda_1, \lambda_1'\}$.

\begin{Theorem}\cite[Thm.~1.7]{BKS-asymptotics}\label{thm:BKS}
  Let $\lambda^{(1)}, \lambda^{(2)}, \ldots$ be a sequence of partitions where
  $|\lambda^{(N)}| \to \infty$ as $N \to \infty$.
  \begin{enumerate}[(i)]
    \item $\cX_{\lambda^{(N)}}[\maj]^* \Rightarrow \cN(0, 1)$ if and only if
      $\aft(\lambda^{(N)}) \to \infty$.
    \item $\cX_{\lambda^{(N)}}[\maj]^* \Rightarrow \cIH_M^{*}$ if and only if
      $\aft(\lambda^{(N)}) \to M < \infty$.  
  \end{enumerate}
\end{Theorem}

\Cref{thm:BKS} shows that the set $\bZ_{\geq 1} \cup \{\infty\}$
parameterizes the set of all possible limit distributions
associated to the $q$-hook length formulas and the standardized 
random variables $\cX_\lambda[\maj]^*$.  If we instead parameterize the limit
distributions by $\left\{\frac{1}{n} : n \in \mathbb{Z}_{\geq
1}\right\} \cup \{0\},$ we get a parameter space and a distribution
space which are homeomorphic as topological spaces.  Hence, we
introduce the notion of a metric space of standardized distributions.


\begin{Definition} The   \textit{metric space of Irwin--Hall
distributions} is 
 \[ \bfM_{\cIH} \coloneqq \{\cIH_M^* : M \in \mathbb{Z}_{\geq 1}\}, 
 \]
and the \textit{metric space of $\SYT$ distributions} is 
\[ \bfM_{\SYT} \coloneqq \{\cX_\lambda[\maj]^* : \lambda \in
\Par, f^\lambda > 1\}. \]   Endow
$\bfM_{\cIH}$ and $\bfM_{\SYT}$ with the topology inherited from the
topology of distributions of real-valued random variables under the
L\'evy metric, which is characterized by convergence in distribution
\cite[Ex.~14.5]{MR1324786}.  
\end{Definition}

By the Central Limit Theorem, $\overline{\bfM_{\cIH}}=\bfM_{\cIH}\cup
\{\cN(0, 1)\}$.  In light of \Cref{thm:BKS}, we have the following
very precise description of the minimal compactification of the metric
space of $\SYT$ distributions.

\begin{Corollary}\label{cor:cM_SYT_structure}
  In the L\'evy metric, \begin{equation}\label{eq:cM_SYT_overline}
\overline{\bfM_{\SYT}} = \bfM_{\SYT} \sqcup \overline{\bfM_{\cIH}}, \end{equation}
which is compact. Moreover,
the set of limit points of $\bfM_{\SYT}$ is exactly $\overline{\bfM_{\cIH}}$.
\end{Corollary}



When the set parametrizing our combinatorial statistics has a natural topology,
one might hope that it is homeomorphic to the space of distributions. For example,
let \[ \bfP_{\cIH} \coloneqq \left\{\frac{1}{n} : n \in
\mathbb{Z}_{\geq 1}\right\} \] be the \textit{Irwin--Hall parameter
space}.  We endow $\bfP_{\cIH} \subset [0, 1]$ with the topology of
pointwise convergence, so $\overline{\bfP_{\cIH}}=\bfP_{\cIH} \sqcup
\{0 \}.$ Since $\cIH_M^* \Rightarrow \cN(0, 1)$ as $M \to \infty$, the
bijection $\overline{\bfP_{\cIH}} \to \overline{\bfM_{\cIH}}$ given by
$\frac{1}{M} \mapsto \cIH_M^*$ and $0 \mapsto \cN(0, 1)$ is a
homeomorphism.  It is less clear how to impose a topology on standard
Young tableaux, but a characterization of the multiset of hook lengths
would be a key consideration.  See \cite[Thm.~7.1]{BKS-asymptotics}.

\begin{Remark}
  Recent work of Kim--Lee identified certain normal \cite{MR3998822}
and bivariate normal \cite{1811.04578} distributions as limits of
normalizations of $\des$ and $(\des, \maj)$ over conjugacy classes in
the symmetric group. In their context, the space of limit
distributions is parameterized by real numbers in $[0, 1]$.
\end{Remark}



\subsection{Semistandard tableaux and plane partitions}\label{ssec:intro:SSYT}

Stanley's \textit{hook-content formula} is a rational product formula
for counting the set $\SSYT_{\leq m}(\lambda)$ of \textit{semistandard
tableaux} of shape $\lambda$ with entries at most $m$. He gave a
natural $q$-analogue of this formula, which is in fact the polynomial
generating function for the \textit{rank} statistic on $\SSYT_{\leq
m}(\lambda)$.  A second rational product formula for $\rank$ on $\SSYT_{\leq
m}(\lambda)$ with important representation-theoretic meaning is given
by the type $A$ case of the $q$-Weyl dimension formula.  Explicitly,

\begin{equation}\label{eq:ssyt.product.formulas}
\sum_{T \in \SSYT_{\leq m}(\lambda)} q^{\rank (T)} = q^{r(\lambda)}
\prod_{u \in \lambda} \frac{[m+c_u]_q}{[h_u]_q} = q^{r(\lambda)}
\prod_{1 \leq i < j \leq m} \frac{[\lambda_i - \lambda_j + j -
i]_q}{[j-i]_q}.
\end{equation}
See \Cref{ssec:back:ssyt} for more details.

Let $\cX_{\lambda; m}[\rank]$ denote the random variable associated
with the $\rank$ statistic on $\SSYT_{\leq m}(\lambda)$, sampled
uniformly at random.  In \Cref{sec:background}, we derive simple
explicit cumulant formulas from these rational expressions which allow
us to study the possible limiting distributions for $\cX_{\lambda;
m}[\rank]^*$.  While the closures $\bfM_{\SYT}$ and $\bfM_{\cIH}$ are
completely characterized above, the closure of the \textit{metric
space of $\SSYT$ distributions}, \[ \bfM_{\SSYT} \coloneqq
\{\cX_{\lambda; m}[\rank]^* : \lambda \in \Par, 
\ell(\lambda)  \leq m \}, \] is much more complicated. In particular, we show
that the following generalization of the Irwin--Hall distributions are
related to limit laws for $\cX_{\lambda; m}[\rank]^*$. 

\begin{Definition}\label{def:uniform-sums}
Given a finite multiset $\mathbf{t}$ of
non-negative real numbers, let 
\begin{equation}\label{eq:t_uniform_sum}
  \cS_{\mathbf{t}} \coloneqq \sum_{t \in \mathbf{t}} \cU\left[-\frac{t}{2}, \frac{t}{2}\right],
\end{equation}
where we assume the summands are independent and $\cU[a, b]$ denotes
the continuous uniform distribution supported on $[a, b]$. If
$\mathbf{t}$ consists of $M$ copies of $1$, then $\cS_{\mathbf{t}} +
\frac{M}{2} = \cIH_M$.  By convention, we consider the multiset
$\mathbf{t}$ as a weakly decreasing sequence of real numbers
$\mathbf{t} = \{t_1 \geq t_2 \geq \cdots \geq t_m\}$ where $t_m \geq 0$.
We call the distribution associated to $\cS_{\mathbf{t}}$ a
\textit{finite generalized uniform sum distribution}.
\end{Definition}

Certain sequences of random variables $\cX_{\lambda; m}[\rank]^*$
which converge to a finite generalized uniform sum distribution are
completely characterized by an auxiliary multiset called the distance
multiset.  This auxiliary set also comes up in the \textit{Turnpike
Reconstruction Problem}, which is essentially the problem of
identifying all possible sequences $\mathbf{t}$ from the
following multiset $\Delta\mathbf{t}$, which has applications in
DNA sequencing and X-ray
crystallography \cite[Sect. 10.5.1]{WeissMarkDsaa}.
The Turnpike Reconstruction Problem is a potential candidate
for being in \textbf{NP-Intermediate}. See \cite{MR2038493}
for further computational complexity considerations.

\begin{Definition}\label{def:distance_multiset}
The \textit{distance multiset} of $\mathbf{t}=\{t_1 \geq t_2 \geq
\cdots \geq t_m \}$ is the multiset \[ \Delta\mathbf{t} \coloneqq
\{t_i - t_j : 1 \leq i < j \leq m\}. \] 
\end{Definition}

To avoid highly cluttered notation coming from the terms in a sequence
indexed by a parameter $N=1,2,\ldots$, we will often drop the explicit
dependence on $N$.  For example, let $\lambda$ and $m$ denote a
sequence of partitions $\lambda^{(1)}, \lambda^{(2)}, \ldots$ and a
sequence of values $m^{(1)}, m^{(2)}, \ldots $ respectively.  If we
assume $\ell(\lambda^{(N)})<m^{(N)}$ for each $N$, we will simply
write $\ell(\lambda) < m$.  Also, $|\lambda|=n$ means there is another
sequence $n^{(1)}, n^{(2)}, \ldots$ such that
the size of the partition $|\lambda^{(N)}|=n^{(N)}$, thus $|\lambda| \to \infty$ and $n \to
\infty$ both imply $|\lambda^{(N)}|\to \infty$ as $N \to \infty$.
Similarly, let $\cX_{\lambda; m}[\rank]$ denote the sequence of
uniform random variables associated with $\SSYT_{\leq
m^{(N)}}(\lambda^{(N)})^{\rank}(q)$. 

\begin{Theorem}\label{thm:SSYT_Delta}
  Let $\lambda$ be an infinite sequence of partitions with
$\ell(\lambda) < m$ where $\lambda_1/m^3 \to \infty$. Let
$\mathbf{t}(\lambda )=(t_{1},\ldots,t_m) \in [0,1]^{m}$ be the finite
multiset with $t_k \coloneqq \frac{\lambda_k}{\lambda_1}$ for $1 \leq
k \leq m$.  Then $\cX_{\lambda; m}[\rank]^*$ converges in distribution
if and only if the multisets $\Delta\mathbf{t}(\lambda )$ converge
pointwise.  In that case, the limit distribution is $\cN(0,1)$ if $m
\to \infty $ and $\cS_{\mathbf{d}}^*$ where $\Delta\mathbf{t}(\lambda
) \to \mathbf{d}$ if $m$ is bounded.
\end{Theorem}

\Cref{thm:SSYT_Delta} suggests we consider the \textit{metric space of
distance distributions} 
\begin{equation}\label{eq:dist.multiset.modulispc}
\bfM_{\DIST} \coloneqq \bigcup_{m \geq 2}\{\cS_{\Delta\mathbf{t}}^* : \mathbf{t} = \{1=t_1 \geq \cdots \geq t_m=0\}\}
\end{equation}
and its associated parameter space $\bfP_{\DIST}$ defined
in \Cref{sub:dist}.
By padding with $0$'s, we consider $\bfP_{\DIST} \subset
\mathbb{R}^{\bN}$ as a sequence space with the topology of pointwise
convergence.  The metric space of distance distributions is
significantly more complex than the metric space of Irwin--Hall
distributions. Nonetheless, a careful analysis involving the topology
of the parameter space of distance multisets done in \Cref{sub:dist}
yields the following results.  We will show that both $\bfP_{\DIST}$
and $\bfM_{\DIST}$ have natural one point compactifications,
\[
\overline{\bfP_{\DIST}} = \bfP_{\DIST} \sqcup
\{\mathbf{0}\} \text{ and } \overline{\bfM_{\DIST} } = \bfM_{\DIST}
\sqcup \{\cN(0, 1)\},
\]
where $\mathbf{0}$ is the infinite sequence of $0$'s.  Furthermore, in
analogy with \Cref{cor:cM_SYT_structure}, we will show that the map
$\overline{\bfP_{\DIST}} \to \overline{\bfM_{\DIST}}$ given by $\mathbf{d} \mapsto
\cS_{\mathbf{d}}^*$ and $\mathbf{0} \mapsto \cN(0, 1)$ is a
homeomorphism between sequentially compact spaces. See
\Cref{thm:PDelta_MDelta}.  Therefore, \Cref{thm:SSYT_Delta} and
\Cref{thm:PDelta_MDelta} combine to give the following complete
characterization of the possible limit laws for a particular family of
semistandard tableaux in analogy with \Cref{cor:cM_SYT_structure}.

\begin{Corollary}\label{cor:cM_SSYT_structure}
  For any fixed $\epsilon>0$, 
let
\[ \bfM_{\epsilon\SSYT}
\coloneqq \{\cX_{\lambda; m}[\rank]^* :
\ell(\lambda) < m \text{ and } \lambda_1/m^3 > (|\lambda| +
m)^\epsilon\} \subset \bfM_{\SSYT}.
\]
Then
\begin{equation}\label{eq:cM_SSYT_overline}
\overline{\bfM_{\epsilon\SSYT}} = \bfM_{\epsilon\SSYT} \sqcup \overline{\bfM_{\DIST}},
\end{equation}
which is compact. Moreover, the set of limit points of
$\bfM_{\epsilon\SSYT}$ is $\overline{\bfM_{\DIST}}$.
\end{Corollary}

\Cref{cor:cM_SSYT_structure} already indicates that the limiting
distributions associated to semistandard tableaux are much more varied
than the case of standard Young tableaux.  See \Cref{sum:SSYT} for a
synopsis of all of the asymptotic limits we have identified for
$\cX_{\lambda, m}[\rank]^*$.  This includes several ``generic''
asymptotic normality criteria and a partial analogue of $\aft$, called
$\weft$, which controls asymptotic normality in many cases of
interest. A complete description of the closure of $\bfM_{\SSYT}$ akin
to \Cref{thm:BKS} and \Cref{cor:cM_SYT_structure} remains open.

\begin{Problem}
  Describe $\overline{\bfM_{\SSYT}}$ in the L\'evy metric. What are all possible
  limit points? 
\end{Problem}

By studying one more special family of semistandard tableaux, we will
show that the Irwin--Hall distributions are also among the limit
points.  Thus, the strongest statement we have shown for the metric
space of limit laws for Stanley’s $q$-hook-content formula is
\[
\bfM_{\SSYT} \cup \bfM_{\DIST} \cup \bfM_{\cIH} \cup \{\cN(0, 1)\}  \subset \overline{\bfM_{\SSYT}}.
\]

Using a well-known bijection, the two product formulas in
\eqref{eq:ssyt.product.formulas} imply product formulas for the
generating function of the \textit{size} statistic on the set $\PP(a
\times b \times c)$ of \textit{plane partitions} fitting in a box.
See the second and third rows of \Cref{tab:formulas}.  Let $\cX_{a
\times b \times c}[\size]$ similarly denote the random variable
associated with the $\size$ statistic on $\PP(a \times b \times c)$.
In the theorem below, we give a complete characterization of the limit
laws for plane partitions and $\{\cX_{a \times b \times
c}[\size]^*\}$.  This leads to an analog of
\Cref{cor:cM_SYT_structure} for the \textit{metric space of plane
partition distributions}, denoted $\bfM_{\PP} \coloneqq \{\cX_{a
\times b \times c}[\size]^*\}$.



\begin{Theorem}\label{thm:PP_median.intro}
Let $a,b,c$ each be a sequence of positive integers.  
   \begin{enumerate}[(i)]
     \item $\cX_{a \times b \times c}[\size]^{*} \Rightarrow \cN(0, 1)$ if and only if
       $  \median\{a, b, c\}  \to \infty$.
     \item $\cX_{a \times b \times c}[\size]^{*} \Rightarrow \cIH_M$ if
       $  ab  \to M < \infty$ and $c \to \infty$.
   \end{enumerate}
\end{Theorem}

\begin{Corollary}\label{cor:cM_PP_structure}
In the L\'evy metric, \begin{equation}\label{eq:cM_PP_overline}
\overline{\bfM_{\PP}} = \bfM_{\PP} \sqcup \overline{\bfM_{\cIH}}, \end{equation}
which is compact. Moreover, the set of limit points of
$\bfM_{\PP}$ is exactly $\overline{\bfM_{\cIH}}$.
\end{Corollary}

\subsection{Linear extensions of forests}\label{ssec:intro:forests}

Knuth \cite[p.~70]{Knuth} gave a rational product formula for counting
the set $\cL(P)$ of linear extensions of a forest $P$, analogous to
the Frame--Robinson--Thrall hook length formula. Using a fixed
bijection $w \colon P \to [n]$, one may interpret $\cL(P)$ as a set of
permutations $\cL(P, w) \subset S_n$ and consider the distribution of
the major index or inversion number statistics on these
permutations. Stanley \cite{MR0332509} and Bj\"orner--Wachs
\cite{MR1022316} gave $q$-analogues of Knuth's formula for major index and
number of inversions using certain labelings $w$. All of these statistics agree up to an
overall shift. See the fourth row of \Cref{tab:formulas} and
\Cref{ssec:forests.background} for details.

Let $\cX_P$ denote the random variable associated with the $\maj$ or
$\inv$ statistic on $\cL(P, w)$ where $w$ is order-preserving. The
distribution of $\cX_P^{*}$ is independent of the choice of statistic
and the choice of $w$. Let \[ \bfM_{\Forest} \coloneqq \{\cX_P^* : \text{$P$ is a
forest}\} \] be the \textit{metric space of forest distributions}.  We
show that the behavior of the possible limiting distributions for
$\cX_P^{*}$ breaks into two distinct regimes.  The first ``generic''
regime exhibits classic asymptotic normality, while the second
``degenerate'' regime allows even more continuous limit laws than have
appeared in the theory for standard or semistandard tableaux.

Let $\rank(P)$ denote the length of a maximal chain in $P$. Let $|P|$
denote the number of vertices. For example, the rank of a complete
binary tree with $2^n-1$ vertices is $n$, so $\rank(P) \approx \log_2
|P|$.  Typically, $\rank(P)$ is much smaller than $|P|$, so the
following theorem covers the ``generic'' regime.

\begin{Theorem}\label{thm:generic.forest.sum}
  Given a sequence of forests $P$, the corresponding
  sequence of random variables $\cX_P^*$ is asymptotically
  normal if
    \[ |P| \to \infty \qquad \text{and} \qquad \limsup \frac{\rank(P)}{|P|} < 1. \]
\end{Theorem}

In the ``degenerate'' regime, $\rank(P) \sim |P|$, so the number of
vertices not in a chosen maximal  chain is relatively small. We completely
describe the possible limit distributions when $|P| - \rank(P) =
o(|P|^{1/2})$.  To do so, we generalize both the distance
distributions and the Irwin--Hall distributions to the distributions
associated to countable sums of independent, continuous, uniform
random variables with finite mean and variance. We call these
\textit{generalized uniform sum distributions}.  Again we can reduce
to sums of independent centralized random variables $\cS_\mathbf{t}$
exactly as in \eqref{eq:t_uniform_sum}, except now we consider
countably infinite multisets $\mathbf{t}=\{t_{1}\geq
t_{2}\geq \dots \}$ of nonnegative real numbers.  See \Cref{sub:inf.sums} for details such as
cumulants, the density function, and the relation to pointwise
convergence in $\mathbb{R}^{\mathbb{N}}$.

The variance of a uniform sum random variable $\cS_\mathbf{t}$ is
closely related to the \textit{$2$-norm} of $\mathbf{t}$,
\[
|\mathbf{t}|_2 \coloneqq \left(\sum_{t \in
\mathbf{t}} t^2\right)^{1/2}.
\]
In this notation, $\Var[\cS_{\mathbf{t}}]= \frac{B_2}{2} |\mathbf{t}|_2^2,$ where $B_{2}=\frac{1}{6}$ is a Bernoulli number.  Thus, in order for $\cS_\mathbf{t}$ to be well
defined, it must have finite variance, so $|\mathbf{t}|_2<\infty$ is required.
Let $\widetilde{\ell}_2 \coloneqq \{\mathbf{t} = (t_1, t_2, \ldots) :
t_1 \geq t_2 \geq \cdots \geq 0, |\mathbf{t}|_2 < \infty\}$. The
standardized general uniform sum distributions are indexed by the
decreasing sequences $\mathbf{t} \in\widetilde{\ell}_2$ such that
$1=\Var[\cS_{\mathbf{t}}]= \frac{B_2}{2}
|\mathbf{t}|_2^2, $ so $|\mathbf{t}|_2^2=\frac{2}{B_2}=12$.  Thus, we
will see the number 12 coming up in several places.  In particular,
define the \textit{hat-operation} on $\mathbf{t}
\in\widetilde{\ell}_2$ with positive $2$-norm by 
\begin{equation}\label{eq:hat.operation}
\widehat{\mathbf{t}} \coloneqq
\frac{\sqrt{12} \cdot \mathbf{t}}{|\mathbf{t}|_2},
\end{equation}
so that $\Var[\cS_{\widehat{\mathbf{t}}}] = 1$ and
$\cS_{\widehat{\mathbf{t}}}=\cS_{\widehat{\mathbf{t}}}^{*}$.

Now, we can return to the limiting distributions of forests in the
``degenerate'' regime.  We show in \Cref{rem:standard.trees} that it
suffices to consider only standardized trees in order to characterize
all of $\bfM_{\Forest}$.  In \Cref{def:elevation}, we associate to
each tree $P$ an \textit{elevation multiset} $\mathbf{e}$ depending on
a maximal chain in $P$.  These multisets determine a new type of
limiting distribution related to the generalized uniform sum
distributions, but with another normal summand.

\begin{Theorem}\label{thm:forest_infty}
  Let $P$ be an infinite sequence of standardized trees with $|P| -
\rank(P) = o(|P|^{1/2})$. Then $\cX_P^*$ converges in distribution if
and only if the multisets $\widehat{\mathbf{e}}$ converge pointwise to
some element $\mathbf{t} \in \widetilde{\ell}_2$. In that case, the
limit distribution is $\cS_{\mathbf{t}}+\cN(0,\sigma^{2})$ where
$|\mathbf{t}|_2^2/12 + \sigma^2 = 1$.
\end{Theorem}

Inspired by \Cref{thm:forest_infty}, we begin the study of
\textit{DUSTPAN distributions} associated to random variables of the
form $\cS_{\mathbf{t}}+\cN(0,\sigma^{2})$, assuming the two random
variables are independent, $\mathbf{t} \in\widetilde{\ell}_2$, and
$\sigma \in \mathbb{R}_{\geq 0}$.  The nomenclature DUSTPAN refers
to a \textit{\underline{d}istribution associated to a \underline{u}niform
\underline{s}um for $\mathbf{\underline{t}}$ \underline{p}lus
\underline{a}n independent \underline{n}ormal distribution}.
The generalized uniform sum distributions with
variance 1 are the special case when $\sigma =0$.  Let
\begin{equation}\label{eq:dustpan.param}
 \bfP_{\DUSTPAN}
\coloneqq \left\{\mathbf{t} \in \widetilde{\ell}_2 :
|\mathbf{t}|_2^{2} \leq 12\right\}
\end{equation}
be the \textit{standardized DUSTPAN parameter space}, considered as a sequence
space with the topology of pointwise convergence. Define  the \textit{metric space of
standardized DUSTPAN distributions} to be 
\begin{equation}\label{eq:dustpan.moduli}
\bfM_{\DUSTPAN} \coloneqq \{\cS_\mathbf{t} + \cN(0, \sigma^{2}) :
|\mathbf{t}|_2^2/12 + \sigma^2 = 1\}.
\end{equation}

The standardized DUSTPAN parameter space $\bfP_{\DUSTPAN}$ is a closed subset of
the sequence space $\widetilde{\ell}_2 \subset \bR^{\mathbb{N}}$ considered as a
\textit{Fr\'echet space} (rather than a Banach space). See
e.g.~\cite[Ex.~5.18(1)]{MR1483073} for more details on this structure.
In fact, $\bfM_{\DUSTPAN}$ is closed as well, and we will show we have
the following homeomorphism of compact spaces.

\begin{Theorem}\label{thm:Pinfty_Minfty}
  The map $\Phi \colon \bfP_{\DUSTPAN} \to \bfM_{\DUSTPAN}$ given by
$\mathbf{t} \mapsto \cS_\mathbf{t} + \cN(0, \sigma^{2})$ where $\sigma
\coloneqq \sqrt{1 - |\mathbf{t}|_2^2/12}$ is a homeomorphism between
compact spaces.
\end{Theorem}

\begin{Corollary}\label{cor:limit.laws.sums}
The limit laws for all possible standardized general uniform
sum distributions $\bfM_{\SUMS} \coloneqq \{\cS_\mathbf{t}^{*}:
\mathbf{t} \in
\widetilde{\ell}_2 \}$
is exactly the 
\textit{metric space of DUSTPAN distributions,} \[
\overline{\bfM_{\SUMS}} = \bfM_{\DUSTPAN}. \]
\end{Corollary}

\begin{Corollary}\label{cor:forest_infty}
  For any fixed $\epsilon>0$, let ${\epsilon\Tree}$ be the set of
standardized trees $P$ for which $|P| - \rank(P) <
|P|^{\frac{1}{2}-\epsilon}$.  Let $\bfM_{\epsilon\Tree}
\coloneqq \{\cX_P^* : P \in {\epsilon\Tree} \} \subset \bfM_{\Forest}
$ be the corresponding metric space of distributions. Then
\begin{equation} \overline{\bfM_{\epsilon\Tree}}
= \bfM_{\epsilon\Tree} \sqcup \bfM_{\DUSTPAN},
\end{equation}
which is compact. Moreover, the set of limit points of
$\bfM_{\epsilon\Tree}$ is $\bfM_{\DUSTPAN}$.
\end{Corollary}

\begin{Remark}\label{rem:manifolds}
  The foundational idea of information geometry is to endow
  spaces of distributions with the structure of
  Riemannian manifolds. Consequently, one may be tempted
  to recast \Cref{thm:Pinfty_Minfty} in the context of manifold
  theory. However, the infinite-dimensional case is generally
  ``not mathematically easy'' \cite[\S2.5, p.39]{MR3495836}.
  Here, $\widetilde{\ell}_2$ is a Hilbert manifold and a Banach
  manifold under the $\widetilde{\ell}_2$-norm, as well as a Fr\'echet
  manifold under pointwise convergence. There does not
  appear to be a generally agreed-upon Hilbert, Banach,
  or Fr\'echet manifold structure  which the closed subset
  $\bfP_{\DUSTPAN}$ inherits from $\widetilde{\ell}_2$, though it could
  perhaps be thought of as a manifold with corners. In any case,
  the inherited Hilbert and Banach topology on $\bfP_{\DUSTPAN}$
  disagrees with the Fr\'echet topology, so for our
  purposes, \Cref{thm:Pinfty_Minfty} requires us to use the
  Fr\'echet structure of pointwise convergence. It is
  consequently unclear if a useful differentiable structure
  exists
  for  $\bfP_{\DUSTPAN}$
\end{Remark}

As with $\bfM_{\SSYT}$, it remains an open problem to completely
classify all possible limit points of $\bfM_{\Forest}$.  The
strongest results we have proven for $q$-hook length formulas for
forests show $\bfM_{\Forest} \cup \bfM_{\DUSTPAN} \subset
\overline{\bfM_{\Forest}}$, implying there are an uncountable number
of possible limit laws for distributions associated to forests.  In
the case of forests, the underlying distributions are always symmetric
and unimodal, in contrast to $\bfM_{\SYT}$ which are not always
unimodal, see \cite[Conj.~8.1]{BKS-asymptotics}.  So,
$\overline{\bfM_{\Forest}}$ does not contain $\bfM_{\SYT}$.

More generally, it is natural to ask which limit laws are possible for
the coefficients of arbitrary $q$-hook-type formulas, namely
polynomials with nonnegative integer coefficients of the form
$\prod_{i=1}^{n}[a_{i}]_q/[b_{i}]_q$.  In \cite{BS-cgfs}, we call such
$q$-integer quotients \textit{cyclotomic generating functions} (CGF's)
and study their properties from a variety of algebraic and
probabilistic perspectives. Let $\bfM_{\CGF}$ denote the corresponding
metric space of standardized distributions. By Prohkorov's Theorem,
$\overline{\bfM_{\CGF}}$ is compact.

\begin{Problem}
  Describe $\overline{\bfM_{\CGF}}$ in the L\'evy metric. What are all
possible limit points?  Is $\bfM_{\CGF} \cup
\bfM_{\DUSTPAN}$ \underline{the} metric space of limit laws for $q$-hook formulas,
referring back to the title of this article?
\end{Problem}

\subsection{Paper organization}

The rest of the paper is organized as follows. In
\Cref{sec:background}, we provide background for the hook and cumulant
formulas summarized in \Cref{tab:formulas}. In \Cref{sec:sums}, we
analyze the metric space of generalized uniform sum distributions and
its variations in order to prove \Cref{thm:Pinfty_Minfty} and its
analog for the distance distributions. The analysis of $\bfM_{\SSYT}$
and $\bfM_{\PP}$ is in \Cref{sec:ssyt}. The analysis of
$\bfM_{\Forest}$ is in \Cref{sec:forests}.  Some additional open
questions and avenues for future work are listed in \Cref{sec:future}.

\section{Background}\label{sec:background}

In this section, we briefly recall statements from the literature we
will need related to asymptotic distributions, semistandard tableaux
and forests. All of our arguments for determining asymptotic
distributions use the Method of Moments/Cumulants. Using work of
Hwang--Zacharovas, we explain a key insight for this paper, namely
that rational product formulas such as appear in \Cref{tab:formulas}
give rise to explicit formulas for cumulants of the corresponding
distributions.  See \cite[\S2-3]{BKS-asymptotics} for a more extensive
exposition aimed at an audience familiar with enumerative
combinatorics.  See \cite{MR1324786} for background in probability.

\subsection{Asymptotic distributions}\label{sub:asymptotic}

Let $\cX$ be a real-valued random variable.  For $d \in \bZ_{\geq 0}$,
the \textit{$d$th moment} $\cX$ is \[ \mu_d \coloneqq \bE[\cX^d].\]
The \textit{moment-generating function} of $\cX$ is \[
M_{\cX}(t)\coloneqq\bE[e^{t\cX}] = \sum_{d=0}^\infty \mu_d
\frac{t^d}{d!}, \] which for us will always have a positive radius of
convergence.  The \textit{characteristic function} of $\cX$ is \[
\phi_{\cX}(t)\coloneqq\bE[e^{it\cX}], \] which exists for all $t \in
\bR$ and which is the Fourier transform of the density or mass
function associated to $\cX$.  We will need the following technical
details for the proofs in future sections.

\begin{Remark}\label{rem:analytic_cf}
The characteristic function $\phi_{\cX}(s) \coloneqq \bE[e^{is\cX}]$
in general converges only for $s \in \bR$. However, if there is a
complex analytic function $\psi(s)$ defined in an open ball $|s|<\rho$
such that $\phi_{\cX}(s) = \psi(s)$ for $-\rho < s < \rho$, then
$\phi_{\cX}(s)$ exists and is analytic in some strip $-\beta <
\mathrm{Im}(s) < \alpha$ where $\alpha, \beta \geq \rho$. Moreover,
for $|s| < \rho$, $\phi_{\cX}(s) = \psi(s)$.  In particular, the
moment-generating function $\bE[e^{t\cX}]$ converges for $-\rho < t <
\rho$, so $\cX$ has moments of all orders and is determined by its
moments. See e.g.~\cite[Thm.~7.1.1, pp.191-193]{MR0346874} and
\cite[Thm.~30.1]{MR1324786} for details.
\end{Remark}

The \textit{cumulants} $\kappa_1, \kappa_2, \ldots$ of $\cX$ are
defined to be the coefficients of the exponential generating
function \[ K_{\cX}(t) \coloneqq \sum_{d=1}^\infty \kappa_d
\frac{t^d}{d!}  \coloneqq \log M_{\cX}(t) = \log \bE[e^{t\cX}]. \]
Hence, they satisfy the recurrence
\begin{align}\label{eq:moment_to_cumulant}
  \mu_d = \kappa_d + \sum_{m=1}^{d-1}
    \binom{d-1}{m-1} \kappa_m \mu_{d-m}, 
\end{align}
so the moments can similarly be recovered from the cumulants and vice
versa.  In particular, \eqref{eq:moment_to_cumulant} implies
$\kappa_{1}=\mu_{1}=\mu=\bE[\cX]$ and $\kappa_{2}=\Var[\cX]=\sigma^{2}$. The
cumulants also satisfy
\begin{enumerate}
\item \textit{(Homogeneity)}: the $d$th cumulant of $c\cX$ is
$c^d\kappa_d$ for $c \in \bR$, and 
\item \textit{(Additivity)} the cumulants of the sum of
\textit{independent} random variables are the sums of the cumulants.
\end{enumerate}
For $d\geq 4$, the moments of independent random variables are not
necessarily sums of the moments, so cumulants work much better for our
purposes.  By homogeneity and additivity, the associated \textit{standardized random
variable} $\cX^* \coloneqq (\cX - \mu)/\sigma$ has cumulants
$\kappa_1^{\cX^*}=0$, $\kappa_2^{\cX^*}=1$, and 
\begin{equation}\label{eq:rescaled.cumulants}
\kappa_d^{\cX^*} =
\frac{\kappa_d^{\cX}}{\sigma^{d}}=\frac{\kappa_d^{\cX}}{(\kappa_2^{\cX})^{d/2}}
\hspace{.1in} \text{ for $d\geq 2$. }
\end{equation}

\begin{Example}\label{ex:normal.cum}
The \textit{normal distribution} $\cN(0,1)$ is the unique distribution
with $\kappa_1=0$, $\kappa_2=1$, and $\kappa_d=0$ for $d\geq 3$.
Therefore, $\cN(\mu ,\sigma^{2})$ is the unique distribution with
cumulants $\kappa_1=\mu $, $\kappa_2=\sigma^{2}$, and $\kappa_d=0$ for
$d\geq 3$.
\end{Example}

\begin{Example}\label{ex:Ucont}
  Let $\cU = \cU[0, 1]$ be the continuous uniform random variable
whose density takes the value $1$ on the interval $[0,1]$ and $0$
otherwise.  Then the moment generating function is $M_{\cU}(t) =
\int_{0}^1 e^{tx}dx = (e^{t} - 1)/t$, so the cumulant generating
function $\log M_{\cU}(t)$ coincides with the exponential generating
function for the \textit{divided Bernoulli numbers} $\frac{B_d}{d}$
for $d \geq 1$. Their exponential generating function $E_D(t)$
satisfies \[ E_D(t) \coloneqq \sum_{d \geq 1} \frac{B_d}{d}
\frac{t^d}{d!}  = \log\left(\frac{e^t - 1}{t}\right). \] Hence, the
$d^{th}$ cumulant for $\cU$ is $\kappa_d^{\cU} = B_d/d$ for $d \geq
1$.  Recall from \Cref{sec:intro}, $\cIH_m$ is the
\textit{Irwin--Hall} distribution obtained by adding $m$ independent
$\cU[0, 1]$ random variables.  By additivity, the $d$th
cumulant of $\cIH_m$ is $mB_d/d$.  More generally, let $\cS \coloneqq
\sum_{k=1}^m \cU[\alpha_k, \beta_k]$ be the sum of $m$ independent
uniform continuous random variables. Then the $d$th cumulant of $\cS$
for $d \geq 2$ is
\begin{equation}\label{eq:Ucont.cumulants} \kappa_d^{\cS} =
\frac{B_d}{d} \sum_{k=1}^m (\beta_k - \alpha_k)^d
\end{equation}
by the homogeneity and additivity properties of cumulants.
\end{Example}

The Method of Moments/Cumulants is based on the following theorem.  All random
variables we encounter will have moments of all orders.

\begin{Theorem}[Frech\'et--Shohat Theorem,
{\cite[Theorem~30.2]{MR1324786}}]\label{thm:moments} Let $\cX_1,
\cX_2, \ldots$ be a sequence of real-valued random variables, and let
$\cX$ be a real-valued random variable. Suppose the moments of $\cX_n$
and $\cX$ all exist and the moment generating functions all have
positive radius of
convergence. If \begin{equation}\label{eq:moments_criterion} \lim_{n
\to \infty} \mu_d^{\cX_n} = \mu_d^{\cX} \hspace{.5cm} \forall d \in
\bZ_{\geq 1}, \end{equation} then $\cX_1, \cX_2, \ldots$ converges in
distribution to $\cX$.  Similarly, if
\begin{equation}\label{eq:cumulants_criterion}
  \lim_{n \to \infty} \kappa_d^{\cX_n} = \kappa_d^{\cX}  \hspace{.5cm} \forall d \in
\bZ_{\geq 1},
  \end{equation}
then $\cX_1, \cX_2, \ldots$ converges in
distribution to $\cX$.
\end{Theorem}

\begin{Corollary}\label{cor:cumulants}
  A sequence $\cX_1, \cX_2, \ldots$ of real-valued
  random variables on finite sets is asymptotically normal
  if for all $d \geq 3$ we have
  \begin{equation}\label{eq:cumulants_criterion2}
 \lim_{n \to \infty} \kappa^{\cX_n^{*}}_d =  \lim_{n \to \infty} \frac{\kappa^{\cX_n}_d}{(\sigma^{\cX_n})^d} = 0.
  \end{equation}
\end{Corollary}

For a positive integer $n$, define the associated \textit{
$q$-integer} to be the polynomial $[n]_{q}=1+q+q^{2}+\dots +q^{n-1} =
(1-q^{n})/(1-q)$.  The $q$-integers factor into cyclotomic polynomials
over the integers.  Therefore, the hook length formulas considered in
this paper are all products of cyclotomic polynomials.  Because these
rational product formulas are polynomial, all cancellation can be done
efficiently by taking the multiset difference between the numerator
and denominator of the cyclotomic factors.

In a forthcoming paper \cite{BS-cgfs}, we investigate general
properties of generating functions which are products of cyclotomic
polynomials with nonnegative coefficients.  For this paper, we just
need two facts.  The first theorem first appeared explicitly in the
work of Hwang--Zacharovas \cite[\S4.1]{MR3346464} building on the work
of Chen--Wang--Wang \cite[Thm. 3.1]{Chen-Wang-Wang.2008}, who in turn
used an argument going back at least to Sachkov
\cite[\S1.3.1]{MR1453118}.

\begin{Theorem}\cite[\S4.1]{MR3346464}\label{thm:CGF_cumulant}
  Suppose $\{a_1, \ldots, a_m\}$ and $\{b_1, \ldots, b_m\}$ are
multisets of positive integers such that 
\begin{equation}\label{eq:cgf.form}
P(q) = \prod_{k=1}^m
\frac{1 - q^{a_k}}{1 - q^{b_k}}= \prod_{k=1}^m
\frac{[a_{k}]_{q}}{[b_k]_{q}} = \sum_k c_k q^k \in \bZ_{\geq 0}[q].
\end{equation}
Let $\cX$ be a discrete random variable with $\bP[X=k] =
c_k/P(1)$. Then the $d$th cumulant of $\cX$
is \begin{equation}\label{eq:CGF_cumulant} \kappa_d^\cX =
\frac{B_d}{d} \left(\sum_{k=1}^m a_k^d - b_k^d\right) \end{equation}
where $B_d$ is the $d$th Bernoulli number (with $B_2 = \frac{1}{2}$).
\end{Theorem}

The following corollary is proved in \cite{BS-cgfs}. It also
follows from the tail decay bound in \cite[Lemma~2.8]{MR3346464}.
We need this for our current investigations for hook length formulas.

\begin{Lemma}[Converse of Frech\'et--Shohat for CGF's]\label{cor:frechet_shohat_converse}
  Suppose $\cX_1, \cX_2, \ldots$ is a sequence of random variables
corresponding to polynomials of the same form as \eqref{eq:cgf.form}.
If $\cX_n^* \Rightarrow \cX$ for some random variable $\cX$, then
$\cX$ is determined by its cumulants and, for all $d \in \bZ_{\geq
1}$, \[ \lim_{n \to \infty} \kappa_d^{\cX_n^*} = \kappa_d^{\cX}. \]
\end{Lemma}

\subsection{Semistandard Young tableaux and plane partitions}\label{ssec:back:ssyt}

We briefly recall the definition and notation for Schur functions,
semistandard tableaux and plane partitions.  For more information on
symmetric functions and their connection with the enumeration of plane
partitions and tableaux, see \cite[Ch. 7]{ec2}.

  A partition $\lambda =(\lambda_{1}\geq \lambda_{2}\geq \dots \geq
\lambda_{k} )$ is a finite decreasing sequence of positive integers.
Let $\ell(\lambda ) =k$ denote the \textit{length} of $\lambda$. We
think of $\lambda$ in terms of its Young diagram, which is a left
justified array of $\ell(\lambda )$ rows with $\lambda_{i}$ cells on
row $i$ and index the cells in matrix notation.

A \textit{semistandard Young tableau}, or just semistandard tableau
for short, of shape $\lambda$ is a filling
of the cells of $\lambda$ with positive integer labels, possibly
repeated, such that the labels weakly increase to the right in rows
and strictly increase down columns. The set of semistandard Young
tableaux of shape $\lambda$ is denoted $\SSYT(\lambda)$.  The subset
of $\SSYT(\lambda)$ filled with integers no greater than $m$ is
denoted $\SSYT_{\leq m}(\lambda)$, which is a finite set.  The
\textit{type} of a semistandard tableau $T$ is the composition
$\alpha(T) = (\alpha_1, \alpha_2, \ldots)$ where $\alpha_i$ is the
number of times $i$ appears in $T$.  The \textit{Schur function}
\[
s_{\lambda}(x_1, x_2,\ldots) \coloneqq  \sum_{T \in \SSYT(\lambda)}  \mathbf{x}^{\alpha(T)}
\]
is the \textit{type generating function} for all semistandard tableaux
of shape $\lambda$, where $\mathbf{x}^{\alpha} \coloneqq
x_1^{\alpha_1} x_2^{\alpha_2} \cdots$.

The \textit{rank} of a semistandard tableau $T$ is a nonnegative
integer statistic depending only on the type.  It is defined by \[
\rank(T) \coloneqq \rank(\alpha) \coloneqq \sum_{i
\geq 1} (i-1) \alpha_i. \] For example, for a fixed partition
$\lambda$, the smallest possible rank of any $T \in \SSYT(\lambda)$
occurs for the tableau with all 1's in the first row, all 2's in the
second row, etc. in the diagram of $\lambda$.  Therefore, the minimal
rank is $\sum (i-1)\lambda_{i}$, which we denote as $\rank(\lambda )$.
The \textit{rank generating function} for $\SSYT(\lambda)$ is given by
the \textit{principal specialization of the Schur function},
\begin{align*}
  s_{\lambda}(1, q, q^2, \ldots)
    &= \SSYT(\lambda)^{\rank}(q) =  \sum_{T \in \SSYT(\lambda)}  q^{\rank (T)} \\
  s_{\lambda}(1, q, q^2, \ldots, q^{m-1}) &= \SSYT_{\leq
m}(\lambda)^{\rank}(q)  \sum_{T \in \SSYT(\lambda)_{\leq m}}  q^{\rank (T)}.
\end{align*}
The motivation for considering this particular specialization comes
from the $q$-analog of the Weyl dimension formula in representation
theory. Stembridge \cite[\S2.2-2.3, Prop.~2.4]{MR1262215} put a ranked
poset structure on the weights of a semisimple Lie algebra, which in
type $A$ reduces to $\rank(\alpha)$.  The following rational product
formula for $s_\lambda(1, q, q^2, \ldots, q^{m-1})$ follows easily
from the classical ratio of determinants definition of Schur
polynomials.

\begin{Theorem}[{\cite[\S7.1]{MR0002127}}, {\cite[(7.105)]{ec2}}]\label{thm:SSYT_rank_gf.1}
  For any partition $\lambda$ and positive integer  $m \geq  \ell(\lambda)$, 
  \begin{equation}\label{eq:SSYT_rank_gf.1}
    s_\lambda(1, q, q^2, \ldots, q^{m-1})
      = q^{\rank(\lambda)}
         \prod_{1 \leq i < j \leq m} \frac{[\lambda_i - \lambda_j + j - i]_q}
         {[j - i]_q}.
  \end{equation}
\end{Theorem}

Stanley gave an alternate rational product formula for $ s_\lambda(1,
q, \ldots, q^{m-1})$, which is called the \textit{$q$-hook-content
formula}. Here the \textit{content} of a cell $u$ in row $i$, column
$j$ in $\lambda$ is defined as $c_u \coloneqq j-i$. Also, the
\textit{hook length} of cell $u$, denoted $h_{u}$, is the number of
cells directly east of $u$, plus the number of cells directly south of
$u$ in the diagram of $\lambda$.

\begin{Theorem}[{\cite[Thm.~7.21.2]{ec2}}]\label{thm:SSYT_rank_gf.2}
  For any partition $\lambda$ and positive integer  $m \geq
  \ell(\lambda)$, 
  \begin{equation}\label{eq:SSYT_rank_gf.2}
    s_\lambda(1, q, \ldots, q^{m-1})
      = q^{\rank(\lambda)} \prod_{u \in \lambda} \frac{[m + c_u]_q}{[h_u]_q}.
  \end{equation}
\end{Theorem}

The two product formulas for $s_\lambda(1, q, \ldots, q^{m-1})$ are
each useful in different circumstances. The product in
\eqref{eq:SSYT_rank_gf.1} involves $\binom{m}{2}$ terms, whereas the
product in \eqref{eq:SSYT_rank_gf.2} involves $|\lambda|$ terms.
One can observe from these formulas that $ s_\lambda(1, q, \ldots,
q^{m-1})$ is symmetric about the mean nonzero coefficient.  From the
representation theory of $\GL_{2}(\mathbb{C})$, it is known that $
s_\lambda(1, q, \ldots, q^{m-1})$ is also unimodal. See
\cite{GoodmanOharaStanton1992} for a combinatorial proof relying on
the unimodality of the Gaussian polynomials.

Recently, Huh--Matherne--M\'esz\'aros--St.Dizier \cite{huh2019logarithmic}
showed that Schur
polynomials are strongly log-concave.  However, we note that
$s_\lambda(1, q, \ldots, q^{m-1})$ is not always log-concave.  For
example,
\[
s_{(3,1)}(1, q, q^{2}, q^{3})= 
q^{10} + 2q^{9} + 4q^{8} + 5q^{7} + 7q^{6} + 7q^{5} + 7q^{4} + 5q^{3} + 4q^{2} + 2q^{1} + 1,
\]
which is not log-concave since $5^{2}<4\cdot 7$.

Combining \Cref{thm:CGF_cumulant} and \Cref{thm:SSYT_rank_gf.2}, we
get an exact formula for the cumulants of the random variable
associated to the rank function on semi-standard Young tableaux on the
alphabet $[m]$ chosen uniformly.  This cumulant formula is the key to
analyzing the asymptotic distributions.  

\begin{Corollary}\label{cor:ssyt.cgf} Fix a partition $\lambda$.  If
$\kappa_d^{\lambda; m}$ is the $d$th cumulant of the random variable
associated to $\rank$ on $\SSYT_{\leq m}(\lambda)$, then, for $d >
1$, \begin{align} \label{eq:SSYT_cumulants.1} \kappa_d^{\lambda; m} &=
\frac{B_d}{d} \left(\sum_{1 \leq i < j \leq m}
           (\lambda_i - \lambda_j + j - i)^d - (j - i)^d\right) \\
    \label{eq:SSYT_cumulants.2}
      &= \frac{B_d}{d} \left(\sum_{u \in \lambda} (m+c_u)^d - h_u^d\right).
  \end{align}
\end{Corollary}

Observe, the summands in \eqref{eq:SSYT_cumulants.2} can be negative,
but the summands in \eqref{eq:SSYT_cumulants.1} are each clearly
positive.  Thus, $\kappa_d^{\lambda; m}$ has the same sign as the
Bernoulli number $B_{d}$, namely it is negative if and only if $d$ is
divisible by 4, and $\kappa_d^{\lambda; m}=B_{d}=0$ if and only if
$d>1$ and odd.

\begin{Definition}
  A \textit{plane partition} is a finite collection of unit cubes in
the positive orthant of $\bR^3$ stacked towards the origin. More
formally, it is a finite lower order ideal in $\bZ_{\geq 1}^3$ under
the component-wise partial order. We may imagine a plane partition
$\rho$ as a matrix with entry $\rho_{ij}$ recording the number of
cells with $x$-coordinate $i$ and $y$-coordinate $j$.  The
\textit{size} of a plane partition $\rho$ is the number of cubes,
denoted $|\rho|=\sum \rho_{ij}$.  We write $\PP(a \times b \times c)$
for the set of all plane partitions fitting inside an $a$ by
$b$ by $c$ rectangular prism.
\end{Definition}

There is a straightforward bijection between plane partitions and
semistandard Young tableaux of rectangular shape,	      
\begin{equation}
  \begin{split}
    \PP(a \times b \times c) &\too{\sim} \SSYT_{\leq a+c}((b^a)) \\
    \rho &\mapsto T \text{ where }T_{ij} = c - \rho_{ij} + i.
  \end{split}
\end{equation}
All $T$ and $\rho$ in the bijection are rectangular arrays with $a$
rows and $b$ columns with entries labeled using matrix indexing
conventions.  Letting $|T| \coloneqq \sum_{i,j} T_{ij}$, note that
$|T| = \rank(T) + ab$ and $|T| + |\rho| = abc + b\binom{a+1}{2}$ is
constant. Hence, the unique element of minimal size in $\PP(a \times b
\times c)$, namely $\varnothing$, maps to the unique maximal rank
tableau in $\SSYT_{\leq a+c}((b^a))$ with values $c-i$ in row $i$ for
each $1\leq i\leq a$.   

By \Cref{thm:SSYT_rank_gf.1} and \Cref{thm:SSYT_rank_gf.2}, we know
$\SSYT(\lambda)^{\rank}(q)$ is symmetric up to an overall $q$-shift.
Similarly, $\PP(a \times b \times c)$ is closed under box
complementation, so it follows from the bijection and
\eqref{eq:SSYT_rank_gf.1} that
\begin{align}\label{eq:ppformula}
  \PP(a \times b \times c)^{\size}(q)
    &= q^{-\rank(\lambda)} \SSYT_{\leq a+c}((b^a))^{\rank}(q)
    \\
    &= \prod_{i=1}^a \prod_{j=1}^b
       \frac{[a+c+ j-i]_q}{[a+b-i-j+1]_q}
       \\
    &= \prod_{i=1}^a \prod_{j=1}^b
       \frac{[i+j+c-1]_q}{[i+j-1]_q}
       \\
    &= \prod_{i=1}^a \prod_{j=1}^b \prod_{k=1}^c
       \frac{[i+j+k-1]_q}{[i+j+k-2]_q}.
\end{align}
The later two product formulas are originally due to MacMahon.  See
the proof of \cite[Thm. 7.21.7]{ec2} for more details and
\cite[pp.~402-403]{ec2} for historical references. In particular, the
cumulants of $\size$ on $\PP(a \times b \times c)$ are given by
\eqref{eq:SSYT_cumulants.1} or \eqref{eq:SSYT_cumulants.2} where
$\lambda = (b^a)$ and $m = a+c$.

\subsection{Linear extensions of forests}\label{ssec:forests.background}

Next, we summarize the relevant terminology and results from
\cite{MR1022316}.  Briefly recall, a \textit{tree} is a finite,
connected simple graph with no cycles.  A \textit{forest} is a finite
disjoint union of trees.  A tree is \textit{rooted} if it has a
distinguished vertex, called the root. A forest is
\textit{rooted} if each of its trees is rooted.  The \textit{Hasse
diagram} of a partially ordered set (\textit{poset}) $P$ is the
graph with vertex set $P$ where there is an edge between $x$
and $y$ if $y$ \textit{covers} $x$, i.e.~$x <_P y$ and there does not exist $u \in P$
such that $x <_P u <_P y$. We refer to a poset as a \textit{forest} if
its Hasse diagram is a forest with roots as maximal elements, or
equivalently if every element of $P$ is covered by at most one
element.

\begin{Definition}
  Let $P$ be a finite partially ordered set. The \textit{rank} of $P$
is the maximum number of elements in any chain $u_1 < u_2 < \cdots <
u_{k}$ in $P$. For instance, if $P$ is a singleton, its rank is $1$.
Note that this definition is one larger than the standard definition in
\cite[Ch.3]{ec1}, but it is more convenient for our purposes.
\end{Definition}

\begin{Definition}
  Let $P$ be a poset. A \textit{labeling} of $P$ is a bijection $w
\colon P \to [n]$, and a \textit{labeled poset} is a pair $(P, w)$
where $w$ is a labeling of $P$.  A labeling $w$ of $P$ for which $w(p)
\leq w(q)$ whenever $p \leq_P q$ is called a \textit{natural
labeling}.  A labeling $w$ of $P$ is \textit{regular} if for all $x
<_P z$ and $y \in P$, if $w(x) < w(y) < w(z)$ or $w(x) > w(y) > w(z)$
then $x <_P y$ or $y <_P z$.
Regular labelings of forests include the postorder, preorder, and
inorder labelings, which are commonly used in computer science.
\end{Definition}

\begin{Definition}
  A \textit{linear extension} of $P$ is  an
  ordered list   $p_1, \ldots, p_n$ of the elements of $P$
 such that $i \leq j$ whenever $p_i \leq_P p_j$.
  If $(P, w)$ is a labeled poset, a linear extension can be thought
  of as the permutation $i \mapsto w(p_i)$ of $[n]$.
  The set $\cL(P, w)$ is the set of all permutations
  obtained in this fashion from linear extensions of the
  labeled poset $(P, w)$.
\end{Definition}

It is often convenient to use a natural labeling $w$ of $P$ so that
$\id \in \cL(P, w)$.  Choosing labelings which are not natural forces
inversions to appear in any $\sigma \in \cL(P, w)$.  Finding the
minimum number of inversions in any linear extension of an arbitrarily
labeled poset motivates the following analogues related to inversions
and descents in permutations.  

\begin{Definition}\label{def:inv.maj.posets}
  Let $(P, w)$ be a labeled poset. Set
  \begin{align*}
    \Inv(P, w)
      &\coloneqq \{(w(x), w(y)) : x <_P y\text{ and }w(x) > w(y)\}
      &\text{(\textit{inversion set})} \\
    \inv(P, w)
      &\coloneqq |\Inv(P, w)|
      &\text{(\textit{inversion number})} \\
    \Des(P, w)
      &\coloneqq \{w(x) : w(x) > w(y), y\text{ covers }x  \in P \}
      &\text{(\textit{descent set})} \\
    \maj(P, w)
      &\coloneqq \sum_{x \in \Des(P, w)} h_x
      &\text{(\textit{major index})} 
  \end{align*}
where the \textit{hook length} of an element $x \in P$ is
\begin{equation}\label{def:poset_hooklength}
h_x \coloneqq \#\{t \in P : t \leq_P x\}. 
\end{equation}
\end{Definition}

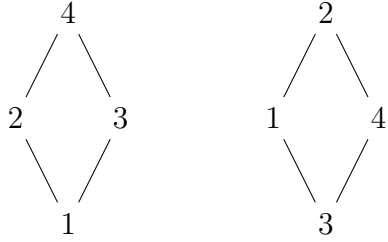
\begin{figure}
\begin{center}
\begin{tikzpicture}[scale=.7]
  \node (one) at (0,2) {$4$};
  \node (b) at (-1,0) {$2$};
  \node (c) at (1,0) {$3$};
  \node (zero) at (0,-2) {$1$};
  \draw (zero) -- (b) -- (one) -- (c) -- (zero);
\end{tikzpicture}
\hspace{.5in}
\begin{tikzpicture}[scale=.7]
  \node (one) at (0,2) {$2$};
  \node (b) at (-1,0) {$1$};
  \node (c) at (1,0) {$4$};
  \node (zero) at (0,-2) {$3$};
  \draw (zero) -- (b) -- (one) -- (c) -- (zero);
\end{tikzpicture}

\end{center}
\caption{A naturally labeled poset $(P,v )$ on the left and another
labeling of the same diamond poset $(P,w )$ on the right which is
not natural or regular.}  \label{fig:diamond}
\end{figure}

\begin{Example}\label{example:diamond}
For the first labeled poset $(P,v )$ in \Cref{fig:diamond}, we have
$\cL(P, v)=\{1234, 1324 \}$,\ $ \Inv(P, v) =\Des(P, v)= \emptyset $, and
$ \inv(P, v) = \maj(P, v)=0$.  For the second labeled poset $(P,w )$
in \Cref{fig:diamond}, we have $\cL(P, w)=\{3142, 3412 \}$,\ $ \Inv(P,
w) =\{(3,1),(3,2),(4,2)\}$, \ $\Des(P, w)= \{3,4 \} $, $ \inv(P, w) =
\maj(P, w) =3$.   The hook lengths of the diamond poset are
$1,2,2,4$.
\end{Example}

\begin{Remark}
  One can consider a partition $\lambda$ as a poset on its cells where $(u, v) \leq
(x, y)$ if and only if $u \leq x$ and $v \leq y$.  However, the hook
lengths of $\lambda$ do \textit{not} agree with
\eqref{def:poset_hooklength} except when $\lambda$ is a single row or
column.  For example, the hook lengths for the partition $(2,2)$ are
$1,2,2,3$, while the hook lengths for the diamond poset are $1,2,2,4$.
\end{Remark}

Mallows and Riordan first studied the inversion enumeration on labeled rooted
trees \cite{RiordanMallows.1968}, and connected it to cumulants of the
lognormal distribution.
Knuth gave a hook length formula for $|\cL(P, w)|$ \cite[p.~70]{Knuth}
for posets which are forests.  Bj\"orner--Wachs \cite{MR1022316} and
Stanley \cite{MR0332509} generalized Knuth's result to $q$-hook length
formulas using the $\inv$ and $\maj$ statistics on $\cL(P,
w)$. Stanley considered only the case when $w$ is natural, i.e.~when
$\inv(P, w) = \maj(P, w) = 0$, for the $\maj$ generating function.
Recently Zaguia has studied linear extensions of forests and proved
the ``1/3-2/3 Conjecture''  holds on such posets \cite{Zaguia.2018}.

\begin{Theorem}[{\cite[Thm.~1.1-1.2, Cor.~3.1,
Thm.~6.1-6.2]{MR1022316}}]
\label{thm:bw_invmaj_forest}
  Let $(P, w)$ be a labeled poset with $n$ elements.  Then \[ \cL(P,
w)^{\maj}(q) \coloneqq  \sum_{\pi \in \cL(P, w)} q^{\maj (\pi )}=q^{\maj(P, w)}
\frac{[n]_q!}{\prod_{u \in P} [h_u]_q} \] if and only if $P$ is a
forest. Similarly, 
\[ \cL(P, w)^{\inv}(q) \coloneqq 
\sum_{\pi \in \cL(P, w)} q^{\inv (\pi )}= q^{\inv(P, w)}
\frac{[n]_q!}{\prod_{u \in P} [h_u]_q} \] if and only if $(P,w)$ is a
regularly labeled forest.  Moreover, if $P$ is a forest,
$\frac{[n]_q!}{\prod_{u \in P} [h_u]_q}$ has symmetric and unimodal
coefficients.
\end{Theorem}

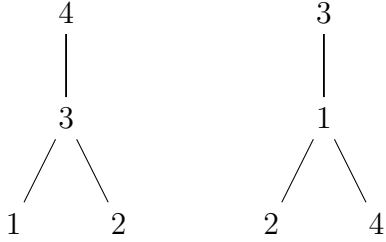
\begin{figure}
\begin{center}
\begin{tikzpicture}[scale=.7]
  \node (one) at (0,4) {$4$};
  \node (c) at (0,2) {$3$};
  \node (a) at (-1,0) {$1$};
  \node (b) at (1,0) {$2$};
  \draw (b) -- (c) -- (one) -- (c) -- (a);
\end{tikzpicture}
\hspace{.5in}
\begin{tikzpicture}[scale=.7]
  \node (one) at (0,4) {$3$};
  \node (c) at (0,2) {$1$};
  \node (a) at (-1,0) {$2$};
  \node (b) at (1,0) {$4$};
  \draw (b) -- (c) -- (one) -- (c) -- (a);
\end{tikzpicture}
\end{center}
\caption{A naturally labeled poset $(P,v )$ on the left and another
labeling of the same forest poset $(P,w )$ on the right which is
not natural.}  \label{fig:d3}
\end{figure}

\begin{Example}\label{example:diamond.3}
For the first labeled poset $(P,v )$ in \Cref{fig:d3}, we have $\cL(P,
v)=\{1234, 2134 \}$,\ $ \Inv(P, v) =\Des(P, v)= \emptyset $, and $
\inv(P, v) = \maj(P, v)=0$.  For the second labeled poset $(P,w )$ in
\Cref{fig:d3}, we have $\cL(P, w)=\{2413, 4213 \}$,\ $ \Inv(P, w)
=\{(2,1),(4,1),(4,3)\}$, \ $\Des(P, w)= \{2,4 \} $, $ \inv(P, w) = 3$,
$\maj(P, w) = 2$.  Note $\cL(P, w)^{\maj} = q^{2}+q^{3}$, and $\cL(P,
w)^{\inv} = q^{3}+q^{4}$.  The hook lengths of the underlying poset
are $1,1,3,4$.  One can verify the formulas in
\Cref{thm:bw_invmaj_forest} hold in each of these cases, but they
don't hold for the diamond poset.  
\end{Example}

Given a forest $P$, define the polynomial 
\begin{equation}\label{eq:def.cgf.forests}
\cL_{P}(q) \coloneqq [n]_q!/\prod_{u \in P} [h_u]_q, 
\end{equation}
and let $\cX_P$ the associated random variable.  Note, the
distribution of $\cX_P$ does not depend on the choice of labeling of
the vertices of $P$ since $\cL_{P}(q)$ depends only on the unlabeled
poset structure.  We also get simple formulas for the associated
cumulants in the next two statements.

\begin{Remark}\label{rem:forest_degree}
By the unimodality result in
\Cref{thm:bw_invmaj_forest}, we know $\cL_{P}(q):=[n]_q!/\prod_{u \in
P} [h_u]_q$ has nonzero coefficients in an interval, so it has no internal
zeros.  The degree of $\cL_{P}(q)$ is
\[
\sum_{k=1}^n k - \sum_{u \in P} h_u, 
\]
and the mean of $\cX_P$ is half the degree.    
\end{Remark}

\begin{Corollary}\label{cor:forest_cumulants}
  Let $P$ be a forest with $n$ elements.  Suppose $d \in \bZ_{\geq
2}$. Let $\kappa_d^P$ denote the $d$th cumulant of the random variable
$\cX_{P}$.  Then,
\[ \kappa_d^P = \frac{B_d}{d} \left(\sum_{k=1}^n k^d
- \sum_{u \in P} h_u^d\right). \]
\end{Corollary}

\begin{Remark}\label{rem:standard.trees}
In order to characterize all possible limit laws for the standardized
random variables associated with $\maj$ and $\inv$ on labeled forests,
we only need to consider the set of all distributions associated with
standardized trees as follows.  Given any forest $P$, we may turn $P$
into a tree by adding a new vertex covering the roots of all the trees
of $P$. It is easy to see that the quotient in
\eqref{eq:def.cgf.forests} is unchanged, so the cumulants and the
corresponding distributions are the same.  Similarly, if $P$ is a tree
and the root has exactly one child, we may delete the root while
preserving the fact that $P$ is a tree, and the quotient in
\eqref{eq:def.cgf.forests} is again unchanged.  Consequently, we say a
forest is \textit{standardized} if it is a tree and the root has at
least two children. Therefore,
\[
\bfM_{\Forest} \coloneqq \{\cX_P^* : \text{$P$ is a forest}\} =
\{\cX_P^* : \text{$P$ is a standardized tree}\}.
\]
\end{Remark}

\subsection{Riemann integral estimates}\label{ssec:Riemann}

Many of our theorems depend on approximations using a mixture of
combinatorics and analysis.  In particular, we return to certain basic
sums over and over again.  Let ${\mathbf h}_{d}(a,b) =
\sum_{j=0}^{d}a^{j}b^{d-j}$ denote the complete homogeneous symmetric
function on two inputs.  

\begin{Lemma}\label{lem:riemann_estimate}
  For positive integers $a,b$, and $d>1$,   we have
    \[ \frac{1}{d} \left[(a+b)^d - a^d\right]
        < \sum_{j=a+1}^{a+b} j^{d-1}
        < \frac{1}{d} \left[(a+b)^d - a^d\right]
           + (a+b)^{d-1} - a^{d-1}. \]
  Equivalently,
    \[ \frac{b}{d} \, \mathbf{h}_{d-1}(a+b, a) < \sum_{j=a+1}^{a+b} j^{d-1}
         < \frac{b}{d}\,  \mathbf{h}_{d-1}(a+b, a) + b\, \mathbf{h}_{d-2}(a+b, a). \]

  \begin{proof} Use a Riemann integral estimate.
  \end{proof}
\end{Lemma}

\subsection{Standard notation for approximations}\label{sub:notation}

We use the following standard Bachmann--Landau asymptotic notation
without further comment. We write $f(n)=\Theta(g(n))$ to mean
there exist constants $a,b>0$ such that for $n$ large enough, we have
$ag(n)\leq f(n) \leq bg(n)$. If $f(n)=O(g(n))$, then there exists a
constant $c > 0$ such that for all $n$ large enough,
we have $ f(n) \leq cg(n)$.  On the other hand, if $f(n)=o(g(n))$, then
as $n \to \infty$, we have $ \frac{f(n)}{g(n)} \to 0$.  Similarly,
$f(n)=\omega(g(n))$ implies $ \frac{f(n)}{g(n)} \to \infty$ as $n \to
\infty$, and $f(n) \sim g(n)$ implies $ \frac{f(n)}{g(n)} \to 1$ as
$n \to \infty$.

\section{Metric spaces related to  uniform sum distributions}\label{sec:sums}

Motivated by applications to $\bfM_{\SSYT}$ and $\bfM_{\Forest}$ in
the next two sections, we first analyze the distributions of finite
and infinite sums of uniform continuous random variables. We
parameterize these distributions using certain sequence spaces and
precisely relate weak convergence of the underlying distributions to
pointwise convergence of the parametrizing sequences.  The closure of
the space of all possible distributions associated to standardized sums
of independent uniform random variables leads us to define the
metric space of DUSTPAN distributions.  We also describe a closed
subset of the DUSTPAN distributions related to distance multisets,
which appear in the study of $\bfM_{\SSYT}$.  

\subsection{Generalized uniform sum distributions and decreasing
sequence space}\label{sub:inf.sums}

The Irwin--Hall distributions, also known as \textit{uniform sum
distributions}, are the distributions associated to finite sums of
independent, identically distributed, uniform random variables
supported on $[0,1]$. First, we relax the requirement that they be
identically distributed, and then we relax the requirement that they
are finite sums.

\begin{figure}[ht]
 \includegraphics[height=3cm]{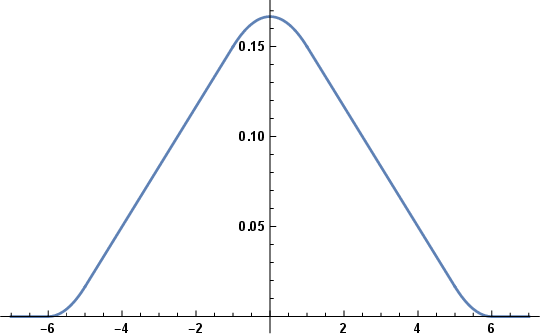}
\hspace{.5in}
 \includegraphics[height=3cm]{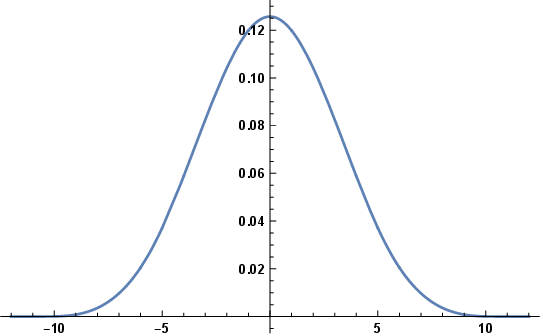}
    \caption{Plots of density functions for the distributions $\cS_{\mathbf{t}}$
with $\mathbf{t}=(6,5,1)$ and $\mathbf{t}=(6,5,5,5,1)$.}
    \label{fig:sums}
\end{figure}

Consider a random variable defined as the sum of $m$ independent
uniform continuous random variables of the form $\cS \coloneqq
\sum_{k=1}^m \cU[\alpha_k, \beta_k]$ with $\alpha_k\leq \beta_k$ for
each $k$.  We call the distribution of $\cS$ a \textit{generalized
uniform sum distribution}.  See \Cref{fig:sums} for example
density functions.  We note that each of the generalized uniform sum
distributions is non-normal, though the histograms may look quite
similar.  By \Cref{ex:Ucont}, the $d$th cumulant of $\cS$ for $d \geq
2$ is
\begin{equation}\label{eq:Ucont.cums}
\kappa_d^{\cS} = \frac{B_d}{d} \sum_{k=1}^m (\beta_k - \alpha_k)^d,
\end{equation}
which only depends on the differences $t_k \coloneqq \beta_k -
\alpha_k$.  It is useful to compare \eqref{eq:Ucont.cums} to the
cumulants in \eqref{eq:CGF_cumulant}.

The random variable $\cS$ can be expressed as a constant overall shift
$c=\frac{1}{2}\sum_{k=1}^{m} (\alpha_k + \beta_k)$ plus a
\textit{uniform sum random variable associated to} $\mathbf{t}$
\begin{equation}\label{eq:t_uniform_sum.b}
  \cS_{\mathbf{t}} \coloneqq \sum_{k=1}^{m}
     \cU\left[-\frac{t_{k}}{2}, \frac{t_{k}}{2}\right],
\end{equation}
where $\mathbf{t} = \{t_{1} \geq t_{2} \geq \ldots \geq t_m\}$ is a
multiset of non-negative real numbers written in decreasing order.
Thus, up to an overall constant shift, in order to classify all
possible finite generalized uniform sum distributions, it suffices to
classify finite sums of independent central continuous uniform random
variables of the form \eqref{eq:t_uniform_sum.b}.

\begin{figure}[ht]
 \includegraphics[height=3cm]{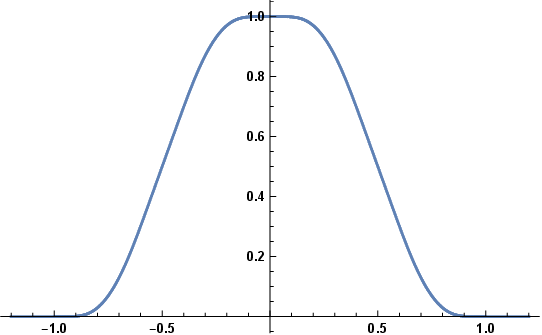}
\hspace{.5in}
 \includegraphics[height=3cm]{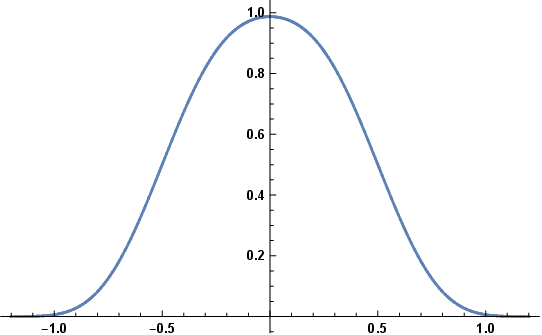}
\caption{Plots of density functions for the distributions
$\cS_{\mathbf{t}}$ with $\mathbf{t}=(1, 1/2, 1/4, 1/8, 1/16, 1/32, 1/64, 1/128, 1/256)$
and $\mathbf{t}=(1, 1/2, 1/3, 1/4, 1/5, 1/6, 1/7, 1/8, 1/9)$.
\label{fig:powersofhalf}}
\end{figure}

\begin{Example}\label{ex:dustpan2n}
Consider the $1/2$-power sequence $\mathbf{t}=(1,1/2,1/4,1/8,\dots )$.
The density function for the distribution $\cS_{\mathbf{t}}$
in \Cref{fig:powersofhalf} has a rather flat top like the sum of two
uniform distributions, in contrast to the harmonic sequence
$\mathbf{t}=(1, 1/2, 1/3, 1/4, 1/5, \ldots)$.
\end{Example}

We will show below that a similar classification holds for the
distributions associated to countable sums of independent continuous
uniform random variables, which are defined provided the expectation
and variance are finite.  Again we have a nice formula for the
cumulants of infinite sums of uniform random variables simply by
letting $m \to \infty$.  Observe that \eqref{eq:Ucont.cums}
is very similar to the definition of the
$p$-norm for a real vector space.

\begin{Definition}
  Let $\mathbf{t} = (t_1, t_2, \ldots)$ be a sequence of
  non-negative real numbers. For $p \in \bR_{\geq 1}$, the
  \textit{$p$-norm} of $\mathbf{t}$ is
  $|\mathbf{t}|_p \coloneqq \left(\sum_{k=1}^\infty t_k^p\right)^{1/p}$.
  We also set $|\mathbf{t}|_\infty \coloneqq \sup_k t_k$.
\end{Definition}

The $p$-norm has many nice properties.  In particular for $d\geq 2$
and $\mathbf{t} = (t_{1},\ldots,t_m)$, we have
\begin{equation}\label{eq:Ucont.cumulants.2}
\kappa_d^{\cS_{\mathbf{t}}} = \frac{B_d}{d} \sum_{k=1}^m (t_k)^d =
\frac{B_d}{d} |\mathbf{t}|_d^d.
\end{equation}
It is well-known (e.g.~\cite[Ex.~7.3, p.58]{MR1483073})
that if $1 \leq p \leq q \leq \infty$, then
$|\mathbf{t}|_p \geq |\mathbf{t}|_q$, and that if $|\mathbf{t}|_p <
\infty$, then $\lim_{p \to \infty} |\mathbf{t}|_p =
|\mathbf{t}|_\infty$. Thus, if $\mathbf{t}$ is
weakly decreasing, $|\mathbf{t}|_\infty = \sup_k t_k = t_1$.

The \textit{sequence space with finite $p$-norm} $\mathbf{\ell}_p
\coloneqq \{\mathbf{t} = (t_1, t_2, \ldots) \in \bR_{\geq 0}^{\bN} : |\mathbf{t}|_p <
\infty\}$ is commonly used in functional analysis and statistics.
Here we define a related concept for analyzing sums of central
continuous uniform random variables.

\begin{Definition}
  The \textit{decreasing sequence space with finite $p$-norm} is \[
\widetilde{\ell}_p \coloneqq \{\mathbf{t} = (t_1, t_2,
\ldots) : t_1 \geq t_2 \geq \cdots \geq 0, |\mathbf{t}|_p <
\infty\}. \] 
\end{Definition}

The elements of $\widetilde{\ell}_p$ may equivalently be thought of as
the set of \textit{countable multisets of non-negative real numbers
with finite $p$-norm}. Any finite multiset of non-negative real
numbers can be considered as an element of $\widetilde{\ell}_p$ with
finite support by sorting the multiset and appending $0$'s.  The
multisets in $\widetilde{\ell}_{p}$ are uniquely determined by their
$p$-norms.  In fact, any subsequence of $p$-norm values injectively
determines the multiset provided the sequence goes to infinity.

\begin{Lemma}\label{lem:norm_injective}
  Let $\mathbf{t}, \mathbf{u} \in \widetilde{\ell}_p$ for some $1 \leq p
\leq \infty$.  Suppose $|\mathbf{t}|_{p_j} = |\mathbf{u}|_{p_j}$ for
some sequence $p_j \to \infty$. Then $\mathbf{t} =
\mathbf{u}$.

  \begin{proof}
    We have
      \[ t_{1} = \sup_k t_k = |\mathbf{t}|_\infty = \lim_{j \to \infty}
          |\mathbf{t}|_{p_{j}} = \lim_{j \to \infty} |\mathbf{u}|_{p_{j}} =
          |\mathbf{u}|_\infty = \sup_k u_k = u_{1}.
      \]
    We may remove the first elements from both $\mathbf{t}$ and
    $\mathbf{u}$ to obtain multisets $(t_{2}, t_{3},\ldots)$ and $(u_{2},
    u_{3},\ldots)$ which are both in $\widetilde{\ell}_p$ and have equal
    $p_{j}$-norms again.  While removing these largest elements alters the
    $p_{j}$-norms, it does so by the same amount for both $\mathbf{t}$ and
    $\mathbf{u}$. Repeating the argument, $t_i = u_i$ for all $i$, so
    $\mathbf{t} = \mathbf{u}$.
  \end{proof}
\end{Lemma}

\begin{Theorem}\label{prop:multiset_finite}
  Finite generalized uniform sum distributions are bijectively
parameterized by \[ \mathbb{R} \times \{\mathbf{t} \in \widetilde{\ell}_2 :
\mathbf{t}\text{ has finite support}\}. \]
  
  \begin{proof} As noted above, every such distribution is defined by
a random variable of the form $c+\cS_{\mathbf{t}}$ for some $c \in
\mathbb{R}$ and $\mathbf{t} = (t_{1},\dots , t_{m}, 0,0,\dots ) \in \widetilde{\ell}_2 $.  To show
uniqueness, suppose $\cS_{\mathbf{t}} = \cS_{\mathbf{u}}$.  By
\eqref{eq:Ucont.cumulants.2}, we know \[ \frac{B_d}{d}
|\mathbf{t}|_{d}^d = \frac{B_d}{d} \sum_{k=1}^m t_k^d =
\kappa_d^{\cS_{\mathbf{t}}} = \kappa_d^{\cS_{\mathbf{u}}} =
\frac{B_d}{d} \sum_{k=1}^m u_k^d = \frac{B_d}{d}
|\mathbf{u}|_{d}^d.  \] Therefore, since the even Bernoulli numbers
are non-zero, we have $|\mathbf{t}|_d = |\mathbf{u}|_d$ for each $d$
even, which is a sequence approaching infinity. Hence, by
\Cref{lem:norm_injective}, $\mathbf{t} = \mathbf{u}$.  \end{proof}
\end{Theorem}

The probability density functions (PDF) for any finite generalized uniform sum
distributions can be determined as a convolution.  We will not need
this formula in the rest of this paper, but we note it here for
completeness. It was used to generate \Cref{fig:powersofhalf}.

\begin{Lemma}\label{lem:density}
 Let $\mathbf{t} = \{t_1\geq  \ldots\geq  t_m > 0\}$. Then
 $\mathrm{PDF}(\mathcal{S}_{\mathbf{t}}; x)$ is given by 
\[
 \frac{1}{2(m-1)! t_1 \cdots t_m} \sum_{\epsilon_1, \ldots, \epsilon_m \in \{\pm 1\}} \epsilon_1 \cdots \epsilon_m \left(x + \frac{\epsilon_1 t_1 + \cdots + \epsilon_m t_m}{2}\right)^{m-1} \mathrm{sgn}\left(x + \frac{\epsilon_1 t_1 + \cdots + \epsilon_m t_m}{2}\right).
\]
\end{Lemma}

\begin{proof} For the case $m=1$, 
\[
\mathrm{PDF}(\mathcal{U}[-t_1/2, t_1/2]; x) = \frac{1}{2t_1} \left(\mathrm{sgn}\left(x+\frac{t_1}{2} \right) - \mathrm{sgn}\left(x-\frac{t_1}{2} \right)\right).
\]
Let $\ast$ denote convolution.
One can check that for all $u>0,$ we have the convolution identity 
\begin{align*}
x^k \mathrm{sgn}(x)\ \ast \ &\frac{1}{2} \left(\mathrm{sgn}(x+u) - \mathrm{sgn}(x-u)\right) \\
&= \frac{1}{k+1} \left((x+u)^{k+1} \mathrm{sgn}(x+u) - (x-u)^{k+1} \mathrm{sgn}(x-u)\right).
\end{align*}
The probability density function of the sum of independent random
variables is the convolution of their density functions. Therefore,
the general case of the lemma now follows by applying the $m=1$ case
and the convolution identity inductively.
\end{proof}

\begin{Remark}\label{rem:IHdensity}
 When $t_1 = \cdots = t_m = 1$, the formula in \Cref{lem:density}  collapses to
\[
\mathrm{PDF}(\mathcal{IH}_m - m/2; x) = \frac{1}{2(m-1)!} \sum_{k=0}^m (-1)^k \binom{m}{k} \left(x + \frac{m-k}{2} - \frac{k}{2}\right)^{m-1} \mathrm{sgn}\left(x + \frac{m-k}{2} - \frac{k}{2}\right).
\]
Hence we recover the density formula for the Irwin--Hall
distributions \cite[p. 296]{JohnsonKotzBalakrishnan}
\[
\mathrm{PDF}(\mathcal{IH}_m; x) = \frac{1}{2(m-1)!} \sum_{k=0}^m
(-1)^k \binom{m}{k} \left(x - k\right)^{m-1} \mathrm{sgn}\left(x -
k\right).
\]
\end{Remark}

\begin{Remark}\label{rem:sum.cdf}
A similar formula for the cumulative distribution function of
$\mathcal{S}_{\mathbf{t}}$ as a sum over the vertices of the hypercube
is given in \cite{MR1538918}. See also \cite[p. 298-300]{JohnsonKotzBalakrishnan}
for relevant discussion.
\end{Remark}


We now turn to infinite sums of independent uniform continuous random
variables.  Our next goal is to generalize \Cref{prop:multiset_finite}
to this setting. To do so, we must first extend the uniform-sum
distributions $\cS_\mathbf{t}$ to countably infinite multisets
$\mathbf{t}$, and discuss the basic properties of these random
variables including existence, characteristic functions, and
cumulants.  Existence depends on the following well-known result,
which often appears in treatments of the law of large numbers.  See,
for example, \cite[Thm.~2.5.3]{MR2722836}.

\begin{Theorem}[Kolmogorov's Two-Series Theorem]\label{thm:kolmogorov}
  Let $\cX_1, \cX_2, \ldots$ be a sequence of independent
  real-valued random variables. Suppose $\bE[\cX_k] = 0$ and
  $\sum_{k=1}^\infty \Var[\cX_k] < \infty$. Then
  $\sum_{k=1}^\infty \cX_k$ converges almost surely.
\end{Theorem}

Almost sure convergence implies convergence in distribution.
Therefore, by Kolmogorov's Two-Series Theorem, we are lead to the
following definition.  

\begin{Definition}\label{def:generalized.uniform.sum}
A \textit{generalized uniform sum distribution} is any distribution
associated to a random variable with finite mean and variance given as a
countable sum of independent continuous uniform random variables.
As in the finite case, such random variables are given by a constant
overall shift plus a \textit{uniform sum random variable}
\[ \cS_{\mathbf{t}} \coloneqq  \cU\left[-\frac{t_1}{2},
\frac{t_1}{2}\right] + \cU\left[-\frac{t_2}{2}, \frac{t_2}{2}\right] +
\cdots \] for some $\mathbf{t}=(t_1, t_2, \ldots) \in
\widetilde{\ell}_2$. Kolmogorov's Theorem applies since
$\Var[\cU[-t/2, t/2]] = \frac{B_2}{2} t^2$ and
$\sum_{k=1}^\infty \Var[\cU[-t_k/2, t_k/2]]
= \frac{B_2}{2}|\mathbf{t}|_2^2 < \infty$.
\end{Definition}

Conversely, Kolmogorov's stronger Three-Series Theorem
\cite[Thm.~2.5.4]{MR2722836} shows that if $\sum_{i=1}^\infty t_i^2 =
\infty$, then $\sum_{i=1}^\infty \cU[-t_i/2, t_i/2]$ diverges with
positive probability, so the assumption $|\mathbf{t}|_2^2 < \infty$ is
essential. In this way we also see that uncountably many non-zero
summands of independent continuous uniform random variables must
diverge.  Thus, we cannot extend \Cref{def:generalized.uniform.sum}
beyond countable sums.

We claim that each uniform sum random variable $\cS_{\mathbf{t}}$ for
$\mathbf{t} \in \widetilde{\ell}_2$ gives rise to a distinct
distribution.  In order to prove the claim, we need to verify the
relationship between the $p$-norms and the cumulants of the infinite
sums is as expected.  To do so, we describe the characteristic and
moment-generating functions of $\cS_\mathbf{t}$.

\begin{Lemma}\label{lem:inf_unif_props}
  Let $\mathbf{t} =(t_1, t_2, \ldots) \in \widetilde{\ell}_{2}$. Then
$\cS_{\mathbf{t}}$ exists, has moments of all orders, and is
determined by its moments. The characteristic function is the entire
function \begin{equation}\label{eq:char.inf.sum}
\phi_{\cS_{\mathbf{t}}}(s) = \prod_{k=1}^\infty \newsinc{st_k/2},
\qquad s \in \bC.  \end{equation} Moreover, $\bE[\cS_{\mathbf{t}}] =
0$, $\Var[\cS_{\mathbf{t}}]< \infty$, and for each $d \in \bZ_{\geq
2}$, \begin{equation}\label{eq:cum.inf.sum}
\kappa_d^{\cS_{\mathbf{t}}} = \frac{B_d}{d} \sum_{k=1}^\infty t_k^d =
\frac{B_d}{d} |\mathbf{t}|_d^d.  \end{equation}
  
  \begin{proof} As mentioned above, the assumption $\mathbf{t} \in
\widetilde{\ell}_{2}$ and \Cref{thm:kolmogorov} together imply
$\cS_{\mathbf{t}}$ exists.  The characteristic function of $\cU[-x,
x]$ is \begin{equation}\label{e:char.uniform}
\phi_{\cU[-x, x]}(s) = \frac{1}{2x} \int_{-x}^x e^{ist}\,dt = \frac{e^{isx}-e^{-isx}}{2isx} =
\frac{\sin(sx)}{sx} \coloneqq \newsinc{sx}, \end{equation} where
$\newsinc{0} \coloneqq 1$.  Consequently, the $n$th partial sum
$\cS_n=\sum_{k=1}^{n}\cU\left[-\frac{t_k}{2}, \frac{t_k}{2}\right]$
has characteristic function $\phi_{\cS_n}(s) = \prod_{k=1}^n
\newsinc{st_k/2}$. Almost sure convergence implies convergence in
distribution, so $\cS_n \Rightarrow \cS_{\mathbf{t}}$.
Thus, by L\'evy's Continuity Theorem, we have for
each $s \in \bR$ that \[ \phi_{\cS_{\mathbf{t}}}(s) =
\prod_{k=1}^\infty \newsinc{st_k/2}. \] By \Cref{lem:inf_unif_cf}
below, the product form for $\phi_{\cS_{\mathbf{t}}}$ is entire and hence
complex analytic on an open ball, so \eqref{eq:cum.inf.sum} follows
from \Cref{rem:analytic_cf}. Likewise, $\cS_{\mathbf{t}}$ has moments
of all orders and $\cS_{\mathbf{t}}$ is determined by its moments.
    
    Since the entire functions $\phi_{\cS_n}(s)$ converge uniformly on
compact subsets of $\bC$ to $\phi_{\cS_{\mathbf{t}}}(s)$, it follows
that the $d^{th}$ moment can be determined by the constant term of the
$d^{th}$ derivative of the characteristic function
\begin{align*} \lim_{n \to \infty} \bE[\cS_n^d] = \lim_{n \to
\infty} i^{-d} \phi_{\cS_n}^{(d)}(0) = i^{-d}
\phi_{\cS_{\mathbf{t}}}^{(d)}(0) =
\bE[\cS_{\mathbf{t}}^d] \end{align*} for all $d\geq 1$.  The moments
of any random variable determine its cumulants and vice
versa. Therefore, the cumulant formula now follows from
\eqref{eq:Ucont.cumulants.2}, including the first two
moments.  \end{proof}
\end{Lemma}


\begin{Lemma}\label{lem:inf_unif_cf}
Let $\mathbf{t}=(t_1, t_2, \ldots) \in \widetilde{\ell}_{2}$.  As a
function of $s$, the infinite product \[ \prod_{i=1}^\infty
\newsinc{st_i/2} \] converges to an entire function in the complex plane.
Moreover, for $|s| < 1/|\mathbf{t}|_2$,
  \[ \left|\prod_{i=1}^\infty \newsinc{st_i/2}\right|
      \leq e. \]

  \begin{proof}
    For each $D>0$, the entire function $\frac{1-\newsinc{z}}{z^2}$ is bounded
    on $|z| < D$ by some constant $C>0$. Thus
      \[ \left|1 - \newsinc{z}\right| \leq C|z|^2 \qquad \text{for }|z|<D. \]
    Consequently, for $|s| < 2D/\sup\{t_i\}$, we have
      \[ \left|1 - \newsinc{st_i/2}\right| < \frac{C}{4} |s|^2 t_i^2. \]
    Hence
      \[ \sum_{i=1}^\infty \left|1 - \newsinc{st_i/2}\right|
          \leq \frac{C}{4} |s|^2 |\mathbf{t}|_2^2 < \infty. \]
    Thus, the sum converges uniformly on compact subsets of $\{|s| < 2D/\sup\{t_i\}\}$.
    Taking $D \to \infty$, the sum converges uniformly on compact subsets of all of $\bC$.
    The result now follows by standard criteria for infinite product convergence such
    as \cite[Thm.~15.6]{rudin_real_complex}.
    
    For the growth rate bound, it is straightforward to check that when $D=1/2$,
    we may use $C=4$. Since $|\mathbf{t}|_2 \geq |\mathbf{t}|_\infty = \sup\{t_i\}$,
    for $|s| < 1/|\mathbf{t}|_2$, we have
    \begin{align*}
      \left|\prod_{i=1}^\infty \newsinc{st_i/2}\right|
        &= \prod_{i=1}^\infty \left|1 - \left(1 - \newsinc{st_i/2}\right)\right|
        \leq \prod_{i=1}^\infty \left(1 + \left|1 - \newsinc{st_i/2}\right|\right) \\
        &\leq \prod_{i=1}^\infty \left(1 + |s|^2 t_i^2\right)
        \leq \prod_{i=1}^\infty \exp\left(|s|^2 t_i^2\right) \\
        &= \exp\left(|s|^2 |\mathbf{t}|_2^2\right)
        \leq \exp\left(1\right) = e.
    \end{align*}
  \end{proof}
\end{Lemma}

\begin{Theorem}\label{prop:multiset_infinite}
Generalized uniform sum distributions are bijectively parameterized by 
$\mathbb{R} \times \widetilde{\ell}_2.$ 
In particular, if
$\mathbf{t}, \mathbf{u} \in \widetilde{\ell}_2$ with $\mathbf{t} \neq
\mathbf{u}$, then $\cS_{\mathbf{t}} \neq \cS_{\mathbf{u}}$.
Furthermore, $\cS_{\mathbf{t}}^* = \cS_{\mathbf{u}}^*$ if and only if
$\mathbf{t}, \mathbf{u}$ differ by a scalar multiple.

  \begin{proof} The first and second claims follow exactly as in
\Cref{prop:multiset_finite} using the cumulant formula in
\Cref{lem:inf_unif_props}.  For the third claim, we can assume
$|\mathbf{t}|_\infty = |\mathbf{u}|_\infty$ by rescaling if necessary
and $\cS_{\mathbf{t}}^* = \cS_{\mathbf{u}}^*$. From
\Cref{lem:inf_unif_props} and the general properties of cumulants, it
follows that for all $d$ even, \[
|\mathbf{t}|_d^d/|\mathbf{t}|_2^{d/2} =
|\mathbf{u}|_d^d/|\mathbf{u}|_2^{d/2}. \] Taking $d$th roots and the
limiting sequence of positive even integers $d$, this
implies \begin{align*}
\frac{|\mathbf{t}|_\infty}{|\mathbf{t}|_2^{1/2}} = \lim_{d \to \infty}
\frac{|\mathbf{t}|_d}{|\mathbf{t}|_2^{1/2}} = \lim_{d \to \infty}
\frac{|\mathbf{u}|_d}{|\mathbf{u}|_2^{1/2}} =
\frac{|\mathbf{u}|_\infty}{|\mathbf{u}|_2^{1/2}}.  \end{align*} Since
$|\mathbf{t}|_\infty = |\mathbf{u}|_\infty$, we have $|\mathbf{t}|_2 =
|\mathbf{u}|_2$, which hence gives $|\mathbf{t}|_d = |\mathbf{u}|_d$
for all $d$ even. Again by \Cref{lem:norm_injective}, we have
$\mathbf{t} = \mathbf{u}$.  \end{proof}
\end{Theorem}

\begin{Example}
  Infinite sums of independent continuous uniform random variables have
appeared elsewhere in the literature, though rarely. For instance, when $\mathbf{t} =
(1, 1/2, 1/4, 1/8, \ldots)$, the cumulative distribution function of
$\cS_{\mathbf{t}} = \sum_{k=1}^\infty \cU[-1/2^k, 1/2^k]$ is the
so-called \textit{Fabius function}, \cite{MR0197656}, which is a known
example of a $C^\infty$-function on an interval which is nowhere
analytic. The characteristic function is nonetheless entire by
\Cref{lem:inf_unif_cf}.
\end{Example}

\begin{Example}
  Another interesting case arises from
  $\mathbf{t} = (1, 1/2, 1/3, 1/4, \ldots)$. Since
  $|\mathbf{t}|_2 = \sum_{k=1}^\infty 1/k^2 < \infty$,
  $\cS_{\mathbf{t}} = \sum_{k=1}^\infty \cU[-1/(2k), 1/(2k)]$
  converges almost surely. For $d\geq 1$,  we have
  \begin{align*}
    \kappa_{2d}
      &= \frac{B_{2d}}{2d} \sum_{k=1}^\infty \frac{1}{k^{2d}}
        = \frac{B_{2d}}{2d} \zeta(2d).
  \end{align*}
  Using the known identity
    \[ \zeta(2d) = (-1)^{d+1} \frac{B_{2d} (2\pi)^{2d}}{2(2d)!}, \]
  it follows that
  \begin{align*}
    \log\phi_{\cS_{\mathbf{t}}}(s)
      &= \sum_{k=1}^\infty \kappa_k \frac{s^k}{k!} \\
      &= -\sum_{d=1}^\infty \frac{\zeta(2d)^2}{d} \left(\frac{s}{2\pi}\right)^{2d},
  \end{align*}
  which is valid in a complex neighborhood of $s=0$.
  This last expression is similar to the left-hand side of the known identity
    \[ \sum_{d=0}^\infty \zeta(2d) s^{2d} = -\frac{\pi s}{2} \cot(\pi s). \]
\end{Example}

\begin{Example}
  Let $\alpha \in \bR_{>0}$ and set $\mathbf{t}^{(N)} = (1/N^\alpha,
  1/N^\alpha, \ldots, 1/N^\alpha, 0, \ldots)$ where there are $N$
  non-zero terms. Then $|\mathbf{t}^{(N)}|_p =
  N^{\frac{1}{p}-\alpha}$. So, for $1 \leq p < \infty$,
    \[ \lim_{N \to \infty} |\mathbf{t}^{(N)}|_p =
        \begin{cases}
          0 & \text{if }p > 1/\alpha \\
          1 & \text{if }p = 1/\alpha \\
          \infty & \text{if }p < 1/\alpha.  \end{cases} \] On the
other hand, for each $k$ we have $\lim_{N \to \infty} t_k^{(N)} = 0$,
independent of $\alpha$. Hence we have a large family of sequences
which each converges pointwise to $(0, 0, \ldots)$, but which have
different limiting $p$-norms. In particular, when $\alpha = 1/2$ we
have $\lim_{N \to \infty} |\mathbf{t}^{(N)}|_2 = 1 \neq 0 = |(0, 0,
\ldots)|_2$, so the limit of the $2$-norms is not the $2$-norm of the
limit. The interplay between convergence in $\widetilde{\ell}_2$ and
convergence of generalized uniform sum distributions is consequently
somewhat subtle, which we treat in the next subsection.
\end{Example}

\subsection{Pointwise convergence and convergence in even norms}

The decreasing sequence space $\widetilde{\ell}_2$ has a natural
notion of pointwise convergence.  In this subsection, we relate pointwise convergence to
convergence of $p$-norms for all positive even $p\geq 4$, assuming the
$2$-norms are bounded.

\begin{Lemma}\label{lem:even_to_first}
Fix $M \in \mathbb{R}$. Let $\mathbf{t}^{(N)} \in
\widetilde{\ell}_2$ be a countable sequence of sequences such that
$|\mathbf{t}^{(N)}|_2^2 \leq M$ for each $N$ and  \[
\lim_{N \to \infty} |\mathbf{t}^{(N)}|_{2d} = \tau_{2d} \] exists for
all $d \in \bZ_{\geq 2}$. Then

\begin{enumerate}[(i)]
\item $\lim_{d \to \infty} \tau_{2d}$ exists,
\item $\lim_{N \to \infty} t_1^{(N)}$ exists, 
\item $\lim_{d \to \infty} \tau_{2d} = \lim_{N \to \infty} t_1^{(N)} =
\lim_{N \to \infty} |\mathbf{t}^{(N)}|_\infty$, and

\item $\mathbf{t}^{(N)}$ converges
  pointwise to $\mathbf{t}=(t_{1},t_{2},\dots ) \in
  \widetilde{\ell}_2$ where $t_i = \lim_{N \to \infty} t_i^{(N)}$.
\end{enumerate} 

  \begin{proof} For (i), if $d \leq e \leq \infty$ then
$|\mathbf{t}^{(N)}|_{2d} \geq |\mathbf{t}^{(N)}|_{2e}$ by properties
of the $p$-norm.  Therefore, $\tau_{2d} \geq \tau_{2e} \geq 0$ and
$\lim_{d \to \infty} \tau_{2d}$ exists.

    For (ii), observe that since $\mathbf{t}^{(N)}$ is a decreasing sequence in
    $\widetilde{\ell}_2$, we know  $|\mathbf{t}^{(N)}|_\infty=t_1^{(N)} \geq t_i^{(N)}$ for all
    $i$.  Therefore, for
all $d \in \bZ_{\geq 1}$, we have 
    \begin{align*}
      |\mathbf{t}^{(N)}|_{2d}^{2d}
        &= \sum_i (t_i^{(N)})^{2d} \\
        &\leq \sum_i (t_1^{(N)})^{2(d-1)} (t_i^{(N)})^2 \\
        &\leq (t_1^{(N)})^{2(d-1)} \cdot M.
    \end{align*}
Combining this with the fact that $t_{1}^{(N)} \leq
|\mathbf{t}^{(N)}|_{2d}$ by definition of the $p$-norm,
one has \begin{equation}\label{eq:2d_infty_bounds} t_{1}^{(N)} \leq
|\mathbf{t}^{(N)}|_{2d} \leq (t_{1}^{(N)})^{1-\frac{1}{d}} \cdot
M^{\frac{1}{2d}}.  \end{equation} Taking $N \to \infty$ in
\eqref{eq:2d_infty_bounds}
gives \begin{equation}\label{eq:2d_infty_bounds.2} \limsup_{N \to
\infty} (t_{1}^{(N)}) \leq \tau_{2d} \leq \liminf_{N \to \infty}
(t_{1}^{(N)})^{1 - \frac{1}{d}} \cdot
M^{\frac{1}{2d}}.  \end{equation} Taking $d \to \infty$ in
\eqref{eq:2d_infty_bounds.2} gives 
\begin{equation}\label{eq:pinched}
 \limsup_{N \to \infty}
(t_{1}^{(N)}) \leq \lim_{d \to \infty}\tau_{2d} \leq \liminf_{N \to
\infty} (t_{1}^{(N)}), 
\end{equation}
so $\lim_{N \to \infty} t_1^{(N)} = \lim_{d \to \infty} \tau_{2d}$
which implies the limit exists by (i).  Part (iii) also follows from
\eqref{eq:pinched} and the fact that
$|\mathbf{t}^{(N)}|_\infty=t_1^{(N)}$.

Part (iv) follows by an inductive argument.  By (ii), $t_{1}=\lim_{N
\to \infty} t_1^{(N)} $ exists. Define another sequence of sequences
$\mathbf{u}^{(N)} \coloneqq \{t_2^{(N)} \geq t_3^{(N)} \geq \cdots\}$,
so that $|\mathbf{u}^{(N)}|_{2}^{2} = |\mathbf{t}^{(N)}|_{2}^{2} -
(t_1^{(N)})^{2} \leq M$ and
\begin{align*}
|\mathbf{u}^{(N)}|_{2d}^{2d} = |\mathbf{t}^{(N)}|_{2d}^{2d} -
(t_1^{(N)})^{2d} \qquad \Rightarrow \qquad \lim_{N \to \infty}
|\mathbf{u}^{(N)}|_{2d} = \left(\tau_{2d}^{2d} -
t_1^{2d}\right)^{\frac{1}{2d}}\text{ exists}
\end{align*}
by the hypotheses on $\mathbf{t^{(N)}}$.  
By (iii) applied to $\mathbf{u}^{(N)}$, $t_2
\coloneqq \lim_{N \to \infty} u_1^{(N)} = \lim_{N \to \infty}
t_2^{(N)}$ exists. Repeating the argument, $\mathbf{t}^{(N)}$
converges pointwise to  $(t_{1},t_{2},\dots)$.  
\end{proof}
\end{Lemma}

\begin{Lemma}\label{lem:pointwise_to_even}
  Suppose $\mathbf{t}^{(N)} \in \widetilde{\ell}_2$ with
$|\mathbf{t}^{(N)}|_2^2 \leq M$ converges pointwise to
$\mathbf{t}\in \widetilde{\ell}_2$. Then
$|\mathbf{t}|_2^2 \leq M$ and for all $d \geq 2$, 
\[
|\mathbf{t}|_{2d} = \lim_{N \to \infty} |\mathbf{t}^{(N)}|_{2d}.
\]


  \begin{proof}
By Fatou's Lemma applied to the counting measure on
$\bZ_{\geq 1}$, \[ |\mathbf{t}|_2^2 \leq \liminf_{N \to \infty}
|\mathbf{t}^{(N)}|_2^2 \leq M. \]
Fix $d \geq 2$. For each $N$, we have $t_1^{(N)} \geq
    t_2^{(N)} \geq \cdots \geq t_i^{(N)} \geq \cdots$. Thus
\[
      M \geq (t_1^{(N)})^2 + \cdots + (t_i^{(N)})^2   \geq  i(t_i^{(N)})^2,
\]
which implies 
\[
	       (t_i^{(N)})^2 \leq \frac{M}{i} \
        \Rightarrow \ (t_i^{(N)})^{2d} \leq \left(\frac{M}{i}\right)^d.
\]
    Since $\sum_{i=1}^\infty \frac{1}{i^d}$ converges for $d \geq 2$,
    the sequence $(t_i^{(N)})^{2d}$
    is dominated by the integrable function
    $\left(\frac{M}{i}\right)^d$ over the positive integers.  By Lebesgue's Dominated Convergence 
    Theorem, since $\lim_{N \to \infty} (t_i^{(N)})^{2d} = t_i^{2d}$,
    we have
      \[ \lim_{N \to \infty} |\mathbf{t} - \mathbf{t}^{(N)}|_{2d} = 0
         \qquad \Rightarrow \qquad
         \lim_{N \to \infty} |\mathbf{t}^{(N)}|_{2d} = |\mathbf{t}|_{2d}. \]
  \end{proof}
\end{Lemma}

\begin{Corollary}\label{cor:pointwise.convergence.equiv}
  Suppose $\mathbf{t}^{(N)} \in \widetilde{\ell}_2$ with
$|\mathbf{t}^{(N)}|_2^2 \leq M$.    Then
$\mathbf{t}^{(N)}$ converges pointwise to $\mathbf{t}$ if and only if
$|\mathbf{t}|_{2d} = \lim_{N \to \infty} |\mathbf{t}^{(N)}|_{2d}$  for
all $ d \geq 2.$

\begin{proof}
The proof follows directly from \Cref{lem:even_to_first} and
\Cref{lem:pointwise_to_even}.
\end{proof}
\end{Corollary}

Observe that \Cref{cor:pointwise.convergence.equiv} says nothing about
the $2$-norm of the sequences.  It is possible for $\mathbf{t}^{(N)}
\to \mathbf{t}$ pointwise, even if $ |\mathbf{t}|_2^2 \neq \lim_{N \to
\infty}|\mathbf{t}^{(N)}|_2^2$, as the next example and lemma
illustrate.

\begin{Example}
  In the Irwin--Hall case, we have $\cIH_N = \cS_{\mathbf{t}^{(N)}} + N/2$
  where
    \[ \mathbf{t}^{(N)} = (\underbrace{1, \ldots, 1}_{N\text{ copies}}, 0, \ldots). \]
  Since $|\mathbf{t}^{(N)}|_2^2 = N$, after standardizing,
  $\cIH_N^* = \cS_{\widehat{\mathbf{t}^{(N)}}}$ where
    \[ \widehat{\mathbf{t}^{(N)}}
       = (\underbrace{\sqrt{12/N}, \ldots, \sqrt{12/N}}_{N\text{ copies}}, 0, \ldots), \]
  which converges pointwise to $\mathbf{t} = (0, 0, \ldots)$.
  Nonetheless,
  $|\mathbf{t}|_2^2 = 0 < 12 = |\widehat{\mathbf{t}^{(N)}}|_2^2$
  and $\cIH_N^* \Rightarrow \cN(0, 1)$.
\end{Example}

\begin{Lemma}\label{cor:weirdness.at.2norm}
For every $\mathbf{t}=(t_{1},t_{2},\dots ) \in \widetilde{\ell}_2$ and
every $M\geq |\mathbf{t}|_2^2$, there exists a sequence $\mathbf{t}^{(N)}$
of finitely supported decreasing sequences such that
$|\mathbf{t}^{(N)}|_2^2=M$ and $\mathbf{t}^{(N)} \to \mathbf{t}$
pointwise.

\begin{proof}
Define a sequence of sequences $\mathbf{t}^{(N)} \in
\widetilde{\ell}_2$ with $|\mathbf{t}^{(N)}|_2^2 = M$ as follows.  Let
\[
\epsilon_N \coloneqq \sqrt{M - \sum_{i=1}^N t_{i}^2}.
\]
For each $N\geq 1$, choose $m_N \in \bZ_{\geq 1}$ large enough so that
$\epsilon_N/m_N \leq \frac{1}{N}$.  Set
\[ \mathbf{t}^{(N)} = (t_{1}, t_{2}, \ldots, t_{N},
\underbrace{\epsilon_N/m_N, \ldots, \epsilon_N/m_N}_{m_N^2\text{
copies}}, 0, 0, \dots ).
\]
As claimed, $\mathbf{t}^{(N)} \to \mathbf{t}$ pointwise and \[
|\mathbf{t}^{(N)}|_2^2 = \sum_{i=1}^N t_i^2 + m_N^2 \cdot
\left(\frac{\epsilon_N}{m_N}\right)^2 = M. \]
\end{proof}
\end{Lemma}

\begin{Example}\label{ex:weirdness.at.2norm}
Consider again $\mathbf{t} = (1,1/2, 1/3, \ldots)$ so $|\mathbf{t}|_2^2 =
\sum_{i=1}^\infty(\frac{1}{i})^2 = \pi^{2}/6 \approx 1.6449$.  Let
$\epsilon_N \coloneqq \sqrt{2 - \sum_{i=1}^N (\frac{1}{i})^2}$. For
each $N\geq 1$, set
\[ \mathbf{t}^{(N)} = (1, 1/2, \ldots, 1/N,
\underbrace{\epsilon_N/N, \ldots, \epsilon_N/N}_{N^2\text{
copies}}, 0, 0, \dots ).
\]
Clearly $\mathbf{t}^{(N)} \to \mathbf{t}$ pointwise and \[
|\mathbf{t}^{(N)}|_2^2 = \sum_{i=1}^N t_i^2 + N^2 \cdot
\left(\frac{\epsilon_N}{N}\right)^2 = 2. \] However, $
|\mathbf{t}|_2^2 = \pi^{2}/6 \neq 2 = \lim_{N \to \infty}|\mathbf{t}^{(N)}|_2^2$.  
\end{Example}

\begin{Lemma}\label{cor:pointwise.convergence.equiv.converg.dis}
  Suppose $\mathbf{t}^{(N)} \in \widetilde{\ell}_2$ converges
  pointwise to $\mathbf{t} \in \widetilde{\ell}_2$ with
  $|\mathbf{t}^{(N)}|_2^2 \to \tau_2 < \infty$. Then
    \[ \cS_{\mathbf{t}^{(N)}} \Rightarrow \cS_{\mathbf{t}}
        + \cN(0, \sigma^2) \]
  where $\sigma = \sqrt{(\tau_2 - |\mathbf{t}|_2^2)/12}$ and
  the sum is independent.

  \begin{proof}
    By \Cref{lem:pointwise_to_even},
    $\lim_{N \to \infty} |\mathbf{t}^{(N)}|_{2d}
    = |\mathbf{t}|_{2d}$ for all $d \in \bZ_{\geq 2}$, so for all $d \geq 3$,
      \[ \kappa_d^{\cS_{\mathbf{t}^{(N)}}}
          \to \kappa_d^{\cS_{\mathbf{t}}}
          = \kappa_d^{\cS_{\mathbf{t}} + \cN(0, \sigma^2)} \]
    since $\kappa_d^{\cN(0, \sigma^2)} = 0$. As
    for $d=2$,
      \[ \kappa_2^{\cS_{\mathbf{t}^{(N)}}} \to \frac{\tau_2}{12}
         = \frac{|\mathbf{t}|_2^2}{12} + \sigma^2
         = \kappa_2^{\cS_{\mathbf{t}} + \cN(0, \sigma^2)}. \]
   The result follows by the Method of Moments/Cumulants.
\end{proof}
\end{Lemma}

In light of \Cref{cor:pointwise.convergence.equiv.converg.dis},
pointwise convergence in $\widetilde{\ell}_2$ leads to us to study an
additional family of sums of random variables.  Note, the sum of two
generalized uniform sum random variables is another generalized
uniform sum of random variables.  Also, the sum of two normal
distributions is normal, so we have reached a natural limit to the
generalizations.

\begin{Definition}\label{def:uniform+nornmal}
A \textit{DUSTPAN distribution} is a
\textit{\underline{d}istribution
associated to a \underline{u}niform
\underline{s}um for $\mathbf{\underline{t}}$ \underline{p}lus
\underline{a}n independent \underline{n}ormal distribution}
$\cS_{\mathbf{t}}+\cN(0, \sigma^{2})$, assuming the two random
variables are independent, $\mathbf{t} \in\widetilde{\ell}_2$, and
$\sigma \in \mathbb{R}_{\geq 0}$.
\end{Definition}

\begin{figure}[ht]
 \includegraphics[height=6cm]{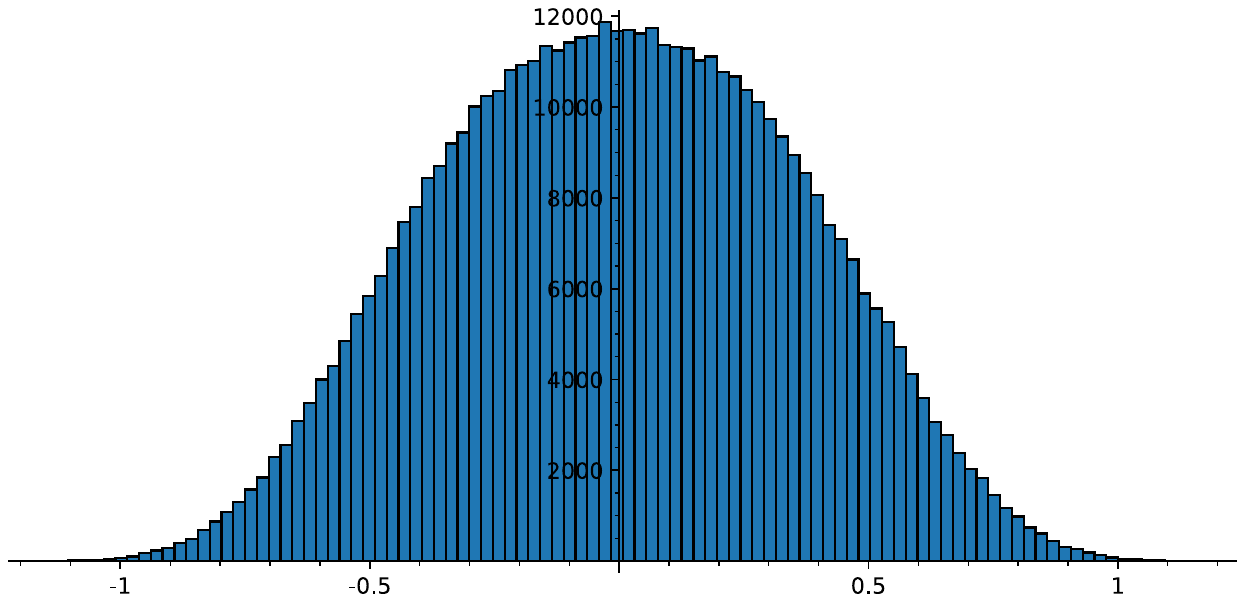}
\hspace{.1in}
 \includegraphics[height=6cm]{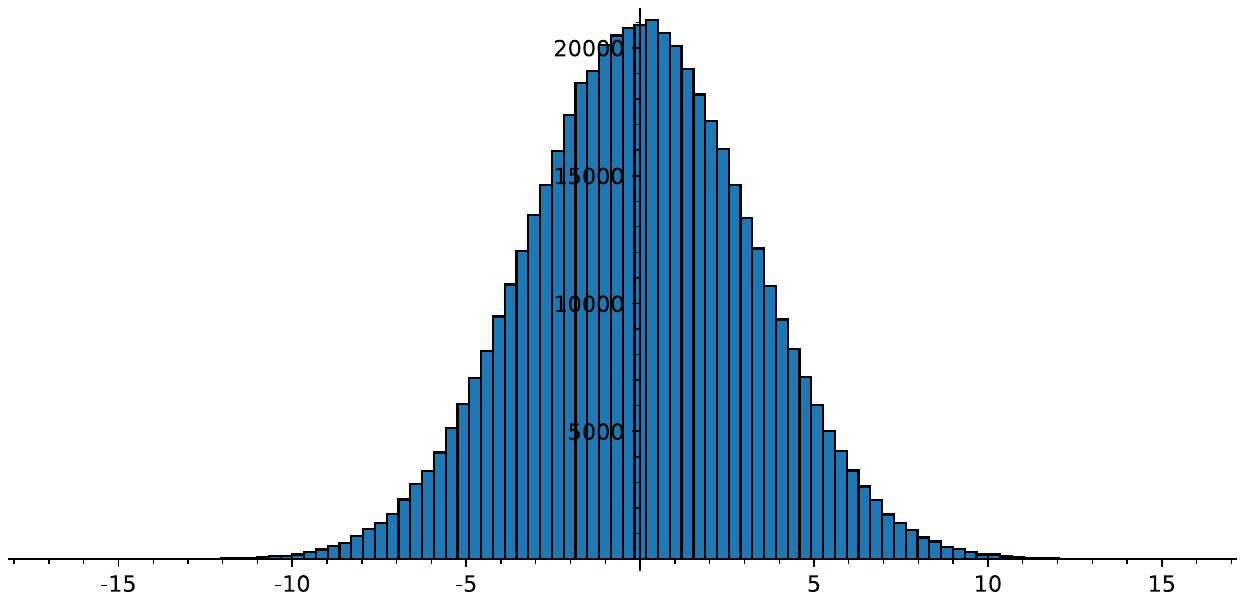}
\hspace{.1in}
 \includegraphics[height=6cm]{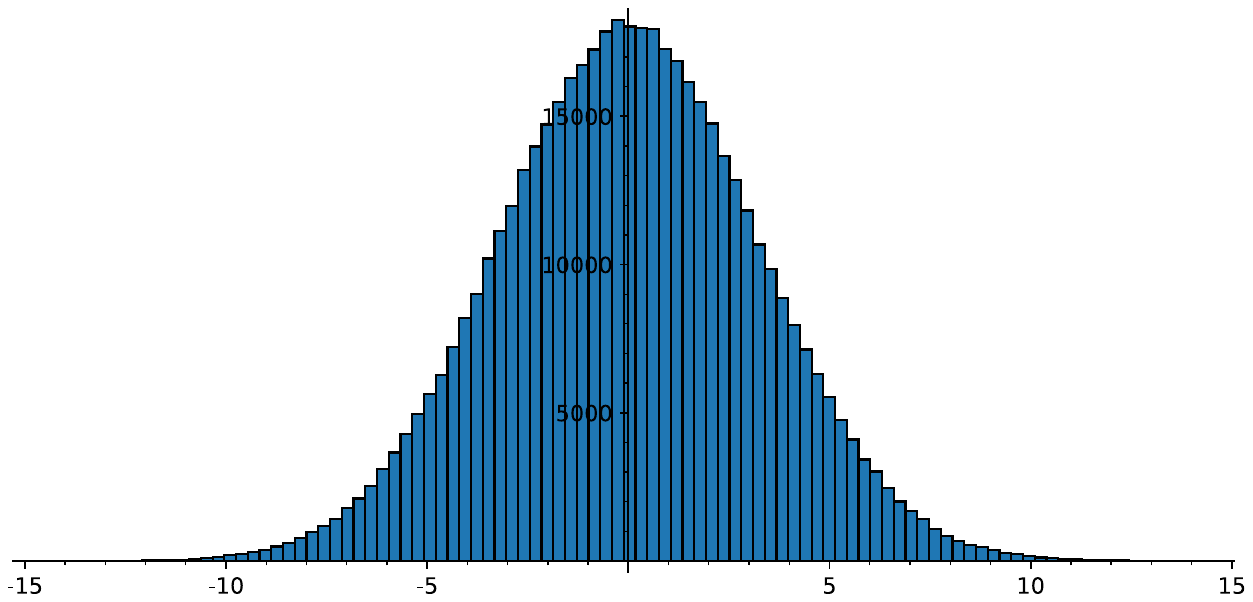}
    \caption{Histograms of the distributions $\cS_{\mathbf{t}}$,
    $\cN(0,\sigma )$, and $\cS_{\mathbf{t}}+\cN(0,\sigma )$
with $\mathbf{t}=(1,\frac{1}{2}, \frac{1}{3}, \frac{1}{4},
\frac{1}{5}, \frac{1}{6}, \frac{1}{7}, \frac{1}{8})$ and $\sigma \approx 3.22$.}
    \label{fig:1overn}
\end{figure}

\begin{Example}\label{ex:dustpan1overn}
Consider the $1/n$-sequence $\mathbf{t}=(1,1/2,1/3,\dots )$ again.  Let
$\sigma = \sqrt{12-\pi^{2}/6}$.  The distribution $\cS_{\mathbf{t}}$
has a small variance compared to $\cN(0, \sigma^{2})$, so
$\cS_{\mathbf{t}}+\cN(0, \sigma^{2})$ looks  like a fat  normal
distribution.  See the approximation in \Cref{fig:1overn}.
\end{Example}

\subsection{The metric space of DUSTPAN distributions}\label{sub:dustpan}

Recall the \textit{metric space of DUSTPAN distributions},
  \[ \bfM_{\DUSTPAN}
      \coloneqq \{\cS_\mathbf{t}
      + \cN(0, \sigma) : |\mathbf{t}|_2^2/12 + \sigma^2 = 1\}, \]
along with the \textit{DUSTPAN parameter space}
  \[ \bfP_{\DUSTPAN}
      \coloneqq \left\{\mathbf{t} \in \widetilde{\ell}_2
        : |\mathbf{t}|_2^{2} \leq 12\right\}. \]
We will show below that $\bfP_{\DUSTPAN}$ and
$\bfM_{\DUSTPAN}$ are homeomorphic closed sets in their respective
topologies of pointwise convergence and convergence in distribution,
thus completing the task of completely characterizing all possible
limit laws of standardized general uniform sum distributions.

From \Cref{def:uniform+nornmal}, it follows that the characteristic
functions of DUSTPAN distributions have nice properties.  Recall that
a \textit{normal family} of holomorphic functions in some open set $U
\subset \bC$ is one where every infinite sequence has a subsequence
which converges uniformly on compact subsets of $U$.

\begin{Lemma}\label{lem:normal_family}
  The set of characteristic functions $\{\phi_{\cS}(s) : \cS \in
\bfM_{\DUSTPAN}\}$ is a normal family of entire functions.
  
  \begin{proof} Let $\cS = \cS_{\mathbf{t}} + \cN(0, \sigma) \in
\bfM_{\DUSTPAN}$. By definition, the characteristic function of a DUSTPAN
distribution is the product of the corresponding characteristic
functions for the normal and generalized uniform sum distributions,
\[ \phi_{\cS}(s) = \exp(-\sigma^2/2)
\prod_{i=1}^\infty \newsinc{st_k}/2. \] By the growth bound in
\Cref{lem:inf_unif_cf}, for $|s| < \frac{1}{12}$, we have \[
|\exp(-\sigma^2/2) \prod_{i=1}^\infty \newsinc{st_k}/2| \leq
\exp(1). \] Thus $\{\phi_{\cS}(s) : \cS \in \bfM_{\DUSTPAN}\}$ is a family
of bounded analytic functions on $|s| < \frac{1}{12}$. By Montel's
Theorem, it is a normal family in that domain.  The bound in
\Cref{lem:inf_unif_cf} may be extended to any bounded domain using the
same argument, so it is in fact a normal family of entire
functions.  \end{proof}
\end{Lemma}

\begin{Lemma}[Converse of Frech\'et--Shohat for DUSTPAN's]\label{prop:MoM_converse}
  Suppose a sequence of DUSTPAN distributions $\cX_N \coloneqq
\cS_{\mathbf{t}^{(N)}} + \cN(0, \sigma^{(N)}) \in \bfM_{\DUSTPAN}$
converges in distribution to some $\cX$. Then $\bE[\cX^d] < \infty$
exists for all $d \in \bZ_{\geq 1}$, $\cX$ is determined by its
moments, and $\lim_{N \to \infty} \bE[\cX_N^d] =
\bE[\cX^d]$.
  
  \begin{proof} By L\'evy's Continuity Theorem, $\phi_{\cX_N}(s) \to
\phi_{\cX}(s)$ for all $s \in \bR$.  By \Cref{lem:normal_family}, we
may replace $\cX_N$ if necessary with a subsequence for which
$\phi_{\cX_N}(s)$ converges uniformly on compact subsets so that we
can assume $\phi_{\cX}(s)$ is entire.  Therefore, the moment
generating function of $\cX$ has positive radius of convergence,
moments of all order exist, $\cX$ is determined by its moments, and
the limit of the moments is the moment of the limit.  \end{proof}
\end{Lemma}

We may now restate and prove \Cref{thm:Pinfty_Minfty} from the introduction.

\begin{repTheorem}{thm:Pinfty_Minfty}
  The map $\Phi \colon \bfP_{\DUSTPAN} \to \bfM_{\DUSTPAN}$ given by $\mathbf{t}
\mapsto \cS_\mathbf{t} + \cN(0, \sigma)$ where $\sigma \coloneqq
\sqrt{1 - |\mathbf{t}|_2^2/12}$ is a homeomorphism between
sequentially compact spaces.
  
  \begin{proof} The parameter space $\bfP_{\DUSTPAN}$ is closed under
pointwise convergence by \Cref{lem:pointwise_to_even}. Moreover, it
is sequentially compact under pointwise convergence, either by Tychonoff's Theorem
applied to $[0, \sqrt{12}]^\bN$ or by a simple diagonalization
argument. Since
$\bfP_{\DUSTPAN}$ and $\bfM_{\DUSTPAN}$ are metrizable and $\Phi$ is a bijection
by \Cref{prop:multiset_infinite}, we need only show that
\begin{equation*}
  \mathbf{t}^{(N)} \to \mathbf{t} \text{ in }\cP\text{ pointwise }
    \qquad\Leftrightarrow\qquad
    \cS_{\mathbf{t}^{(N)}} + \cN(0, \sigma^{(N)})
      \Rightarrow \cS_{\mathbf{t}} + \cN(0, \sigma).
\end{equation*}
The forwards direction follows from \Cref{lem:pointwise_to_even}
and the Method of Moments/Cumulants exactly as in
the proof of \Cref{cor:pointwise.convergence.equiv.converg.dis}.
The backwards direction follows from \Cref{prop:MoM_converse}
and \Cref{lem:even_to_first}.
\end{proof}
\end{repTheorem}


\begin{Corollary}\label{cor:cM_compact}
The metric space of DUSTPAN distributions $\bfM_{\DUSTPAN}$ is compact,
hence it is closed and bounded in the space of distributions under the
L\'evy metric.

\begin{proof}
$\bfP_{\DUSTPAN}$ is a compact subset of $\widetilde{\ell}_2$ under
pointwise convergence, so
$\bfM_{\DUSTPAN}$ is compact under the L\'evy metric as well by \Cref{thm:Pinfty_Minfty}.
\end{proof}
\end{Corollary}

\begin{Corollary}\label{cor:Minfty_closure}
  The closure of the metric space  $\{\cS_{\mathbf{t}}: \mathbf{t} \in
\widetilde{\ell}_2, |\mathbf{t}|_2^2 = 12, \mathbf{t}\text{ is
finite}\}$ in the L\'evy metric is $\bfM_{\DUSTPAN}$.
  
  \begin{proof}
 Since $\{\cS_{\mathbf{t}}: \mathbf{t} \in \widetilde{\ell}_2,
|\mathbf{t}|_2^2 = 12, \mathbf{t}\text{ is finite}\} \subset
\bfM_{\DUSTPAN}$ by definition and $\bfM_{\DUSTPAN}$ is closed by
\Cref{cor:cM_compact}, we know
\[
\overline{\{\cS_{\mathbf{t}}: \mathbf{t} \in \widetilde{\ell}_2,
|\mathbf{t}|_2^2 = 12, \mathbf{t}\text{ is finite}\}} \subset
\overline{\bfM_{\DUSTPAN}} = \bfM_{\DUSTPAN}.
\]
For the other inclusion, we just need to show each $\mathbf{t} \in
\widetilde{\ell}_2$ with $|\mathbf{t}|_2^2 \leq 12$ is the pointwise
limit of a sequence $\mathbf{t}^{(N)} \in \widetilde{\ell}_2$ with
$|\mathbf{t}^{(N)}|_2^2 = 12$ and $\mathbf{t}$ finite by \Cref{thm:Pinfty_Minfty}.  As noted
above, this follows from \Cref{cor:weirdness.at.2norm}.  \end{proof}
\end{Corollary}

\subsection{The metric space of distance distributions}\label{sub:dist}

For convenience, we recall some of the definitions and notation from
the introduction.  For each $\mathbf{t} \in \widetilde{\ell}_2$ with
$|\mathbf{t}|_2>0$, let
\[
 \widehat{\mathbf{t}}
\coloneqq \frac{\sqrt{12} \cdot \mathbf{t}}{|\mathbf{t}|_2}
\]
be the rescaled sequence such that $|
\widehat{\mathbf{t}}|_{2}^{2}=12$ and
$\cS_{\widehat{\mathbf{t}}}=\cS_{\mathbf{t}}^{*}$.  By definition of
the hat-operation,
 $\widehat{\mathbf{t}} \in \bfP_{\DUSTPAN}$ and
$\Phi(\widehat{\mathbf{t}})=\cS_{\widehat{\mathbf{t}}}+\cN(0,0)=\cS_{\widehat{\mathbf{t}}}=\cS_{\widehat{\mathbf{t}}}^{*}$.
The \textit{distance multiset} of $\mathbf{t}=\{t_1 \geq t_2 \geq \cdots
\geq t_m\}$ is the multiset \[ \Delta\mathbf{t} \coloneqq \{t_i - t_j :
1 \leq i < j \leq m\}, \] and the \textit{metric space of
distance distributions } is 
\begin{equation}\label{eq:M.dist.alt}
\bfM_{\DIST} = \left\{\cS_{\widehat{\Delta\mathbf{t}}} : \mathbf{t} =
\{1=t_1 \geq \cdots \geq t_m=0\} \right\}  .
\end{equation}
Thus, the \textit{parameter space of distance multisets}, mentioned in
\Cref{sec:intro},  is defined as 
\begin{equation}\label{eq:P.dist}
 \bfP_{\DIST} \coloneqq 
        \left\{\widehat{\Delta\mathbf{t}}:  \mathbf{t} = \{1=t_1 \geq \cdots \geq
t_m=0\}  \right\}  
\end{equation}
By padding with $0$'s, consider $\bfP_{\DIST} \subset \bfP_{\DUSTPAN} \subset
\widetilde{\ell}_2  $ as a sequence space with the topology of
pointwise convergence.

\begin{Lemma}\label{lem:Pdist.closure}
The closure of $\bfP_{\DIST}$ is $\bfP_{\DIST} \sqcup \{\mathbf{0}\}$.
\begin{proof}
Let $\mathbf{d}^{(N)} \in \bfP_{\DIST}$ be a sequence converging
pointwise to $\mathbf{d}$.  By \Cref{thm:Pinfty_Minfty}, we can assume
$\mathbf{d} \in \bfP_{\DUSTPAN}$.  By definition, each
$\mathbf{d}^{(N)} = \widehat{\Delta\mathbf{t}^{(N)}}$ for some finite
sequence of real numbers $\mathbf{t}^{(N)} = \{1= t_1^{(N)} \geq
\cdots \geq t_{m^{(N)}}^{(N)} =0\}$.

 Suppose $\limsup_{N \to \infty} m^{(N)} < \infty$. We may pass to a
subsequence for which $m^{(N)} = m$ is constant. We may pass to a
further subsequence for which $\mathbf{t}^{(N)} \in [0, 1]^{m}$
converges pointwise to some $\mathbf{t}=\{1= t_1 \geq \cdots
\geq t_{m} =0\} \in [0, 1]^{m}$ and where
$|\mathbf{t}^{(N)}|_2$ converges.  Clearly the distance multiset
operator $\Delta \colon [0, 1]^m \to [0, 1]^{\binom{m}{2}}$ is
continuous, so $\Delta\mathbf{t}^{(N)} \to \Delta\mathbf{t}$, and
moreover $\mathbf{d}^{(N)}=\widehat{\Delta\mathbf{t}^{(N)}} \to
\widehat{\Delta\mathbf{t}}$, so
$\widehat{\Delta\mathbf{t}} = \mathbf{d}$ which implies $\mathbf{d}
\in \bfP_{\DIST}$.

    Now suppose $\limsup_{N \to \infty} m^{(N)} = \infty$. Again, we
may pass to a subsequence if necessary so we may assume $m^{(N)} \to
\infty$.  Since $t_{1}^{(N)}=1$  and $t_{m^{(N)}}^{(N)}=0$ for each $N$,  we
have 
\begin{align*} |\Delta\mathbf{t}^{(N)}|_2^2
        &= \sum_{1 \leq i < j \leq m^{(N)}} (t_i^{(N)} - t_j^{(N)})^2 \\
        &\geq \sum_{1 < \ell < m^{(N)}} \left[(1 - t_\ell^{(N)})^2 + (t_\ell^{(N)} - 0)^2\right] \\
        &\geq \sum_{1 < \ell < m^{(N)}} \frac{1}{2} = \frac{m^{(N)}}{2} - 1.
    \end{align*}
Therefore, $\lim_{N \to \infty} \frac{\sqrt{12}}{|\Delta\mathbf{t}^{(N)}|_2} \to 0$,
so pointwise $\widehat{\Delta\mathbf{t}^{(N)}} \to \mathbf{0}$.  
\end{proof}

\end{Lemma}

\begin{Corollary}\label{cor:MDelta_closure}
Any pointwise convergent sequence
$\widehat{\Delta\mathbf{t}^{(N)}}$ with $\mathbf{t}^{(N)} = \{1= t_1^{(N)} \geq \cdots \geq
t_{m^{(N)}}^{(N)} =0\}$ converges to $\bf{0}$ if and only if $m^{(N)}
\to \infty$.

\end{Corollary}

\begin{Theorem}\label{thm:PDelta_MDelta}
The map $\Phi_{\DIST} \colon \overline{\bfP_{\DIST}} \to
\overline{\bfM_{\DIST}}= \bfM_{\DIST} \sqcup \{\cN(0, 1)\}$ given by
$\mathbf{d} \mapsto \cS_{\mathbf{d}}$ and $\mathbf{0} \mapsto \cN(0,
1)$ is a homeomorphism between (sequentially) compact spaces.
  
  \begin{proof} First note that $\bfP_{\DIST} \subset \bfP_{\DUSTPAN}$
and $\bfM_{\DIST} \subset \bfM_{\DUSTPAN}$ by construction, so
$\Phi_{\DIST}$ is the restriction of $\Phi: \bfP_{\DUSTPAN}
\longrightarrow \bfM_{\DUSTPAN}$.  Therefore, by
\Cref{thm:Pinfty_Minfty}, $\Phi_{\DIST}$ is a homeomorphism.  Since
closed subsets of compact spaces are compact, we know
$\overline{\bfP_{\DIST}}$ is compact by \Cref{lem:Pdist.closure}.
Furthermore, $\Phi (\overline{\bfP_{\DIST}})=\bfM_{\DIST} \sqcup \{\cN(0, 1)\}$ is closed and
compact. 
  \end{proof}
\end{Theorem}

\section{Metric spaces related to $\SSYT_{\leq m}(\lambda)$
distributions}\label{sec:ssyt}


We next consider the family of generating functions for semistandard
tableaux given by the principal specialization of Schur polynomials,
or equivalently the $\rank$ statistic on $\SSYT_{\leq m}(\lambda)$, as
described in \Cref{ssec:intro:SSYT} and \Cref{ssec:back:ssyt}.  An
interesting special case is given by MacMahon's formula for the
$\size$ statistic on the set $\PP(a \times b \times c)$ of plane
partitions inside an $(a \times b \times c)$ box, given in
\eqref{eq:ppformula}. In particular, we will prove \Cref{thm:SSYT_Delta} and
\Cref{thm:PP_median.intro}.  We provide a wide variety of limit law
classification results for these statistics in various regimes.  The
subsections are divided into four natural special cases: $n/m \to 0$,
$n/m \to \infty $, cases based on the number of distinct parts of
$\lambda$, and plane partitions.  See \Cref{sum:SSYT} for a summary.

\subsection{Limit laws with $|\lambda |/m \to 0$ and uniform sums}\label{ssec:ssyt_limits.nsmall}
We begin classifying the limit
laws for semistandard Young tableaux.  Throughout this section, we
tacitly assume $\ell(\lambda) \leq m$, so $\SSYT_{\leq m}(\lambda) \neq
\varnothing$.  Furthermore, if $\lambda_m>0$, the first $\lambda_m$
columns of $T \in \SSYT_{\leq m}(\lambda)$ are forced to each be $1,
2, \ldots, m$.  Hence, up to a $q$ shift, $\SSYT(\lambda)^{\rank}(q)$
equals $\SSYT(\mu)^{\rank}(q)$ where $\mu_{i}=\lambda_{i}-\lambda_{m}$.
In order to classify limit laws for $\SSYT_{\leq
m}(\lambda)^{\rank}(q)$, it thus suffices to assume throughout that
$\ell(\lambda) < m$ and $\lambda_m=0$.

We begin with a simple analogue of \Cref{thm:BKS}. This will be our only use
of the hook-content-based cumulant formula; all of our other results rely on the
$q$-Weyl dimension-based cumulant formula. 

\begin{Theorem}\label{thm:SSYT_cumulant_nm0}
Let $\lambda$ denote an infinite sequence of partitions with $|\lambda|=n$.
If $\frac{n}{m} \to 0$, then for each fixed $d \in \bZ_{\geq 2}$, the
corresponding sequence of cumulants is 
\begin{equation}\label{eq:cumulant.n/m.to.0}
\kappa_d^{\lambda; m} \sim
\frac{B_d}{d} nm^d.
\end{equation}
Furthermore, we can characterize convergence in distribution in the
case $\frac{n}{m} \to 0$ depending on the limiting value of $n$.
\begin{enumerate}[(i)]
\item If $n$ converges to a finite value $N$,
then $\cX_\lambda[\rank]^*$ converges in distribution to
$\cIH_N^*$.  \item If $n \to \infty$, then $\cX_\lambda[\rank]$ is
asymptotically normal.
\end{enumerate}
  
  \begin{proof}  For each cell $u \in \lambda$, 
    the trivial bounds $0 \leq c_u \leq n$ and $1 \leq h_u \leq n$ give
      \[ m^d - n^d \leq (m+c_u)^d - h_u^d \leq (m+n)^d.
      \]
Summing over all $u \in \lambda$ and dividing through by $nm^d$
gives
\[ 1 - \left(\frac{n}{m}\right)^d \leq \frac{\sum_{u \in
\lambda} (m+c_u)^d - h_{u}^d}{nm^d} \leq \left(1 +
\frac{n}{m}\right)^d.
\]
When $n/m \to 0$, the lower and upper bounds each tend to $1$.  The
cumulant formula \eqref{eq:cumulant.n/m.to.0} now follows from \eqref{eq:SSYT_cumulants.2}.

    By \eqref{eq:cumulant.n/m.to.0}, $(\kappa_d^{\lambda; m})^* \sim
(B_d/d)/(B_2/2)^{d/2}\cdot n^{1-d/2},$ which eliminates $m$ from the
limits.  If $n \to N$, then $(B_d/d)/(B_2/2)^{d/2}\cdot n^{1-d/2}$
approaches the $d^{th}$ cumulant of $\cIH_N^*$ by
\eqref{eq:rescaled.cumulants} and \Cref{ex:Ucont}, proving (i). If $n
\to \infty$, then $(B_d/d)/(B_2/2)^{d/2} \cdot n^{1-d/2}$ tends to $0$
for $d\geq 3$, which are the cumulants of $\cN(0,1)$, hence (ii)
follows from \Cref{cor:cumulants}.
\end{proof}
\end{Theorem}

\begin{Example}
  Consider a constant sequence of partitions $\lambda^{(N)}=\lambda$
and let $m \to \infty$.  By \Cref{thm:SSYT_cumulant_nm0}(i),
$\cX_{\lambda; m}[\rank]^* \Rightarrow \cIH_{|\lambda|}^*$, which
depends only on $|\lambda|$. On the other hand, if the sequence
$\lambda^{(N)}$ is chosen such that 
$|\lambda^{(N)}| \to \infty$ and $m^{(N)} \sim |\lambda^{(N)}|^2$, the limit is $\cN(0, 1)$
by \Cref{thm:SSYT_cumulant_nm0}(ii).
\end{Example}

\begin{Corollary}
  For any fixed $\epsilon > 0$, let 
    \[ \bfM_\epsilon \coloneqq \{\cX_{\lambda; m}[\rank]^* :   |\lambda| <
    m^{1-\epsilon} \}
        \subset \bfM_{\SSYT}. \]
  In the L\'evy metric,
  \begin{equation}\label{eq:bfM_SSYT_n_bounded}
    \overline{\bfM_\epsilon} = \bfM_\epsilon \sqcup \overline{\bfM_{\cIH}},
  \end{equation}
  which is (sequentially) compact. Moreover,
  the set of limit points of $\bfM_\epsilon$ is $\overline{\bfM_{\cIH}}$.
  
  \begin{proof}
    Given a sequence in $\bfM_\epsilon$, if $m$ is bounded, then so is $n=|\lambda |$, so there
    are only finitely many distinct $(\lambda, m)$ in the sequence and convergence
    occurs if and only if the sequence is eventually constant. On the other hand, if
    $m \to \infty$, then $n < m^{1-\epsilon}$ yields $\frac{n}{m} < m^{-\epsilon} \to 0$,
    so \eqref{eq:bfM_SSYT_n_bounded} follows immediately from
    \Cref{thm:SSYT_cumulant_nm0}.  
    
    By these observations, every infinite sequence of distinct points in $\bfM_\epsilon$
    has a limit point in $\overline{\bfM_{\cIH}}$, so $\bfM_\epsilon$ consists entirely of isolated
    points and $\overline{\bfM_{\cIH}}$ consists entirely of limit points. Sequential compactness is
    similarly clear.
  \end{proof}
\end{Corollary}

\subsection{Limit laws with $|\lambda |/m \to \infty$ and distance distributions}\label{ssec:ssyt_limits.nlarge}

At the other extreme, we may consider the case when $|\lambda |/m \to
\infty$. As we will see, the possible behavior is vastly more varied
in this limit.  Among the sequences of partitions $\lambda$ with $|\lambda |/m
\to \infty$, the easiest case to consider is when
$\lambda_{1}/m^{3}\to \infty$.  This includes the case where $m$
converges to a fixed finite value and $|\lambda| \to \infty$.

For a partition $\lambda=(\lambda_{1},\ldots,\lambda_{m})$, recall
from \Cref{thm:SSYT_Delta} that
\[
\mathbf{t}(\lambda)=(t_{1},\ldots,t_m) \in [0,1]^{m}
\]
is the finite multiset with $t_j \coloneqq
\frac{\lambda_j}{\lambda_1}$ for $1 \leq j \leq m$.  By
\Cref{def:distance_multiset}, the corresponding distance multiset is
\[ \Delta\mathbf{t}(\lambda ) \coloneqq
\{t_i - t_j : 1 \leq i < j \leq m\}. \]

\begin{Lemma}\label{lem:SSYT_nm_big_estimate}
Let $\lambda$ denote an infinite sequence of partitions with
$\ell(\lambda) < m$ and $|\lambda|=n$.  If $\frac{n}{m} \to \infty$ in
such a way that $\lambda_1/m^3 \to \infty$, then
for each fixed $d \in \bZ_{\geq 2}$,

    \[ \frac{\kappa_d^{\lambda; m}}{\lambda_1^d} \sim
(B_d/d)|\Delta\mathbf{t}(\lambda )|_d^d, \] which is the $d^{th}$
cumulant of the rescaled uniform sum $\cS_{\Delta\mathbf{t}(\lambda
)}/\lambda_{1}$.

  \begin{proof}
    Note that
      \[ (\lambda_i - \lambda_j)^d - m^d \leq (\lambda_i - \lambda_j + j-i)^d - (j-i)^d
          \leq (\lambda_i - \lambda_j + m)^d. \]
    Divide through by $\lambda_1^d$ and consider the upper bound.
Setting $t_j \coloneqq
\frac{\lambda_j}{\lambda_1}$ for $1 \leq j \leq m$,  
we find
    \begin{align*}
      \frac{(\lambda_i - \lambda_j + m)^d}{\lambda_1^d}
        &= (t_i - t_j + m/\lambda_1)^d \\
        &= \sum_{k=0}^d \binom{d}{k} (t_i - t_j)^{d-k}
             \left(\frac{m}{\lambda_1}\right)^k \\
        &\leq (t_i - t_j)^d + \sum_{k=1}^d \binom{d}{k}
          \left(\frac{m}{\lambda_1}\right)^k.
    \end{align*}
    Summing over all $1 \leq i < j \leq m$ and considering both bounds gives
    \begin{align*}
      \sum_{1 \leq i < j \leq m} (t_i - t_j)^d
        - \binom{m}{2} \cdot \left(\frac{m}{\lambda_1}\right)^d
        &\leq \frac{\sum_{1 \leq i < j \leq m}
                 (\lambda_i - \lambda_j + j-i)^d - (j-i)^d}{\lambda_1^d} \\
        &\leq \sum_{1 \leq i < j \leq m} (t_i - t_j)^d
           + \binom{m}{2} \sum_{k=1}^d \binom{d}{k} \left(\frac{m}{\lambda_1}\right)^k.
    \end{align*}
    Since $\lambda_1/m^3 \to \infty$, we have
    $\binom{m}{2} (m/\lambda_1) \to 0$. Furthermore, for each
    $k \geq 1$, 
      \[ \binom{m}{2} (m/\lambda_1)^k
          \leq m^{2+k}/\lambda_1^k \leq m^{3k}/\lambda_1^k \to 0. \]
     It follows that
      \[ \frac{\sum_{1 \leq i < j \leq m} (\lambda_i - \lambda_j + j - i)^d - (j-i)}{\lambda_1^d}
          \sim \sum_{1 \leq i < j \leq m} (t_i - t_j)^d. \]
    The result now follows from \eqref{eq:SSYT_cumulants.1} and \eqref{eq:Ucont.cumulants.2}.
  \end{proof}
\end{Lemma}

We use the results on generalized uniform sum distributions from
\Cref{sec:sums} to characterize convergence in distribution in the
next theorem. It is a more explicit statement of
\Cref{thm:SSYT_Delta}.

\begin{Theorem}\label{thm:SSYT_mfixed}
  Let $\lambda$ denote an infinite sequence of partitions, with
  $\ell(\lambda) < m$ and $|\lambda|=n$.  If $\frac{n}{m} \to \infty$ in
  such a way that $\lambda_1/m^3 \to \infty$, then for each fixed $d \in
  \bZ_{\geq 2}$, the standardized cumulants are approximately

  \begin{equation}\label{eq:cum_forlambda1overmcubed}
    (\kappa_d^{\lambda; m})^*
      \sim
	 \frac{B_d}{d}
	 |\widehat{\Delta\mathbf{t}(\lambda)}|_d^{d} = \kappa_d^{\cS_{\widehat{\Delta\mathbf{t}(\lambda)}}}.
  \end{equation}
  Furthermore, we  can characterize convergence in distribution when
  it occurs.
  \begin{enumerate}
    \item [(i)] If $m$ is bounded, then $\cX_{\lambda;
      m}[\rank]^*$ converges in distribution if and only if the 
      multisets $\widehat{\Delta\mathbf{t}(\lambda)}$ converge
      pointwise to some multiset $\mathbf{d}$,
      in which case, the limiting distribution is $\cS_{\mathbf{d}}$
      and $\mathbf{d} \in \bfP_{\DIST}$,
    \item [(ii)]  The sequence $m \to \infty$ if and only if $\cX_{\lambda; m}[\rank]^*$ is 
      asymptotically normal.
  \end{enumerate}
  \begin{proof} By hypothesis, $\lambda_1/m^3 \to \infty$, so
\Cref{lem:SSYT_nm_big_estimate} implies $\kappa_d^{\lambda; m} \sim
(B_d/d)|\Delta\mathbf{t}(\lambda)|_d^d$ for all $d\geq 2$.  Thus, the
standardized cumulants are given by \[ (\kappa_d^{\lambda; m})^* \sim
\frac{(B_d/d)|\Delta\mathbf{t}(\lambda)|_d^d} {(
(B_2/2)|\Delta\mathbf{t}(\lambda)|_2^{2})^{d/2}}
=\frac{B_d/d}{(B_2/2)^{d/2}}
         \left(
	 \frac{|\Delta\mathbf{t}(\lambda)|_d}{|\Delta\mathbf{t}(\lambda)|_2}
	 \right)^{d} =
 \frac{B_d}{d}
	 |\widehat{\Delta\mathbf{t}(\lambda)}|_d^{d}
\]
by the definition of the hat-operation \eqref{eq:hat.operation}.  By
\eqref{eq:Ucont.cumulants.2}, $ \frac{B_d}{d}
|\widehat{\Delta\mathbf{t}(\lambda)}|_d^d$ is the $d^{th}$ cumulant for
the uniform sum random variable
$\cS_{\widehat{\Delta\mathbf{t}(\lambda)}}$.

	By the Method of Moments/Cumulants (\Cref{thm:moments})
together with its converse in this context (\Cref{cor:frechet_shohat_converse}), the sequence
$\cX_{\lambda; m}[\rank]^*$ converges in distribution to some $\cX$ if
and only if the limit of the standardized cumulants 
$(\kappa_d^{\lambda; m})^* \to \kappa_d^{\cX}<\infty$ for each $d\geq 1$,
which happens  if and only if
$\kappa_d^{\cS_{\widehat{\Delta\mathbf{t}(\lambda)}}}
\to \kappa_d^\cX$ for each $d \geq 1$. By the Method of Moments/Cumulants
and its converse for DUSTPAN distributions
(\Cref{prop:MoM_converse}), this occurs if
and only if $\cS_{\widehat{\Delta\mathbf{t}(\lambda)}} \Rightarrow \cX$.
Finally, by \Cref{thm:PDelta_MDelta}, this occurs if and only if
$\widehat{\Delta\mathbf{t}(\lambda)}$ converges pointwise to some
$\mathbf{d} \in \overline{\bfP_{\DIST}}$.  The result follows from
\Cref{cor:MDelta_closure}. In particular, if $m$ is bounded (i) holds, and if
$m \to \infty$ (ii) holds.
\end{proof}
\end{Theorem}

\begin{Example}\label{ex:SSYT_cscaling.row}
  Fix a partition $\lambda$ and a positive integer $m>\ell(\lambda)$.
  Pick a sequence $r^{(N)} \to \infty$ of row scale factors, so that
  $\lambda^{(N)}_i = r^{(N)} \lambda_i$ and $m^{(N)} = m$.
  Clearly $\lambda_1^{(N)}/(m^{(N)})^3 \to \infty$, so
  by \Cref{thm:SSYT_mfixed}(i), we have
  $\cX_{r^{(N)}\lambda; m}[\rank]^* \Rightarrow \cS_{\Delta\lambda}^*$.
\end{Example}

\begin{Example}\label{ex:SSYT_geometric}
  Consider the sequence of partitions with $\lambda^{(N)} = (2^{N-1}, 2^{N-2}, \ldots, 1)$
  and $m^{(N)}=N$. Strictly speaking, $\ell(\lambda^{(N)}) = N = m^{(N)}$ here,
  so recall we can delete the first column and consider the auxiliary sequence $\mu^{(N)}
  = (2^{N-1}-1, 2^{N-2}-1, \ldots, 0)$. Now
    \[ \frac{\mu_1^{(N)}}{(m^{(N)})^3} = \frac{2^{N-1}-1}{N^3} \to \infty \]
  and $m^{(N)} = N \to \infty$. Thus
  $\cX_{(2^{N-1}, 2^{N-2}, \ldots, 1); N}[\rank]$ is asymptotically normal
  by \Cref{thm:SSYT_mfixed}(ii).
\end{Example}

\subsection{Limit laws based on distinct values in
$\lambda$ and the $\weft$ statistic}\label{ssec:ssyt_limits.weft}

We now describe a very general test for asymptotic normality of
$\cX_{\lambda; m}[\rank]$ based on a new statistic we call $\weft$ in
analogy with $\aft$ for standard Young tableaux. This test depends on
the number of distinct values in a partition, so we switch to
exponential notation.  Note, throughout the rest of this section $k$
will denote the number of distinct values in $\lambda$.

\begin{Definition}
  We may write a nonempty partition in \textit{exponential notation}
  $\lambda = \ell_1^{e_1} \cdots \ell_k^{e_k}$ where
  $\ell_1 > \cdots > \ell_k \geq 0$ and $e_i > 0$, meaning
  $\lambda$ has $e_i$ rows of length $\ell_i$. In
  our earlier notation, $m = e_1 + \cdots + e_k$ and
  $n = e_1 \ell_1 + \cdots + e_k \ell_k$.
\end{Definition}

\begin{Lemma}\label{lem:SSYT_theta}
  Take a partition $\lambda = (\lambda_{1},\dots
  ,\lambda_{m})=\ell_1^{e_1} \cdots \ell_k^{e_k}$. Then, uniformly for
  all $d \geq 2$,

  \begin{equation}\label{eq:big.sum.exponent.notation}
    \begin{split}
    \sum_{1 \leq i < j \leq m} (\lambda_i - \lambda_j)& (\lambda_i -
    \lambda_j + j - i)^{d-1}\\
      &= \Theta\left(
        \sum_{1 \leq a < b \leq k}   (\ell_a - \ell_b) e_a e_b
        (\ell_a - \ell_b - 1+ e_a + \cdots + e_b)^{d-1}\right).
    \end{split}
  \end{equation}

  \begin{proof}
    Observe that we may restrict the sum in
    \eqref{eq:big.sum.exponent.notation} to just the indices with
    $\lambda_i \neq \lambda_j$. Hence, we group the terms according to
    the distinct values $\lambda_i = \ell_a$ and $\lambda_j = \ell_b$ for
    $1 \leq a < b \leq k$. The contribution to the sum
    in \eqref{eq:big.sum.exponent.notation} for all $\lambda_i =
    \ell_a$ and $\lambda_j = \ell_b$ for a fixed $a<b$ is
    \begin{equation}\label{eq:big.sum.2}
      \sum (\ell_a - \ell_b) (\ell_a - \ell_b + j - i)^{d-1}
    \end{equation}
    where the sum is over $i, j$ such that 
    $e_1 + \cdots + e_{a-1} + 1 \leq i \leq e_1 + \cdots + e_a$
    and $e_1 + \cdots + e_{b-1} + 1 \leq j \leq e_1 + \cdots + e_b$.
    Reindexing with $p = e_a - (i - e_1 - \cdots - e_{a-1}) + 1$
    and $q = j - e_1 - \cdots - e_{b-1}$, the sum in
    \eqref{eq:big.sum.2} becomes
    \begin{align}\label{eq:big.sum.3}
      (\ell_a - \ell_b) \sum_{\substack{1 \leq p \leq e_a\\1 \leq q \leq e_b}}
        (\ell_a - \ell_b + p + q - 1 + e_{a+1} + \cdots + e_{b-1})^{d-1}.
    \end{align}

    Next, note that for fixed $d \geq 2$, $(u+v+w)^d = \Theta(u^d + v^d + w^d)$
    uniformly for all $u, v, w \geq 0$, since then
    \begin{align*}
      u^d + v^d + w^d
        &\leq (u+v+w)^d \\
        &\leq (3\max\{u, v, w\})^d = 3^d \max\{u^d, v^d, w^d\} \\
        &\leq 3^d (u^d + v^d + w^d).
    \end{align*}
    Letting $u=p$, $v=\ell_a-\ell_b+e_{a+1}+\cdots+e_{b-1}-1$, and
    $w=q$, we see the sum in \eqref{eq:big.sum.2} and \eqref{eq:big.sum.3} is $\Theta$ of
    \begin{align*}
      (\ell_a - \ell_b)
        &\sum_{\substack{1 \leq p \leq e_a \\ 1 \leq q \leq e_b}}
        \left[(\ell_a - \ell_b + e_{a+1} + \cdots + e_{b-1} - 1)^{d-1}
        + p^{d-1} + q^{d-1}\right] \\
        &= (\ell_a - \ell_b) \left[e_a e_b
             (\ell_a - \ell_b + e_{a+1} + \cdots + e_{b-1} - 1)^{d-1}
             + e_b \sum_{1 \leq p \leq e_a} p^{d-1}
             + e_a \sum_{1 \leq q \leq e_b} q^{d-1}\right].
    \end{align*}
    Since $d\geq 2$, $\sum_{1 \leq p \leq e_{a}} p^{d-1} = \Theta(e_{a}^d)$ uniformly
    for all $u \in \bZ_{\geq 1}$ by the sum bounds in
    \Cref{lem:riemann_estimate}, and similarly $\sum_{1 \leq p \leq
    e_{b}} q^{d-1} = \Theta(e_{b}^d)$.  
    Consequently, the preceding sum and also the sum in \eqref{eq:big.sum.2}
    are $\Theta$ of
    \begin{align*}
      (\ell_a - \ell_b) &\left[e_a e_b (\ell_a - \ell_b + e_{a+1} + \cdots + e_{b-1} - 1)^{d-1}
      + e_b e_a^d + e_a e_b^d\right] \\
        &= (\ell_a - \ell_b) e_a e_b \left[
             (\ell_a - \ell_b + e_{a+1} + \cdots + e_{b-1} - 1)^{d-1}
             + e_a^{d-1} + e_b^{d-1}\right] \\
        &= \Theta\left((\ell_a - \ell_b) e_a e_b
             (\ell_a - \ell_b - 1 + e_a + \cdots + e_b)^{d-1}\right).
    \end{align*}
    The result follows by summing over all $1 \leq a < b \leq k$, since the
    preceding bounds were all uniform.
  \end{proof}
\end{Lemma}

\begin{Theorem}\label{lem:SSYT_ansuff}
  Let $\lambda=\ell_1^{e_1} \cdots \ell_k^{e_k}$ denote an infinite
  sequence of partitions with $\ell(\lambda) \leq  m$,
  $\ell_{1}>\ell_{2}>\dots >\ell_{k}\geq 0$ and each $e_{i}>0$. Then, for
  $d\geq 2$ even, 

  \begin{equation}\label{eq:gen.cumulant}
    \kappa_d^{\lambda; m}
      = \Theta\left(\sum_{1 \leq a < b \leq k} (\ell_a - \ell_b)
          e_a e_b (\ell_a - \ell_b - 1+  e_a + \cdots + e_b)^{d-1}\right).
  \end{equation}
  Furthermore, $\cX_{\lambda; m}[\rank]$ is
  asymptotically normal if
  \begin{equation}\label{eq:aft-like-stat-ssyt}
    \weft(\lambda) \coloneqq \frac{\sum_{1 \leq a < b \leq k} (\ell_a - \ell_b)
      e_a e_b (\ell_a - \ell_b - 1 + e_a + \cdots + e_b)} {(\ell_1 - \ell_k -1 + m )^2} \to \infty.
  \end{equation}

  \begin{proof}
    In general, $u^d - v^d = (u-v) \sum_{i=0}^{d-1}u^{i}v^{d-i-1}
    =(u-v) \mathbf{h}_{d-1}(u, v)$, so
    \eqref{eq:SSYT_cumulants.1} gives
      \[ \kappa_d^{\lambda; m} = \sum_{1 \leq i < j \leq m}
          (\lambda_i - \lambda_j) \mathbf{h}_{d-1}(\lambda_i - \lambda_j + j-i, j-i). \]
    For fixed $d\geq 2$ even and $u \geq v \geq 0$, we have
    $\mathbf{h}_{d-1}(u, v) = \Theta(u^{d-1})$ since
      \[ u^{d-1} \leq u^{d-1} + u^{d-2} v + \cdots + v^{d-1}
          \leq u^{d-1} + u^{d-1} + \cdots + u^{d-1} = du^{d-1}. \]
    Consequently, 
      \[ \kappa_d^{\lambda; m} = \Theta\left(\sum_{1 \leq i < j \leq m}
          (\lambda_i - \lambda_j) (\lambda_i - \lambda_j + j - i)^{d-1}\right). \]
    Hence, \eqref{eq:gen.cumulant} holds by \Cref{lem:SSYT_theta}.

    We use the cumulant formula in \eqref{eq:gen.cumulant} to prove the
    asymptotic normality result. Write $x_{ab} \coloneqq (\ell_a -
    \ell_b) e_a e_b$ and $y_{ab} \coloneqq \ell_a - \ell_b - 1 + e_a +
    \cdots + e_b$ for $1\leq a<b\leq k$. By
    \eqref{eq:gen.cumulant}, we have for all $d \geq 2$ even

    \begin{align*}
      (\kappa_d^{\lambda; m})^*
        = \Theta\left(\frac{\sum_{1 \leq a < b \leq k} x_{ab} y_{ab}^{d-1}}
           {\left(\sum_{1 \leq a < b \leq k} x_{ab}
           y_{ab}\right)^{d/2}}\right).
    \end{align*}
    Note that $y_{1k} \geq y_{ab}$, so $\hat{y}_{ab} \coloneqq
    y_{ab}/y_{1k} \leq 1$. Hence $\hat{y}_{ab}^{d-1} \leq \hat{y}_{ab}$
    and
      \[ \sum_{1 \leq a < b \leq k} x_{ab} \hat{y}_{ab}^{d-1}
          \leq \sum_{1 \leq a < b \leq k} x_{ab} \hat{y}_{ab}. \]
    Consequently,
    \begin{align*}
      \frac{\sum_{1 \leq a < b \leq k} x_{ab} y_{ab}^{d-1}}
          {\left(\sum_{1 \leq a < b \leq k} x_{ab}
          y_{ab}\right)^{d/2}}
        &= \frac{y_{1k}^{d-1}}{y_{1k}^{d/2}}
          \frac{\sum_{1 \leq a < b \leq k} x_{ab} \hat{y}_{ab}^{d-1}}
          {\left(\sum_{1 \leq a < b \leq k} x_{ab} \hat{y}_{ab}\right)^{d/2}} \\
        &\leq \frac{1}{y_{1k}^{1 - d/2}}
          \left(\sum_{1 \leq a < b \leq k} x_{ab} \hat{y}_{ab}\right)^{1 - d/2} \\
        &= \left(\sum_{1 \leq a < b \leq k} x_{ab}
y_{ab}/y_{1k}^2\right)^{1 - d/2}.  \end{align*} The latter
parenthesized quantity equals the $\weft(\lambda )$ statistic in
\eqref{eq:aft-like-stat-ssyt} by
construction. Thus, \begin{equation}\label{eq:normalized.cum.weft}
(\kappa_d^{\lambda; m})^* = O\left(\frac{1}{\weft(\lambda
)^{d/2-1}}\right).  \end{equation} When $d \geq 4$, $\weft(\lambda )
\to \infty$ implies $(\kappa_d^{\lambda; m})^* \to 0$.  Thus,
asymptotic normality follows from \Cref{cor:cumulants} since all of
the odd Bernoulli numbers for $d\geq 3$ are zero.
\end{proof}
\end{Theorem}

\begin{Example}\label{cor:SSYT_staircase}
  Let $\lambda^{(N)} = \delta_N \coloneqq (N-1,N-2, \ldots, 2, 1,0)$
  be the staircase partition for $N>1$.  We have $e_1 = \cdots = e_N = 1$ and
  $\ell_1 = N-1, \ldots, \ell_N = 0$. In this case, \eqref{eq:aft-like-stat-ssyt} simplifies to
    \[ \weft(\lambda^{(N)})= \frac{\sum_{1 \leq a < b \leq N}
        2(b-a)^{2}}{(2N - 2)^{2}} =N^2 \frac{N+1}{24(N-1)}. \]
  Thus, as $N\to \infty$ this statistic goes to infinity, so
  $\cX_{\delta_N; N}[\rank]$ is asymptotically normal by
  \Cref{lem:SSYT_ansuff}.
\end{Example}

The characterization in \Cref{lem:SSYT_ansuff} is powerful enough to
prove asymptotic normality in many cases of interest.  We will use the
criteria in the next corollary to further simplify the arguments in
the examples below and the applications to plane partitions.  As
mentioned in the introduction to this section, we can assume
$\ell(\lambda)<m$ without loss of generality. We may include the case
$\ell(\lambda)=m$ if desired by replacing $\ell_1$ with $\ell_1 -
\ell_k$ in the following result.

\begin{Corollary}\label{thm:SSYT_technical}
  Let $\lambda=\ell_1^{e_1} \cdots \ell_k^{e_k}$
  denote an infinite sequence of partitions with
  $\ell(\lambda) <  m$, so
  $\lambda_{1}=\ell_{1}>\ell_{2}>\dots >\ell_{k} =  0$,
  and each $e_{i}>0$. Then
  $\cX_{\lambda; m}[\rank]$ is asymptotically normal in the following
  situations.

  \begin{enumerate}[(i)]
    \item $\frac{m^2}{k\ell_1(k+\ell_1)} \to 0$ and $k \to \infty$.
    \item $\frac{e^{[2]}}{k(\ell_1/m + 1)^2} \to \infty$, where
      $e^{[2]}$ denotes the second largest element among $e_1, \ldots, e_k$.
    \item $\frac{\ell_1 e_1 e_k}{\ell_1 + m} \to \infty$.
  \end{enumerate}

  \begin{proof}
    For (i), suppose $\frac{m^2}{k\ell_1(k+\ell_1)} \to 0$ and $k \to \infty$.
    We have
    \begin{align*}
      \sum_{1 \leq a < b \leq k} &(\ell_a - \ell_b) e_a e_b
          (\ell_a - \ell_b - 1 + e_a + \cdots + e_b) \\
        &\geq \sum_{1 < p < k}
           [(\ell_1 - \ell_p) e_1 e_p (\ell_1 - \ell_p - 1 + e_1 + \cdots + e_p)
             + \ell_p e_p e_k (\ell_p - 1 + e_p + \cdots + e_k)] \\
        &\geq \sum_{1 < p < k} [(\ell_1 - \ell_p)(\ell_1 - \ell_p + p-1)
         + \ell_p (\ell_p + k-p)] \\
        &\geq \sum_{\frac{k}{4}+1 < p < \frac{3k}{4}}
             [(\ell_1 - \ell_p)(\ell_1 - \ell_p + k/4) + \ell_p (\ell_p + k/4)] \\
        &= \sum_{\frac{k}{4}+1 < p < \frac{3k}{4}}
             [\ell_1^2 - 2 \ell_1 \ell_p + 2 \ell_p^2 + k \ell_1/4].
    \end{align*}
    Set $x_p \coloneqq \ell_p/\ell_1$ and divide the preceding
    inequality by $\ell_1^2$. Suppose $k \geq 4$.
    The final expression becomes
    \begin{align*}
      \sum_{\frac{k}{4}+1 < p < \frac{3k}{4}}
          [1 - 2x_p + 2x_p^2 + k/(4\ell_1)]
        &\geq \sum_{\frac{k}{4}+1 < p < \frac{3k}{4}}
          [1/2 + k/(4\ell_1)] \\
        &\geq (k/2-2) (1/2+k/(4\ell_1)) \\
        &\geq k/16 \cdot (1 + k/\ell_1).
    \end{align*}
    Consequently,
    \begin{align*}
      \weft(\lambda)
        &\geq \frac{k/16 \cdot (1+k/\ell_1)}{(1+(m-1)/\ell_1)^2} \\
        &\geq \frac{1}{16} \frac{k\ell_1(k+\ell_1)}{(\ell_1 + m)^2} \\
        &\geq \frac{1}{16} \min\left\{\frac{k\ell_1(k+\ell_1)}{(2\ell_1)^2},
            \frac{k\ell_1(k+\ell_1)}{(2m)^2}\right\} \\
        &\geq \frac{1}{64} \min\left\{k, \frac{k\ell_1(k+\ell_1)}{m^2}\right\} \to \infty,
    \end{align*}
since $\frac{m^2}{k\ell_1(k+\ell_1)} \to 0$ and $k \to \infty$ by
hypothesis. The result now follows from \Cref{lem:SSYT_ansuff}.

    For (ii), suppose $\frac{e^{[2]}}{k(\ell_1/m + 1)^2} \to \infty$.
    By definition, $e^{[2]}\leq m-e_{i}$ for all $i$.   
    Thus,
    \begin{align*}
      \sum_{1 \leq a < b \leq k} (\ell_a - \ell_b) &e_a e_b
        (\ell_a - \ell_b - 1 + e_a + \cdots + e_b) \\
        &\geq \sum_{1 \leq a < b \leq k} e_a^2 e_b + e_a e_b^2 \\
        &= \sum_{1 \leq i \leq k} e_i^2 (e_{i+1} + \cdots + e_k)
        + \sum_{1 \leq i \leq k} (e_1 + \cdots + e_{i-1}) e_i^2 \\
        &= \sum_{1 \leq i \leq k} e_i^2 (m - e_i) \\
        &\geq e^{[2]} \sum_{1 \leq i \leq k} e_i^2.
    \end{align*}
    If $f_1 \geq \cdots \geq f_k \geq 0$, then
    \begin{align*}
      (f_1 + \cdots + f_k)^2
        &= f_1^2 + \cdots + f_k^2 + 2\sum_{1 \leq i < j \leq k} f_i f_j \\
        &\leq f_1^2 + \cdots + f_k^2 + 2\sum_{1 \leq i < j \leq k} f_i^2 \\
        &= \sum_{i=1}^k (1 + 2(k-i)) f_i^2 \\
        &\leq 2k \sum_{i=1}^k f_i^2.
    \end{align*}
    This latter bound is independent of the actual order of the $f_i$.
    Consequently,
      \[ e^{[2]} \sum_{1 \leq i \leq k} e_i^2
          \geq e^{[2]} \frac{m^2}{2k}. \]
    Clearly $(\ell_1 - \ell_k - 1 + m)^2 \leq (\ell_1 + m)^2$. Hence
    \begin{align*}
      \weft(\lambda)
        &= \frac{\sum_{1 \leq a < b \leq k} (\ell_a - \ell_b) e_a e_b
        (\ell_a - \ell_b - 1 + e_a + \cdots + e_b)}{(\ell_1 - \ell_k - 1 + m)^2} \\
        &\geq \frac{e^{[2]} m^2}{2k(\ell_1 + m)^2} =
\frac{e^{[2]}}{2k(\ell_1/m + 1)^2} \to \infty,
    \end{align*}
since
$\frac{e^{[2]}}{k(\ell_1/m + 1)^2} \to \infty$. The result again
follows from \Cref{lem:SSYT_ansuff}.

    For (iii), suppose $\frac{\ell_1 e_1 e_k}{\ell_1 + m} \to \infty$. We have
    \begin{align*}
      \weft(\lambda)
        &\geq \frac{(\ell_1 - \ell_k) e_1 e_k (\ell_1 - \ell_k - 1 + m)}{(\ell_1 - \ell_k - 1 + m)^2} \\
        &\geq \frac{\ell_1 e_1 e_k}{\ell_1 + m} \to \infty.
    \end{align*}
    The result again follows from \Cref{lem:SSYT_ansuff}.
  \end{proof}
\end{Corollary}


\begin{Example}\label{ex:SSYT_distinct}
  Suppose $\lambda$ is a sequence of partitions with distinct parts
  and $\ell(\lambda) \to \infty$. Then $\ell_1 \geq k = m$ and
$\frac{m^2}{k\ell_1 (k+\ell_1)} \leq \frac{1}{k} \to 0$. By
\Cref{thm:SSYT_technical}(i), the sequence $\cX_{\lambda; m}[\rank]$
is asymptotically normal.
\end{Example}

\begin{Example}\label{ex:SSYT_square}
  Suppose $\lambda$ is a sequence of partitions with $m = \ell_1$
  and $k \to \infty$. Then $\frac{m^2}{k\ell_1 (k+\ell_1)} \leq \frac{1}{k} \to 0$.
  Again by \Cref{thm:SSYT_technical}(i), the sequence $\cX_{\lambda; m}[\rank]$
  is asymptotically normal.
\end{Example}

\begin{Remark}
  The limit shape of a randomly chosen partition of $n$ as $n \to \infty$
  is well-known to be the curve
    \[ e^{-\frac{\pi}{\sqrt{6}} x} + e^{-\frac{\pi}{\sqrt{6}} y} = 1 \]
  where $(x, y)$ corresponds to $(i/\sqrt{n}, \lambda_i/\sqrt{n})$
  \cite[Thm.~4.4, p.99]{MR1402079}. One consequently expects
  $\lambda_1 \approx \sqrt{n}$, and certainly $k \to \infty$.
  It seems natural to use $m = \lambda_1' \approx \sqrt{n}$, in which
  case
    \[ \frac{m^2}{k\lambda_1(k+\lambda_1)}
        \approx \frac{\sqrt{n}^2}{k\sqrt{n}(k+\sqrt{n})}
        \leq \frac{1}{k} \to 0. \]
  Thus, one heuristically expects $\cX_{\lambda; m}[\rank]$ to be
  asymptotically normal for randomly chosen partitions. We do not attempt
  to make this precise.
\end{Remark}

\begin{Question}\label{question:ktoinf}
  Suppose $\lambda$ is a sequence of partitions with $\ell(\lambda) <
m$ and $k$ is the number of distinct parts of $\lambda$.  Does $k \to
\infty$ ensure $\cX_{\lambda; m}[\rank]$ is asymptotically normal?
\end{Question}

\subsection{Limit laws for plane partitions}\label{ssec:plane.parts}
We may use \Cref{thm:SSYT_technical}(ii) to deduce the complete
characterization of the asymptotic limits for plane partitions in a
box from the introduction.  The following is a restatement of
\Cref{thm:PP_median.intro}.

\begin{Theorem}
  The $\size$ statistic on $\PP(a \times b \times c)$ is
  asymptotically normal if and only if
    \[ \median\{a, b, c\} \to \infty. \]
  If $ab$ converges and $c \to \infty$, the normalized limit law is the Irwin--Hall
  distribution $\cIH_{ab}^*$.
\end{Theorem}

\begin{proof}
  From the discussion in \Cref{ssec:back:ssyt}, we have
    \[ \cX_{\PP(a \times b \times c)}[\size]^*
        = \cX_{\SSYT_{\leq a+c}((b^a))}[\rank]^*. \]
  Let $\lambda = (b^a) = b^a 0^c$, so $n=ab$, $k=2$,
  $\ell_1=b$, $\ell_2=0$, $e_1=a$, $e_2=c$, and $m=a+c$.
  Suppose $\median\{a, b, c\} \to \infty$. Without loss of
  generality, we may suppose $b \leq a \leq c$, so $a \to \infty$.
  In this case, $e^{[2]} = a$ and $b/(a+c) \leq 1/2$. Hence
  $\frac{e^{[2]}}{k(\ell_1/m + 1)^2} = \frac{a}{2(b/(a+c)+1)^2}
  \geq \frac{a}{2(3/2)^2} \to \infty$ and asymptotic
  normality follows from \Cref{thm:SSYT_technical}(ii).

  On the other hand, if $\median\{a, b, c\}$ is bounded, we may
  suppose $a \leq b \leq c$, so that $n = ab$ is bounded. If $c \to
  \infty$, then the standardized limit distribution is $\cIH_{ab}^*$
  provided $ab$ converges by \Cref{thm:SSYT_cumulant_nm0}(i). The result
  follows.
\end{proof}

We conclude this section by giving some sample applications
of the preceding results to three natural scaling limits of partitions
obtained by stretching rows and/or columns by scale factors
tending to $\infty$.

\begin{Example}\label{ex:SSYT_cscaling.col}
  Continuing \Cref{ex:SSYT_cscaling.row}, instead pick a sequence
  $c^{(N)} \to \infty$ of column scale factors, so that
    \[ \lambda^{(N)} = (\underbrace{\lambda_1, \ldots, \lambda_1}_{c^{(N)}},
          \cdots, \underbrace{\lambda_m, \ldots, \lambda_m}_{c^{(N)}}), \]
  $m^{(N)} = c^{(N)} m$, $\ell_1^{(N)} = \lambda_1$, $e_i^{(N)} = c^{(N)}e_i$,
   and $(e^{(N)})^{[2]}
  = c^{(N)} e^{[2]}$. Thus
    \[ \frac{(e^{(N)})^{[2]}}{k^{(N)} (\ell_1^{(N)}/m^{(N)} + 1)^2}
       = \frac{c^{(N)}e^{[2]}}{k(\lambda_1/(c^{(N)}m)+1)^2}
       \sim \frac{c^{(N)}e^{[2]}}{k} \to \infty, \]
  so by \Cref{thm:SSYT_technical}(ii),
  $\cX_{(c^{(N)}\lambda')'; c^{(N)}m}[\rank]$ is asymptotically normal.
\end{Example}

\begin{Example}\label{ex:SSYT_cscaling.both}
  Combining \Cref{ex:SSYT_cscaling.row} and \Cref{ex:SSYT_cscaling.col},
  use both row and column scale factors simultaneously.
  We see
  \begin{align*}
    \frac{\ell_1^{(N)} e_1^{(N)} e_k^{(N)}}{\ell_1^{(N)} + m^{(N)}}
      &= \frac{r^{(N)} (c^{(N)})^2 \lambda_1 e_1 e_k}{r^{(N)}\lambda_1 + c^{(N)}m}
      \to \infty,
  \end{align*}
  so by \Cref{thm:SSYT_technical}(iii),
  $\cX_{r^{(N)}(c^{(N)}\lambda')'; c^{(N)}m}[\rank]$ is asymptotically normal.
  In particular, this includes the case when $c^{(N)} = r^{(N)} \to \infty$ and
  $\lambda^{(N)}$ is obtained from $\lambda$ by replacing each cell with
  a $c^{(N)} \times c^{(N)}$ grid of cells.
\end{Example}

\subsection{Summary}

Here we collect the known cases when $\cX_{\lambda; m}[\rank]^*$ converges
in distribution. Let $n = |\lambda|$, without
loss of generality suppose $\ell(\lambda) < m$,
let $k$ be the number of distinct row lengths of $\lambda$ (including $0$
since $\ell(\lambda) < m)$, let $e_i$ be the multiplicity of the $i$th largest
row length, and let $e^{[2]}$ be the second-largest element
amongst $e_1, e_2, \ldots, e_k$.

\begin{Summary}\label{sum:SSYT}
  \ 
\begin{enumerate}[(i)]
  \item In the following situations,
    $\cX_{\lambda; m}[\rank]^* \Rightarrow \cN(0, 1)$.
    \begin{enumerate}[(a)]
      \item $\frac{n}{m} \to 0$ and $n \to \infty$ (\Cref{thm:SSYT_cumulant_nm0}(ii))
      \item $\lambda_1/m^3 \to \infty$ and $m \to \infty$. Moreover, a
        converse holds. (\Cref{thm:SSYT_mfixed}(ii))
      \item $\lambda^{(N)} = (2^{N-1}, 2^{N-2}, \ldots, 1)$ and $m^{(N)} = N$
        (\Cref{ex:SSYT_geometric})
      \item $\lambda^{(N)} = \delta_N = (N-1, N-2, \ldots, 2, 1, 0)$ and $m^{(N)} = N$
        (\Cref{cor:SSYT_staircase})
      \item $\weft(\lambda) \to \infty$ (\Cref{lem:SSYT_ansuff})
      \item $\frac{m^2}{k\ell_1(k+\ell_1)} \to 0$ and $k \to \infty$
        (\Cref{thm:SSYT_technical}(i))
      \item $\frac{e^{[2]}}{k(\ell_1/m + 1)^2} \to \infty$ (\Cref{thm:SSYT_technical}(ii))
      \item $\frac{\ell_1 e_1 e_k}{\ell_1 + m} \to \infty$ (\Cref{thm:SSYT_technical}(iii))
      \item $e_1 = \cdots = e_k = 1$ and $k \to \infty$ (\Cref{ex:SSYT_distinct})
      \item $m = \lambda_1$ and $k \to \infty$ (\Cref{ex:SSYT_square})
      \item $\lambda = (b^a)$, $m=a+c$, and $\median\{a, b, c\} \to \infty$
        (\Cref{thm:PP_median.intro})
      \item If the sequence $\lambda$ is obtained by successively scaling the
        columns by a factor $c \to \infty$. (\Cref{ex:SSYT_cscaling.col})
      \item If the sequence $\lambda$ is obtained by successively scaling
        the rows and columns by factors of $r, c \to \infty$.
        (\Cref{ex:SSYT_cscaling.both})
    \end{enumerate}
  \item In the following situations, $\cX_{\lambda; m}[\rank]^* \Rightarrow \cIH_M^*$.
    \begin{enumerate}[(a)]
      \item $n/m \to 0$ and $n \to M$. (\Cref{thm:SSYT_cumulant_nm0}(i))
      \item $\lambda = (b^a)$, $m=a+c$, $ab \to M$, and $c \to \infty$
        (\Cref{thm:PP_median.intro})
    \end{enumerate}
  \item In the following situations, $\cX_{\lambda; m}[\rank]^* \Rightarrow
    \cS_{\mathbf{d}}^*$.
    \begin{enumerate}
      \item $\lambda_1 \to \infty$, $m$ is bounded, and
        $\Delta\mathbf{t}(\lambda) \to \mathbf{d}$ where $x_i \coloneqq \lambda_i/\lambda_1$.
        Moreover, a converse holds. (\Cref{thm:SSYT_mfixed}(i))
      \item If the sequence $\lambda$ is obtained by successively scaling the
        rows by a factor $r \to \infty$, and $\mathbf{d} = \Delta\lambda$.
        (\Cref{ex:SSYT_cscaling.row})
    \end{enumerate}
\end{enumerate}
\end{Summary}

\section{Metric spaces related to forest distributions}\label{sec:forests}

In this section, we consider the two $q$-analogs of the number of
linear extensions of posets which come from trees and forests
using variations on the $\inv$ and $\maj$ statistics for permutations as
given by Bj\"orner--Wachs in \cite{MR1022316}.  Recall the background
for these $q$-analogs from \Cref{ssec:forests.background}. As summarized in
\Cref{ssec:intro:forests}, we will show that
the coefficients in the corresponding polynomials ``generically'' are
asymptotically normal, but that the metric space of DUSTPAN
distributions  $\bfM_{\DUSTPAN}$ characterizes all possible limit laws in a
certain degenerate regime. In particular, we prove
\Cref{thm:generic.forest.sum}, \Cref{thm:forest_infty}, and
\Cref{cor:forest_infty}.

\subsection{Generic asymptotic normality for trees and forests}

Recall from \Cref{ssec:forests.background} that for any forest $P$, 
there is an associated $q$-hook length polynomial 
\[
\cL_{P}(q) \coloneqq [n]_q!/\prod_{u \in P} [h_u]_q
\]
and random variable $\cX_P$.  Here we show that the sequences of
random variables $\cX_P$ for forests $P$ are asymptotically normal if
certain numerical conditions hold; see \Cref{thm:generic.forest.sum}.
This covers the ``generic cases''.  We begin by describing a family of
trees which maximize the sum of the hook lengths over all trees of
rank $r$ with $n$ elements.  We use this family of trees to identify
good approximations for the cumulants corresponding with all trees.

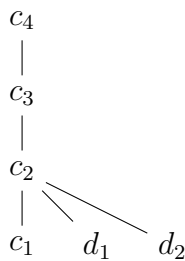
\begin{figure}[h]
\begin{center}
\begin{tikzpicture}
  \node (a1) at (0,0) {$c_1$};
  \node (a2) at (0,1) {$c_2$};
  \node (a3) at (0,2) {$c_3$};
  \node (a4) at (0,3) {$c_4$};
  \node (b1) at (1, 0) {$d_1$};
  \node (b2) at (2, 0) {$d_2$};
  \draw (a1) -- (a2) -- (a3) -- (a4);
  \draw (b1) -- (a2);
  \draw (b2) -- (a2);
\end{tikzpicture}

\end{center}
\caption{The poset $H_{6, 4}$. The chain $C=\{c_{1},c_{2}, c_{3},
c_{4} \}$ of length $4$ has $2$ additional descendants added to the
second-smallest element $c_{2}$.}  \label{fig:Hnk}
\end{figure}

\begin{Definition}
  Suppose $n \in \bZ_{\geq 1}$ and $1 < r \leq n$. Let
  $H_{n, r}$ be the tree obtained by starting with a
  rooted chain $C$ with $r$ elements and adding $n-r$ elements each
  as children of the second-smallest node in the chain. See \Cref{fig:Hnk}.
\end{Definition}

\begin{Lemma}\label{lem:forest_degree}
  Among all trees $P$ with $n$ elements and rank $1<r\leq n$,
  $H_{n, r}$ is the unique maximizer of $\sum_{u \in P} h_u$.
  Consequently,  the degree of $\cL_{P}(q)$ is
  \begin{equation}\label{eq:forest_degree}
    \sum_{k=1}^n k - \sum_{u \in P} h_u
        \geq \sum_{k=1}^n k - \sum_{v \in H_{n, r}} h_v
        = \binom{n-r+1}{2}
  \end{equation}
  
  \begin{proof}
Let $C$ be a
maximal chain of $P$ with $r>1$ elements and second-smallest element
$y$.  If $P \neq H_{n, r}$, let $x \in P - C$ be a leaf of $P$ which
is not a child of $y$. Let $P'$ be the result of moving $x$ to be a
descendant $x'$ of $y$, which preserves the rank and number of
vertices. Since $C$ is maximal, we can easily determine the change in
the sum of the hook lengths: it increases by $\#\{v' \in P' : v' \in
C, v' \geq y\} = |C| - 1 = r-1$ and decreases by $\#\{v \in P : v >
x\} \leq r-1$.  This procedure always weakly increases the sum of the
hook lengths and arrives at $H_{n,r}$ after a finite number of
iterations, so the maximality claim follows.

    Observe the procedure strictly increases the sum of the hook
lengths unless $\#\{v \in P : v > x\} = r-1$.  In this case, let $z$
be the unique cover of $x$ in $P$. By construction,  $z \not \in C$.
After applying the procedure to $x$ to get $P'$, applying the
procedure again to all of $z$'s children and then to $z$ will strictly increase
the sum of hook lengths.  Thus, $P$ has strictly smaller sum of hook
lengths than $H_{n,r},$ and the uniqueness claim follows.

For the equality in \eqref{eq:forest_degree}, we find
      \[ \sum_{v \in H_{n, r}} h_u = 1 \cdot (n-r+1) +
      \sum_{k=n-r+2}^n k = \sum_{k=n-r+1}^n k. \]
Therefore, 
       \[ \sum_{k=1}^n k - \sum_{v \in H_{n, r}} h_v
           = \sum_{k=1}^n k - \sum_{k=n-r+1}^n k
           = \sum_{k=1}^{n-r} k = \binom{n-r+1}{2}. \]
  \end{proof}
\end{Lemma}

\begin{Lemma}\label{lem:forest_generic_theta}
  Suppose $0 \leq \alpha < 1$ and fix $d \in \bZ_{\geq 2}$ even. Uniformly
  for all trees $P$ with $n$ elements and rank $1<r \leq \alpha n$, we have
    \[ |\kappa_d^P| = \Theta(n^{d+1}). \]
  Explicitly,  for a fixed $d \in \bZ_{\geq 1}$, 
\begin{equation}\label{eq:tree.bounds}
 \frac{ba^d}{d} n^{d+1} \leq \sum_{k=1}^n k^d - \sum_{u \in P} h_u^d
\leq \left(\frac{1}{d+1} + \frac{1}{n}\right) n^{d+1} 
\end{equation}
  where $x \coloneqq \left[\left(\frac{2}{1 - \alpha}\right)^2 - 1\right] > 1$,
  $a \coloneqq 1/x$, and $b \coloneqq 1/(x+1)$, so $0<a,b<1$.
  
  \begin{proof} Recall from \Cref{cor:forest_cumulants} that $
\kappa_d^P = \frac{B_d}{d} \left(\sum_{k=1}^n k^d - \sum_{u \in P}
h_u^d\right)$ so $|\kappa_d^P| = \Theta(n^{d+1})$ provided the lower
bound and upper bound in \eqref{eq:tree.bounds} hold.  The upper bound
follows from the upper bound in \Cref{lem:riemann_estimate}.  

    For the lower bound, construct a labeling $w$ of $P$ by iteratively
building up $P$ as follows.
Begin by labeling the root of $P$ with $1$ in $w$. At each step,
increment all existing labels in $w$, pick an element of $P$ which has
not been labeled whose parent has been labeled, and label it with
$1$. Observe that the resulting labeling $w \colon P \to [n]$ is
natural. Consider the quantity $w(u) - h_u$ during this procedure. When
$u$ has initially been labeled, we have $w(u) - h_u = 1 - 1 = 0$. After
$u$ has been labeled, when adding a new vertex $v$, if $v \leq u$ then
both $w(u)$ and $h_u$ are incremented, while if $v \not\leq u$ then
only $w(u)$ in incremented.  Consequently, the final value of $w(u) -
h_u$ counts the number of elements $v$ added after $u$ such that $v
\not\leq u$. In particular, $w(u) - h_u \geq 0$.

Using the real numbers $a,b,x$ defined in the statement of the lemma,
let $M \coloneqq \{u \in P : w(u) - h_u \geq bn\}$.  We claim $\#M
\geq an$.  To prove the claim, suppose to the contrary that $\#M <
an$. By definition, $0 < a, b < 1$. Consequently,

    \begin{align*}
      \sum_{u \in P} w(u) - h_u
        &= \sum_{u \in M} (w(u) - h_u) + \sum_{u \not\in M} (w(u) - h_u) \\
        &\leq \#M \cdot n + (n-\#M) \cdot bn \\
        &= bn^2 + \#M \cdot (1-b)n \\
        &< bn^2 + a(1-b)n^2 = (a+b-ab)n^2.
    \end{align*}
    One may easily check that $a+b-ab = 2/(x+1) = (1-\alpha)^2/2$. Since
    $r \leq \alpha n$, we have $n-r \geq (1-\alpha)n$, so that
    \begin{align*}
      \sum_{u \in P} w(u) - h_u < \frac{(1-\alpha)^2}{2} n^2
         \leq \frac{(n-r)^2}{2}
         \leq \binom{n-r+1}{2},
    \end{align*}
    contradicting \Cref{lem:forest_degree} and verifying the claim.
    Using the claim and the lower bound on the sum in \Cref{lem:riemann_estimate}, we now find
    \begin{align*}
      \sum_{j=1}^n j^d - \sum_{u \in P} h_u^d
        &= \sum_{u \in P} (w(u)^d - h_u^d) \\
        &\geq \sum_{u \in M} (w(u) - h_u) \mathbf{h}_{d-1}(w(u), h_u)
        \geq \sum_{u \in M} (bn) w(u)^{d-1} \\
        &\geq bn \cdot \sum_{j=1}^{\#M} j^{d-1} \\
        &\geq bn \cdot (\#M)^d/d
        \geq \left(ba^d/d  \right)\ n^{d+1}.
    \end{align*}
  \end{proof}
\end{Lemma}

Now, we are prepared to address the question of asymptotic normality
for sequences of random variables associated to trees and forests.
Recall the following theorem from the introduction.

\begin{repTheorem}{thm:generic.forest.sum}
  Given a sequence of forests $P$, the corresponding
  sequence of random variables $\cX_P^*$ is asymptotically
  normal if
    \[ |P| \to \infty \qquad \text{and} \qquad \limsup \frac{\rank(P)}{|P|} < 1. \]


  \begin{proof} By \Cref{rem:standard.trees}, it suffices to assume
$P$ is a tree.  For $d \geq 2$ even, we know $|\kappa_d^P| =
\Theta(n^{d+1})$ by \Cref{lem:forest_generic_theta}, so
$|(\kappa_d^P)^*| = |\kappa_d^P|/|\kappa_2^P|^{d/2} =
\Theta(n^{1-d/2}) \to 0$.  By \Cref{cor:forest_cumulants}, the odd
cumulants vanish. Therefore, the result again follows from
\Cref{cor:cumulants}.
\end{proof}
\end{repTheorem}

\begin{Remark}
One expects most random forest generation techniques to yield a rank
which is logarithmic in the number of nodes with high probability, in
which case \Cref{thm:generic.forest.sum} applies. This is
the sense in which we consider \Cref{thm:generic.forest.sum} to cover ``generic''
trees and forests.
\end{Remark}

\begin{Remark}\label{rem:explicit.bounds.forests}
  More precisely, we may use the explicit bounds in
  \Cref{lem:forest_generic_theta}. Setting $\alpha := r/n$, the lower bound becomes
  $\frac{(1-\alpha)^{2(d+1)}}{4(1+\alpha)^d(3-\alpha)^d} n^{d+1}$.
  Since $0 \leq \alpha \leq 1$, the denominator can be ignored. Considering
  the $d=4$ case for simplicity, we find 
    \[ \kappa_4^*
        = O\left(\frac{n^5}{\left((1-\alpha)^6 n^3\right)^2}\right)
        = O\left(\frac{n^{-1}}{(1-r/n)^{12}}\right)
        = O\left(\frac{n^{11}}{(n-r)^{12}}\right). \]
  Thus asymptotic normality follows when $\frac{n-r}{n^{11/12}} \to \infty$, or
  equivalently when $n-r = \omega(n^{11/12})$. By contrast,
  \Cref{thm:forest_infty} classifies limit laws when $n-r = o(n^{1/2})$.
  Analyzing the possible asymptotic behavior between these extremes is
  still an open problem.
\end{Remark}

\subsection{Degenerate forests and DUSTPAN
distributions}\label{ssec:forests_nn} We now consider sequences of
random variables associated to the ``degenerate'' trees with $n-r =
o(n^{1/2})$.   Note, $n-r =
o(n^{1/2})$ implies $r/n \to 1$, so  these sequences are not covered by \Cref{thm:generic.forest.sum}.
For such trees, we give a simple numerical estimate for the cumulants
in terms of \textit{multisets of elevations}, and use them to
characterize asymptotic normality as well as the other limiting
distributions in terms of the metric space of DUSTPAN distributions
$\bfM_{\DUSTPAN}$.

\begin{Remark}
  To avoid certain redundancies, we restrict to standardized
  trees in the sense of \Cref{rem:standard.trees}. As an example
  of behavior which is prohibited by this assumption, consider the
  trees $H_{n, n-k}$ for fixed $k$, which are not standardized.
  This sequence of trees has rank $r=n-k$, so $\lim r/n = 1$ as $n \to \infty$,
  and \Cref{thm:generic.forest.sum} does not apply.
  Indeed, it is easy to see that
    \[ \cL_{H_{n,n-k}}(q) = \frac{[n]_q!}{\prod_{u \in H_{n, n-k}} [h_u]_q} = [k+1]_{q}!. \]
  Therefore, $\cX_{H_{n, n-k}}^*$ has the same discrete distribution
  for all $n>k$, so the limit distribution is discrete.
  
  On the other hand, if $n-r \to \infty$, the length of the support
  of $\cX_P$ tends to $\infty$ by \Cref{rem:forest_degree} and
  \Cref{lem:forest_degree}. Hence each distribution appears only finitely many times
  in such a sequence. Moreover, since the coefficients are unimodal, any sequence
  $\cX_P^*$ with $n-r \to \infty$ cannot converge to a discrete distribution.
\end{Remark}

We begin with a series of estimates relating the cumulants
$\kappa_d^P$ to the following auxiliary combinatorial quantity on $P$.

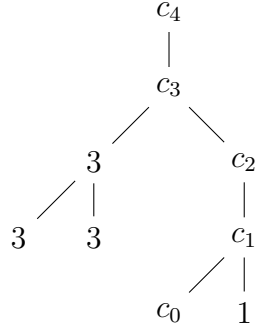
\begin{figure}
\begin{center}
\begin{tikzpicture}
  \node (c0) at (2,0) {$c_0$};
  \node (c1) at (3,1) {$c_1$};
  \node (c2) at (3,2) {$c_2$};
  \node (c3) at (2,3) {$c_3$};
  \node (c4) at (2,4) {$c_4$};
  \node (r0) at (0,1) {$3$};
  \node (r1) at (1,1) {$3$};
  \node (r2) at (1,2) {$3$};
  \node (r3) at (3,0) {$1$};
  \draw (c0) -- (c1) -- (c2) -- (c3) -- (c4);
  \draw (r0) -- (r2);
  \draw (r1) -- (r2);
  \draw (r2) -- (c3);
  \draw (r3) -- (c1);
\end{tikzpicture}

\end{center}
\caption{A tree $P$ with maximal chain $C = \{c_0 < c_1 < c_2 < c_3 < c_4\}$ and
elevations of $P - C$ labeled.}  \label{fig:elevation}
\end{figure}

\begin{Definition}\label{def:elevation}
Let $C$ be a fixed maximal chain in a forest $P$ with $|C| = r$.
For each $u \in P - C$, define the \textit{elevation} of $u$ to be \[
e_u \coloneqq \#\{v \in C : u \not\leq v\}. \] See \Cref{fig:elevation}. Let $s_k(P,C) $ be the
number of elements in $P - C$ with elevation at least $k-n+r$,  
$$s_k(P,C) \coloneqq \#\{u \in P - C : e_u \geq k-n+r\}.$$
\end{Definition}

For example, if $u$ is attached to the root of the tree which is the
maximal element of $C$, then the elevation is $e_u = r-1$. If $u$ is
attached to the second-smallest element of $C$, then $e_u = 1$.  We
see that $e_u = r$ if and only if $u$ is not connected by a path to
$C$.  Thus, if $P$ is a tree, then $1\leq e_{u}< r$, so 
$s_{n-r}(P,C)=n-r$, and $s_{n}(P,C)=0$.

If $P$ is a tree and $C$ is a chain in $P$, then $P - C$ is a forest
so both have associated cumulants. We may relate $\kappa_d^P$ and
$\kappa_d^{P - C}$ using the numbers $s_{k}(P,C)$ as follows.

\begin{Lemma}\label{lem:forest_dc}
  Let $C$ be a maximal chain in a tree $P$ with $n$ elements and
  $|C|= r$.  Then for each $d \in \bZ_{\geq 1}$,
  \begin{equation}\label{eq:forest_dc}
    \kappa_d^P = \kappa_d^{P - C}
      + \frac{B_{d}}{d} \sum_{u \in P - C}
         \sum_{k=n-r+1}^{n-r+e_u} \mathbf{h}_{d-1}(k, k-s_k)
  \end{equation}

  \begin{proof}
    Let $C = v_n > v_{n-1} > \cdots > v_{n-r+1}$.  Note that for
    $u \in P - C$, we have $u < v_k$ if and only if $e_u <
    k-n+r$. Consequently, for all $n-r<k\leq n$  we have 
    \begin{align*}
      h_{v_k}
        &= k - n + r + \#\{u \in P - C : u < v_k\} \\
        &= k - n + r + \#\{u \in P - C : e_u < k-n+r\} \\
        &= k - n + r + (n - r - \#\{u \in P - C : e_u \geq k-n+r\} )\\
        &= k - s_{k}.
    \end{align*}
    Thus, 
    \begin{align*}
      \sum_{k=n-r+1}^n k^d - h_{v_k}^d
      &= \sum_{k=n-r+1}^n (k - h_{v_k}) \mathbf{h}_{d-1}(k, h_{v_{k}}) \\
      &= \sum_{k=n-r+1}^n s_k \mathbf{h}_{d-1}(k, k-s_k) \\
      &= \sum_{k=n-r+1}^n \#\{u \in P - C : e_u \geq k-n+r\}
           \cdot \mathbf{h}_{d-1}(k, k-s_k) \\
      &= \sum_{u \in P - C} \sum_{k=n-r+1}^{n-r+e_u} \mathbf{h}_{d-1}(k, k-s_k).
    \end{align*}
    Therefore, \eqref{eq:forest_dc} follows from the
    cumulant formula in \Cref{cor:forest_cumulants}.  
  \end{proof}
\end{Lemma}

If $P$ is a standardized tree with maximal chain $C$ of size
$|C|=r>1$, it has an element $u \in P-C$ with $e_u = r-1$, so $e_u/r
\sim 1$ for $r$ large.  As we saw in \Cref{sub:inf.sums},
renormalizing a multiset by the maximum value is a useful technique
while not changing the corresponding standardized general uniform
sum distribution.  Consequently, we consider the re-scaled multiset of
elevations $\mathbf{e}/r=\{e_u/r : u \in P - C\},$ which are then
related to the rescaled cumulants $\kappa_d^{P}/r^d$.

\begin{Lemma}\label{lem:cum.sum.to.zero}
  Suppose we have a sequence of standardized trees $P$ such that the
number of elements $n \to \infty$ and the rank $r$ satisfies $n-r =
o(n^{1/2})$, i.e.~$(n-r)/n^{1/2} \to 0$. Let $C$ be a maximal length
chain in $P$. Then, for each $d \in \bZ_{\geq
1}$, $$\frac{\kappa_d^{P-C}}{r^d} =O\left( \frac{(n-r)^{d+1}}{r^d} \right)
\to 0 \hspace{.2in} \text{and}  \hspace{.2in} \frac{\kappa_d^P}{r^{d}} \sim  \sum_{u \in P - C} \left(\frac{e_u}{r}\right)^d=|\mathbf{e}/r|_{d}^{d}.$$


\begin{proof}
Since $n-r = o(n^{1/2})$ and $n \to \infty$, we find $n \sim r$, and
so $n-r = o(r^{1/2})$. Consequently, $(n-r)^2/r \to 0$, and more
generally $(n-r)^{d+1}/r^d \to 0$ for all $d \geq 1$. Therefore, 
\[
 \frac{\kappa_d^{P-C}}{r^d}  = \frac{1}{r^d} \sum_{k=1}^{n-r} k^d -
\frac{1}{r^d} \sum_{u \in P - C} h_u^d =
O\left(\frac{(n-r)^{d+1}}{r^d}\right) \to 0.
\]

Consider the formula for $\kappa_{d}^{P}/r^{d}$ obtained from
\eqref{eq:forest_dc} by dividing both sides by $r^{d}$.  The first
term goes to 0 by the argument above.  The second term is bounded
above and below by 
\begin{equation}\label{eq:forest_sim.1} d \sum_{u \in P - C}
\sum_{k=n-r+1}^{n-r+e_u} (k-s_k)^{d-1} \leq \sum_{u \in P - C}
\sum_{k=n-r+1}^{n-r+e_u} \mathbf{h}_{d-1}(k, k-s_k) \leq d \sum_{u \in
P - C} \sum_{k=n-r+1}^{n-r+e_u} k^{d-1}.
\end{equation}
In \Cref{lem:upper.bound.sum} and \Cref{lem:lower.bound.sum} below, we
will show that, after dividing by $r^d$, both bounds in
\eqref{eq:forest_sim.1} are asymptotic to $\sum_{u \in P - C}
\left(\frac{e_u}{r}\right)^d$. Thus, $\kappa_{d}^{P}/r^{d} \sim \sum_{u
\in P - C} \left(\frac{e_u}{r}\right)^d$.
\end{proof}
\end{Lemma}

\begin{Lemma}\label{lem:upper.bound.sum}
With the same hypotheses as \Cref{lem:cum.sum.to.zero},
\[
\frac{d}{r^{d}}  \sum_{u \in P - C} \sum_{k=n-r+1}^{n-r+e_u} k^{d-1} 
\sim    \sum_{u \in P - C} \left(\frac{e_u}{r}\right)^d.
\]

\begin{proof}
    From \Cref{lem:riemann_estimate}, we have
    \begin{equation}\label{eq:forest_sim.2}
      \begin{split}
        &\sum_{u \in P - C} \left[\left(\frac{n-r}{r} + \frac{e_u}{r}\right)^d
           - \left(\frac{n-r}{r}\right)^d\right]
           \leq \frac{d}{r^d} \sum_{u \in P-C} \sum_{k=n-r+1}^{n-r+e_u} k^{d-1} \\
            &\leq \sum_{u \in P - C} \left[\left(\frac{n-r}{r} + \frac{e_u}{r}\right)^d
              - \left(\frac{n-r}{r}\right)^d\right]
              + \frac{d}{r} \sum_{u \in P - C} \left[\left(\frac{n-r}{r} + \frac{e_u}{r}\right)^{d-1}
              - \left(\frac{n-r}{r}\right)^{d-1}\right].
      \end{split}
    \end{equation}
    Consider the lower bound in \eqref{eq:forest_sim.2}. By
    \Cref{lem:cum.sum.to.zero}, 
    $\sum_{u \in P - C} \left(\frac{n-r}{r}\right)^d = \frac{(n-r)^{d+1}}{r^d} \to 0$
     for all $d \geq 1$. Furthermore, 
    \begin{align*}
      \sum_{u \in P - C} \left(\frac{n-r}{r} + \frac{e_u}{r}\right)^d
        &= \sum_{u \in P - C} \sum_{i=0}^d \binom{d}{i}
           \left(\frac{n-r}{r}\right)^i \left(\frac{e_u}{r}\right)^{d-i} \\
        &= \sum_{i=0}^d \binom{d}{i} \left(\frac{n-r}{r}\right)^i
            \sum_{u \in P - C} \left(\frac{e_u}{r}\right)^{d-i} \\
        &\leq \sum_{u \in P - C} \left[\left(\frac{e_u}{r}\right)^d
          + \sum_{i=1}^d \binom{d}{i} \left(\frac{n-r}{r}\right)^i \cdot 1^d\right] \\
        &\sim \sum_{u \in P - C} \left(\frac{e_u}{r}\right)^d.
    \end{align*}
    The first term in the upper bound in \eqref{eq:forest_sim.2} is dominant by a
    similar argument. Therefore, since the upper and lower bound in \eqref{eq:forest_sim.2}
asymptotically converge to the same sum, 
it follows that
      \[ \frac{d}{r^d} \sum_{u \in P - C} \sum_{k=n-r+1}^{n-r+e_u} k^{d-1}
          \sim \sum_{u \in P - C} \left(\frac{e_u}{r}\right)^d. \]
\end{proof}
\end{Lemma}

\begin{Lemma}\label{lem:lower.bound.sum}
With the same hypotheses as \Cref{lem:cum.sum.to.zero},
\[
\frac{d}{r^{d}} \sum_{u \in P - C} \sum_{k=n-r+1}^{n-r+e_u} (k-s_k)^{d-1}
\sim  \sum_{u \in P - C} \left(\frac{e_u}{r}\right)^d.
\]

\begin{proof}
    Consider the expansion
    \begin{align*}
      \frac{d}{r^d} \sum_{u \in P - C} \sum_{k=n-r+1}^{n-r+e_u} (k-s_k)^{d-1}
        = &\sum_{u \in P - C} \frac{d}{r^d} \sum_{k=n-r+1}^{n-r+e_u} k^{d-1} \\
        &+ \sum_{i=1}^{d-1} (-1)^i \binom{d-1}{i} \sum_{u \in P - C} \frac{d}{r^d}
             \sum_{k=n-r+1}^{n-r+e_u} k^{d-1-i} s_k^i.
    \end{align*}
    Since $s_k$ by definition counts a subset of $P - C$, we have $s_k \leq n-r$. Thus,
    for each $1 \leq i \leq d-1$, we have
    \begin{align*}
      \sum_{u \in P - C} \frac{d}{r^d} \sum_{k=n-r+1}^{n-r+e_u} k^{d-1-i} s_k^{i}
        &\leq  \sum_{u \in P - C} \frac{d(n-r)^i}{r^d} \cdot \sum_{k=n-r+1}^n k^{d-i-1} \\
        &= \frac{d(n-r)^{i+1}}{r^d} \cdot \sum_{k=n-r+1}^n k^{d-i-1} \\
        &= O\left(\frac{(n-r)^{i+1}}{r^d} \cdot r^{d-i}\right)
          = O\left(\frac{(n-r)^{i+1}}{r^i}\right).
    \end{align*}
    By \Cref{lem:cum.sum.to.zero}, $\frac{(n-r)^{i+1}}{r^i} \to 0$,
    and so by \Cref{lem:upper.bound.sum}, it follows that
      \[ \frac{d}{r^d} \sum_{u \in P - C} \sum_{k=n-r+1}^{n-r+e_u} (k-s_k)^{d-1}
          \sim \frac{d}{r^d} \sum_{u \in P - C} \sum_{k=n-r+1}^{n-r+e_u} k^{d-1}
          \sim \sum_{u \in P - C} \left(\frac{e_u}{r}\right)^d. \]
\end{proof}
\end{Lemma}

We may combine the preceding results to prove the following more
explicit form of \Cref{thm:forest_infty} from the introduction.

\begin{Theorem}\label{thm:forest_degenerate_limits}
  Let $P$ denote an infinite sequence of standardized trees with $n$
elements and maximal chains $C$ of rank $r$ such that $n \to \infty$
and $n-r = o(n^{1/2})$.  Let $\mathbf{e}=\{e_u : u \in P - C\}$ be
the multiset of elevations for $P$ and $C$.  Then for each fixed $d
\in \bZ_{\geq 2}$ even, the cumulants of $\cX_P^{*}$ are
approximately
\begin{equation}\label{eq:forest_degenerate_limits.0}
(\kappa_d^P)^* \sim \frac{B_d/d}{(B_2/2)^{d/2}}
\left(\frac{|\mathbf{e}/r|_d}{|\mathbf{e}/r|_2}\right)^d =
\frac{B_d}{d} |\widehat{\mathbf{e}}|_{d}^{d}.
\end{equation}
The sequence of random variables $\cX_P^*$ converges in distribution
if and only if the multisets $\widehat{\mathbf{e}}$ converge pointwise
to some multiset $\mathbf{t} \in \bfP_{\DUSTPAN}$, in which case the
limiting distribution is $\cS_{\mathbf{t}} + \cN(0, \sigma) \in
\bfM_{\DUSTPAN}$ where $\sigma \coloneqq \sqrt{1 -
|\mathbf{t}|_2^2/12}$.  In particular, the sequence of random
variables $\cX_P$ are asymptotically normal if and only
if \begin{equation}\label{eq:forest_degenerate_an.0}
|\mathbf{e}/r|_2^{2} \coloneqq \sum_{u \in P - C}
\left(\frac{e_u}{r}\right)^2 \to \infty.
\end{equation}


  
  \begin{proof} Fix $d\geq 2$ even.  By hypothesis, $n-r =
o(n^{1/2})$, so \Cref{lem:cum.sum.to.zero} shows that
\begin{equation}\label{eq:P-C.approx} \frac{\kappa_d^P}{r^d} \sim
\frac{B_d}{d} |\mathbf{e}/r|_d^{d}.
\end{equation}
Therefore,  by \eqref{eq:rescaled.cumulants} 
\[
(\kappa_d^P)^* \sim \frac{B_d/d}{(B_2/2)^{d/2}}
\left(\frac{|\mathbf{e}/r|_d}{|\mathbf{e}/r|_2}\right)^d.
\]
Since $\mathbf{e}/r$ is finite, $|\mathbf{e}/r|_2$ exists, so the
hat-operation is defined on $\mathbf{e}/r$ and $\widehat{\mathbf{e}/r}
= \widehat{\mathbf{e}}$ after cancellation.  Hence,
\eqref{eq:forest_degenerate_limits.0} follows from the definition of
the hat-operation in \eqref{eq:hat.operation}.

By the Method of Moments/Cumulants (\Cref{thm:moments}) together with
\Cref{cor:frechet_shohat_converse}, the sequence $\cX_P^*$ converges
in distribution to some $\cX$ if and only if $\frac{B_d}{d}
|\widehat{\mathbf{e}}|_{d}^{d}$ converges to $\kappa_d^{\cX}$ for each
$d \in \bZ_{\geq 1}$.  By \Cref{cor:pointwise.convergence.equiv} and
the fact that $|\widehat{\mathbf{e}}|_{2}^{2}=12$ by definition, this
occurs if and only if $\widehat{\mathbf{e}}$ converges pointwise to
some $\mathbf{t}$.  Therefore, by \Cref{thm:PDelta_MDelta}, we have
$\mathbf{t} \in \bfP_{\DUSTPAN}$ and $\cX$ has the associated DUSTPAN
distribution $\Phi(\mathbf{t}) =\cS_{\mathbf{t}} + \cN(0, \sigma).$

In particular, the limiting distribution of $\cX_P^*$ is $\cN(0, 1)$
if and only if $\widehat{\mathbf{e}} \to \mathbf{0}$.  Now
$\widehat{\mathbf{e}} \to \mathbf{0}$ if and only if
$|\mathbf{e}/r|_\infty/|\mathbf{e}/r|_2 = 1/|\mathbf{e}/r|_2 \to 0$
since standardized trees have an element of elevation $r-1$. In
particular, the limit is $\cN(0, 1)$ if and only if $|\mathbf{e}/r|_2
\to \infty$.  \end{proof}
\end{Theorem}

\begin{Remark}
We note that considering only standardized trees in
\Cref{thm:forest_degenerate_limits} is essential for the ``if and only
if'' conditions to hold.  For example, consider a sequence of trees
$H_{n, r}$ with maximal chain $C$ of size $r$ such that $n \to \infty$
and $n-r = o(n^{1/2})$.  Since $\cL_{H_{n, r}}(q) = [n-r+1]_q!$
and $n-r \to \infty$, $\cX_{H_{n, r}}$ is asymptotically normal by
\cite{MR0013252}. However, we have elevation $e_u=1$ for all $u \in
H_{n, r} -C$. Therefore, $\sum_{u \in H_{n, r} - C} (e_u/r)^2 =
(n-r)/r^2 \to 0$ rather than $\infty$.
\end{Remark}

\begin{Remark}\label{rem:easy.construction}
  One can construct sequences of standardized trees with $n-r =
o(n^{1/2})$ where $\mathbf{e}/r$ converges to any prescribed finite
multiset $\mathbf{t} = (t_{1}\geq t_{2}\geq \cdots \geq t_{m}) \in
\widetilde{\ell}_{2}$ with $|\mathbf{t}|_\infty = 1$. For each
$N=m+3,m+4,\ldots$, let $r_{N}=N-m$.  To construct the tree $P_{N}$,
start with a chain $C_{N}=(v_{0}<v_{1}<\cdots < v_{r_{N}-1})$, and for
each nonzero value $1=t_{1}\geq t_{2}\geq \cdots \geq t_{m}$, add a
child to $v_{\lceil (r_{n}-1) t_{i} \rceil}$.  Finally, for each
$t_{i}=0$, add one additional child to $v_{1}$.  As constructed
$n=|P_{N}|=N$, $r=r_{N}-1$ and $n-r=m$ is constant.  Since $t_{1}=1$
by assumption, the root of $P_{N}$ has at least one child so it is a
standard tree. Furthermore, $m=|P_{N}-C_{N}|$ so the elevation
multiset of $P_{N}$ has exactly $m$ elements.  By construction, the
multisets $\mathbf{e}/r=\{e_{N}/r_{N}: u \in P_{N} -C_{N} \}$
approaches $\mathbf{t}$ as $N \to \infty.$ Therefore,
$\cS_{\widehat{\mathbf{t}}}$ is the limiting distribution of
$\cX_{P_{N}}^*$.  By \Cref{cor:Minfty_closure}, we know that the
closure of $\{\cS_{\widehat{\mathbf{t}}}: \mathbf{t} \in
\widetilde{\ell}_2, \mathbf{t}\text{ is finite}\}$ is
$\bfM_{\DUSTPAN}$.  Thus, $\bfM_{\Forest} \cup \bfM_{\DUSTPAN} \subset
\overline{\bfM_{\Forest}}$ as claimed in \Cref{sec:intro}.
\end{Remark}

\begin{repCorollary}{cor:forest_infty}
Let ${\epsilon\Tree}$ be the set of standardized trees $P$ for which $|P|
- \rank(P) \leq  |P|^{\frac{1}{2}}$.  Let
$\bfM_{\epsilon\Tree} \coloneqq \{\cX_P^* : P \in {\epsilon\Tree} \}
\subset \bfM_{\Forest} $ be the corresponding metric space of
distributions. Then
\begin{equation} \overline{\bfM_{\epsilon\Tree}}
= \bfM_{\epsilon\Tree} \sqcup \bfM_{\DUSTPAN},
\end{equation}
which is
(sequentially) compact. Moreover, the set of limit points of
$\bfM_{\epsilon\Tree}$ is $\bfM_{\DUSTPAN}$.

  \begin{proof} By the construction in \Cref{rem:easy.construction},
we know $\overline{\bfM_{\epsilon\Tree}} \supset \bfM_{\DUSTPAN}$, and
$\bfM_{\DUSTPAN}$ is closed by \Cref{cor:cM_compact}.  Furthermore,
we have $(|P| - \rank(P))/|P|^{1/2} < |P|^{-\epsilon} \to 0$, so
\Cref{thm:forest_degenerate_limits} applies. Thus,
for every sequence of trees $P \in \epsilon\Tree$ with $|P|\to \infty $
such that the corresponding random variables $\cX_P^* \in
\bfM_{\epsilon\Tree}$ converge in distribution, we know the distribution
must be a DUSTPAN distribution.
On the other hand, for every sequence of trees $P \in \epsilon\Tree$ such
that the corresponding random variables $\cX_P^* \in
\bfM_{\epsilon\Tree}$ converge in distribution but $|P|$ is bounded, we
must have a subsequence where $n=|P|$ is eventually constant.  There
are only a finite number of standardized trees of size $n$ in
$\epsilon\Tree$, so we can further restrict to a sequence where each $P$
is a particular tree, in which case the limiting of $\cX_P^*$ is
itself $\cX_P^* \in \bfM_{\epsilon\Tree}$.
\end{proof}
\end{repCorollary}

\section{Future work}\label{sec:future}

In addition to the open problems mentioned in \Cref{sec:intro} and
\Cref{question:ktoinf}, we pose the following questions for future study.

\begin{Question}
  Suppose we have a sequence of standardized trees such that $n \to
\infty$ where $n-r$ grows at least as fast as $n^{1/2}$ but no faster
than $n^{11/12}$ in the sense that $n-r \not\in o(n^{1/2})$ and $n-r
\not\in \omega(n^{11/12})$.  When is the corresponding sequence of
distributions asymptotically normal?  What non-normal limit laws are
possible?
\end{Question}

\begin{Question}
  Does $\weft(\lambda) \to \infty$ if \textit{and only if} $\cX_{\lambda; m}[\rank]$ is
  asymptotically normal? See \eqref{eq:aft-like-stat-ssyt}.
\end{Question}

\begin{Question}
  Consider the set of rooted, unlabeled forests with $n$ vertices, sampled uniformly
  at random. What is the expected value of the rank $r$, i.e.~the maximum length of a path
  starting at a root of a tree in the forest? How does $r$ compare to $n$ asymptotically
  as $n \to \infty$?
\end{Question}

See \cite{Pittel.94} for growth rates of the form $r \approx \log n$ for certain
random tree generation techniques. For the number of rooted, unlabeled forests
with $n$ vertices, $t$ trees, and rank $r$, see \cite[A291336]{oeis}.
Broutin--Flajolet \cite[Thm.~3]{MR2956055} showed that
$\mathbb{E}[r] \sim C\sqrt{n}$ for an explicit constant $C>0$ when
considering rooted, unlabeled \textit{binary} trees.
The corresponding problem when order is imposed either by
labeling the vertices (resulting in labeled trees) or by ordering the
children (resulting in planar trees) is older, though the
$\mathbb{E}[r] \sim D\sqrt{n}$ behavior is common throughout;
see \cite[p.1]{MR2956055} for a summary and further references.

\begin{figure}[ht]
 \includegraphics[height=5cm]{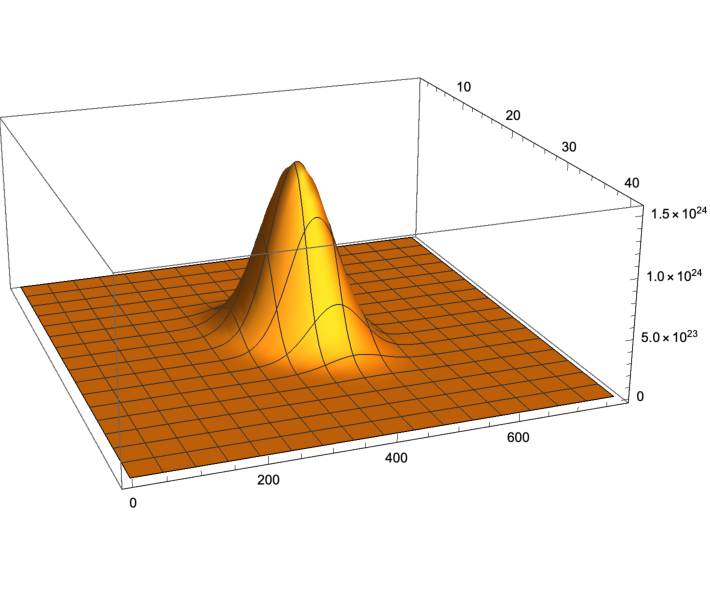} \\
 \includegraphics[height=2cm]{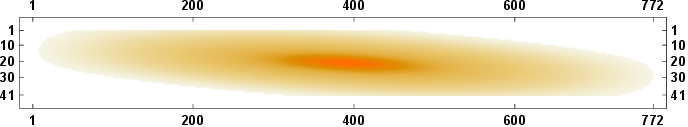}
    \caption{Plots of coefficients of the $q,t$-hook length formula \eqref{eq:qt_SYST}
    with $\lambda = (25, 4, 3, 3, 1, 1, 1, 1, 1)$.}
    \label{fig:qt_example}
\end{figure}

In \cite{Swa20}, the following $q, t$-analogue of the hook length formula
\eqref{eq:q_SYT} is given. Let $(r, c) \in \lambda$ denote a cell in row $r$
and column $c$. Then
\begin{equation}\label{eq:qt_SYST}
  [n]_q! \prod_{(r, c) \in \lambda} \frac{q^{r-1} + tq^{c-1}}{[h(r, c)]_q}
\end{equation}
is the generating function for a pair of statistics $(\maj, \operatorname{neg})$ on
\textit{standard supertableaux} of shape $\lambda$. The $t=0$ case of
\eqref{eq:qt_SYST} yields \eqref{eq:q_SYT}. While \eqref{eq:qt_SYST} is
not literally a quotient of $q$-integers, it is evidently ``nearly'' such a quotient.
Computational evidence suggests the distributions are ``typically'' bivariate normal
with non-trivial covariance, which is strikingly similar to the distributions
encountered by Kim--Lee \cite{1811.04578} for $(\des, \maj)$ on permutations
in fixed conjugacy classes. See \Cref{fig:qt_example} for sample data.

\begin{Question}
  What are the possible limiting distributions of the coefficients of the
  $q,t$-hook length formula \eqref{eq:qt_SYST}? What is the support of
  \eqref{eq:qt_SYST}?
\end{Question}

One referee asked the following natural question, saying ``A result of
this form could give a \textit{conceptual} explanation for some of the
results.''   The authors regard this as an important question, but we do not
expect a simple answer.

\begin{Question}
Can one give a
formula for the statistic rank on $\SSYT_{\lambda,m}$ as a natural sum
of $m$ natural statistics on the tableaux, and then to show that they
are (asymptotically) independent, converging to the uniform law?
\end{Question}

\section*{Acknowledgments}\label{sec:ack}

We would like to thank Persi Diaconis, Matja\v{z} Konvalinka,
Svante Janson, Soumik Pal, Richard Stanley, and John Stembridge for
helpful discussions related to this work.

\newpage

\bibliography{refs}{}
\bibliographystyle{alpha}

\end{document}